\documentclass[11pt,a4paper]{article}

\usepackage[utf8]{inputenc}
\usepackage[T1]{fontenc}
\usepackage[margin=2.5cm]{geometry}
\usepackage{graphicx} 
\usepackage{xcolor}
\usepackage{amsmath, amssymb, bm}
\usepackage{amsthm}
\usepackage{subcaption}
\usepackage{adjustbox}
\usepackage{booktabs}
\usepackage{pgf}
\usepackage{siunitx}
\usepackage[hidelinks]{hyperref}

\numberwithin{equation}{section}

\theoremstyle{remark}
\newtheorem{remark}{Remark}

\newcommand{\joel}[1]{{\color{black}#1}}

\colorlet{joseantocolor}{black}
\newcommand{\joseanto}[1]{{\color{joseantocolor}#1}}

\newcommand{\xv}{\mathbf{x}}

\newcommand{\qv}{\mathbf{q}}
\newcommand{\nv}{\mathbf{n}}

\newcommand{\iv}{\mathbf{i}}
\newcommand{\jv}{\mathbf{j}}
\newcommand{\ev}{\mathbf{e}}
\newcommand{\M}{\mathbf{M}}

\newcommand{\R}{\mathbb{R}}
\newcommand{\N}{\mathbb{N}}

\title{High-order IMEX local discontinuous Galerkin  methods with sparse grids for high-dimensional nonlinear parabolic PDEs}
\author{Joel P. Villarino\thanks{CITIC and Department of Mathematics,
    Universidade da Coru{\~n}a, Campus de Elvi{\~n}a s/n,
    15071 A Coru{\~n}a, Spain}
    \qquad Manuel~J.~Castro\thanks{Departamento de An{\'a}lisis Matem{\'a}tico,
    Estad{\'i}stica e Investigaci{\'o}n Operativa y Matem{\'a}tica Aplicada,
    Universidad de M{\'a}laga, Campus de Teatinos s/n, 29080 M{\'a}laga, Spain}\\
    J.A.~Garc{\'i}a--Rodr{\'i}guez\footnotemark[1]
    \qquad J.G.~L{\'o}pez--Salas\footnotemark[1]}
\date{}

\begin{document}

\maketitle

\begin{abstract}
In this work, we develop a high-order IMEX-LDG scheme on sparse grids for solving nonlinear parabolic PDEs in moderate to high dimensions. The method combines the advantages of LDG methods for convection--diffusion problems, IMEX-RK time integrators for stiff systems, and sparse-grid techniques to mitigate the curse of dimensionality. The resulting scheme is fully parallelized in MPI, is robust, efficient, and scalable, making it suitable for a wide range of applications involving high-dimensional nonlinear parabolic PDEs.
\par\medskip
\noindent\textbf{Keywords:} Local discontinuous Galerkin method; IMEX Runge--Kutta methods; sparse grids; high-dimensional parabolic PDEs.
\par\smallskip
\noindent\textbf{Mathematics Subject Classification:} 65M60; 65M12; 65M70; 65N30.
\end{abstract}

\section{Introduction}

\newcommand{\LDGIMEXCites}{\cite{Wang2015,Wang-Shu-Zhang-2016,Haijin-wang-Shu-M2AN-2016,Haijin-Shu-SIAM-JNA-2016,WANG2018164,Wang2018}}
\newcommand{\ChengSparseRefs}{\cite{Cheng2016,Wang2016,Tao2019,Cheng2019}}
\newcommand{\LDGAppsRefs}{\cite{yan2002localKdV,yan2002localHigher,dong2009analysis}}
\newcommand{\IMEXCoreRefs}{\cite{Ascher1995,Ascher1997,Calvo2001,Wang2015,Kennedy2003}}
\newcommand{\IMEXFVCites}{\cite{Russo05,Boscarino2011,Boscarino2009Uniform,Boscarino2013Diffusion,Boscarino2017FV,BoscarinoRusso2013Flux}}
\newcommand{\SparseFEMRefs}{\cite{Garcke2013,Bungartz_Griebel_2004,Bungartz1995,Griebel2005}}
\newcommand{\SemiLagrangianRefs}{\joseanto{\cite{douglas1982convection,pironneau1982transport,bercovier1983finite,benitez2014pure,bonaventura2014semilagrangian,hecht2024dual}}}

\joseanto{
The main goal of this work is to develop efficient, high-order numerical methods for solving $d$-dimensional scalar nonlinear convection-diffusion problems of the form
    \begin{equation}
        \label{eq:ldg:pde}
        \partial_t u + \nabla_x \cdot \mathbf{f}(u) = \nabla_x \cdot \left(A\,\nabla_x u\right)+s(u),\quad (\xv, t)\in Q_T = \Omega \times (0,T),
    \end{equation}
    with initial state $u(\xv,0)=u_0(\xv)$.  Here $\Omega$ is a bounded rectangular domain in $\R^d$, $\mathbf{f}=\left(f_1,\dots,f_d\right)^T$ is the flux function, $A = (a_{ij})\in\mathcal{M}_d(\R)$ is a symmetric, semi-positive definite (SSP) matrix which could depend on the spatial variable, and $s$ is a function that defines the source term. Without loss of generality, and for greater clarity in the presentation of the method, we assume that the domain is square and the boundary conditions are periodic.

    The focus on high-order accuracy is motivated by the fact that high-order methods are more efficient, i.e., for a prescribed error tolerance, such methods can achieve the target accuracy with fewer degrees of freedom, thereby reducing the overall computational cost.
In line with this objective, in this work we concentrate on implicit--explicit (IMEX) Runge--Kutta (RK) local discontinuous Galerkin (LDG) methods on sparse grids, a combination that provides accuracy, robustness, and favorable computational performance across a broad range of convection--diffusion regimes.

}

\joseanto{
As is well known, designing robust and efficient Galerkin-type numerical methods for nonlinear parabolic PDEs in moderate to high dimensions is a challenging task and remains an active area of research in numerical analysis. 
Several difficulties arise, most notably, the treatment of convection-dominated regimes (requiring the use of numerical fluxes or artificial diffusion techniques), the presence of nonlinearities that may require the solution of nonlinear systems, the temporal stiffness in the time-stepping scheme for the semi-discrete ODE system under mesh refinement, the curse of dimensionality, the complexity of the resulting linear systems and the huge computational cost of the designed numerical methods. Each section of this work addresses these issues and explains the strategies adopted here to overcome them.

To begin with, the first three core difficulties associated with nonlinear parabolic PDEs already arise in any spatial dimension, even in dimension one: the need for robust upwinding of convective terms, the treatment of nonlinear fluxes and source terms, and the stiffness of the ODE system obtained after spatial semi-discretization.

For the convective and nonlinear contributions, the discretization must remain robust in convection-dominated regimes while preserving stability in the presence of strong nonlinear fluxes and source terms. This motivates the use of upwind numerical schemes. 
Finite volume methods are an excellent choice for incorporating upwinding and can also handle strong nonlinearities, but their extension to high orders requires reconstruction techniques that can be computationally expensive and cumbersome in high dimensions. 
Another alternative for handling convective terms are semi-Lagrangian finite element methods \SemiLagrangianRefs, which can attain moderate order, although continuous Galerkin formulations induce global couplings, leading to linear systems with high computational cost and limited suitability for massive parallelization. 
In addition, their extension beyond two or three dimensions is complex and costly due to the searches and interpolations along characteristic curves.
Another drawback is the treatment of nonlinearities in the source terms or in the convection, which may require solving nonlinear systems, resulting in higher computational cost and potential lack of robustness.
}

Discontinuous Galerkin methods, particularly the LDG, are designed for convection--diffusion problems and incorporate upwinding by construction. 
The LDG method was first introduced by Cockburn and Shu \cite{Cockburn1998}. Although it was originally proposed for convection--diffusion problems, it has been applied to dispersive KdV-type equations, bi-harmonic equations, and fifth-order dispersion equations \LDGAppsRefs. 
In LDG, the local polynomial representation on each element yields sparse (block-diagonal) operators, provides spectral-like accuracy for smooth solutions, and allows static condensation of auxiliary gradient variables, see \cite{Cockburn1998,Shu2016,Wang2015,Haijin-Shu-SIAM-JNA-2016,CASTILLO20132279}. 
The method also handles nonlinearities in the source and convection terms naturally, as they are treated explicitly. For these reasons, LDG methods are a very attractive option for solving nonlinear parabolic PDEs in high dimensions, when the non-linearities are located in the source and convection terms.

\joseanto{A further difficulty arises after spatial discretization, since the resulting semi-discrete ODE systems are typically stiff ($\Delta t \approx \min_k h_k^2$), especially under mesh refinement; as a consequence, explicit RK time discretizations for convection-diffusion problems are efficient only when convection dominates and become rapidly impractical in the diffusive regime.
On the other hand,} fully implicit schemes require solving nonlinear systems at each time step whenever convective fluxes or nonlinear source terms are present. 
To address these issues, for the semi-discrete LDG system, we adopt IMEX-RK methods: the parabolic/diffusive terms are treated implicitly and the convective/source terms explicitly, which relaxes the diffusive time-step restriction (elliptic systems amenable to iterative solvers) and avoids treating convection or nonlinear source terms implicitly \IMEXCoreRefs.
IMEX-RK methods have been extensively applied in the context of stiff ODE systems resulting from the semidiscretization of PDEs: for example, coupled with finite differences and finite volumes in \IMEXFVCites \ and with Local Discontinuous Galerkin (LDG) methods in \LDGIMEXCites. 
With this strategy, a double objective is achieved: the diffusive time-step restriction is relaxed, and one can solve semilinear PDEs without having to solve nonlinear problems at each time step, which is crucial both in terms of computational cost and robustness.

An additional difficulty is the curse of dimensionality. 
In full tensor-product meshes (full), the number of degrees of freedom grows exponentially with the dimension, making classical discretizations infeasible. 
\emph{Sparse} grids based on Smolyak indices (\emph{hyperbolic--cross}) mitigate this growth by activating only carefully selected subsets of multi-indices, reducing the complexity by several orders of magnitude while maintaining high order. 
Many works related to sparse-grid finite element formulations can be found in \SparseFEMRefs. 
Some formulations combine these ideas with finite volume schemes, including sparse-grid WENO reconstructions of third and fifth order \cite{LuChenZhang2018WENO,Zhu2021}, while sparse-grid DG formulations have been developed for transport, elliptic, and Vlasov--Maxwell systems \ChengSparseRefs.

\joseanto{Finally,} attending to the computational cost, the algebraic systems arising from the implicit LDG stages are solved via Krylov subspace methods equipped with preconditioners tailored to the block-structured viscous operator. 
The matrix assembly, the explicit convection evaluation, and the  preconditioner and solver are fully parallelized using message passing interface (MPI), ensuring scalability to large clusters and allowing for dealing with large problem sizes.

This work is organized as follows. 
\joseanto{
Section \ref{sec:ldg:num} develops the semi-discrete LDG formulation on tensor grids. Section \ref{ldg:imex} details the IMEX-RK time integrators and the fully discrete numerical scheme. 
Section \ref{ldg:combination} details the sparse-grid combination technique, covering the construction of semi-coarsened meshes, the formulation of local LDG operators, and the coupling strategy among component grids. 
Finally, in Section \ref{ldg:experiments} we report numerical experiments that illustrate accuracy, efficiency, and scalability for general PDEs.}

\section{Semi-discrete LDG discretization} \label{sec:ldg:num}
     Since the DG formulations are grounded in the finite element framework, we \joseanto{start by introducing} the definition of the one-dimensional finite element space utilized in this work.
    Building upon this foundation, the extension to the multidimensional case is straightforward \joseanto{thanks} to the mesh topology employed, essentially structured grids composed of rectangular (tensor-product) cells.

    Consider a general interval $[a, b],\,a<b$. 
    We define a partition composed of $l$ elements denoted by $I_e=[x_{e-1/2}, x_{e+1/2}]$. 
    Hereafter, we assume that the grid points are uniformly spaced, yielding a constant mesh size $h=x_{e+1/2}-x_{e-1/2}\,\,\forall e = 1,\dots,l$. 
    We denote this mesh as $\Omega_l([a,b])$, and define the associated discontinuous finite element space
    \begin{equation}
        \label{eq:ldg:fem-space:1d}
        V_{l}^{\nu} = \{v\in L^2\left((a,b)\right)\,:\,v|_{I_e}\in\mathcal{P}_{\nu}(I_e)\,,\,\forall\,e=1,\dots,l\},
    \end{equation}
    where $\mathcal{P}_{\nu}(I_e)$ denotes the space of polynomials of degree at most $\nu$ defined on $I_e$. 
    It is worth noting that the term discontinuous stems from the fact that continuity is not enforced across element boundaries, that is,
    \begin{equation}
        \label{eq:ldg:fem-space:discontinuity-property}
        v|_{I_e}(x_{e+1/2})\neq v|_{I_{e+1}}(x_{e+1/2}),\quad v\in V_{l}^{\nu},
    \end{equation}
    thereby explicitly allowing discontinuous function representations.
    
    For a $d-$dimensional domain $\Omega=[a,b]^d$, we introduce the multi-index vector $\bm{\alpha} = (\alpha_1, \alpha_2,\dots,\alpha_d)\in \N^d$ to indicate the mesh size in each coordinate direction.
    Regarding operations involving multi-indices, unless stated otherwise, they are to be interpreted as component-wise operations.
    For example, for two multi-indices $\bm{\alpha},\bm{\beta}\in\N^d$, the relation $\bm{\beta}\leq\bm{\alpha}$ is defined as $\bm{\beta}\leq\bm{\alpha} \iff \beta_k\leq\alpha_k,\,\forall\,k=1,\dots,d$.
    Another frequently encountered operation involves the comparison between a natural number and a multi-index.
    In this context, the comparison is also performed component-wise. 
    For instance, given a multi-index $\bm{\beta}\in\N^d$ and a natural number $\hat{b}\in\N$, the relation $\bm{\beta}\leq \hat{b}$ is defined as $\bm{\beta}\leq \hat{b} \iff \beta_k\leq \hat{b},\,\forall\,k=1,\dots,d$.
    
    On the $d$-dimensional domain $\Omega$, we define the tensor-product grid as the tensor product of the underlying one-dimensional meshes,  
    \begin{equation} \label{eq:ldg:grid:nd}
        \Omega_{\bm{\alpha}}([a,b]^d)=\bigotimes_{k=1}^d\Omega_{\alpha_k}([a,b]).
    \end{equation}
    Under this construction we define the mesh size vector as $\mathbf{h}_{\bm{\alpha}} = \left[h_{\alpha_1},h_{\alpha_2},\dots,h_{\alpha_d}\right]^T$,
    and the total number of cells of the grid is given by $N_{\bm{\alpha}}=\prod_{k=1}^d\alpha_k$. 

    Let $\bm{\beta}$ be a multi-index such that $\bm{\beta}\leq\bm{\alpha}$.
    We denote by $Q_{\bm{\beta}}$ the elementary cell associated with the index $\bm{\beta}$, defined by
    \begin{equation}
        Q_{\bm{\beta}} = \{\xv\,\,:\,\,x_i\in I_{\beta_i},\,\,i=1,\dots,d\,\}.
    \end{equation}
    By construction, the full tensor-product grid is the union of all such elementary cells, that is,
    \begin{equation}
        \Omega_{\bm{\alpha}}([a,b]) = \bigcup_{\bm{\beta}\leq\bm{\alpha}}Q_{\bm{\beta}}.
    \end{equation}
    The boundary of this elementary cell, denoted by $\partial Q_{\bm{\beta}}$, decomposes into \joseanto{$2\cdot d$} faces of dimension $(d-1)$. 
    These faces are obtained by fixing one coordinate at an endpoint of the interval $I_{\beta_k}$ while allowing the others to vary.
    Specifically, for each coordinate direction $k=1,\dots,d$, we define the left and right faces as
    \begin{subequations}
    \label{eq:ldg:cell:boundaries}
    \begin{align}
        \Gamma_{k,-}^{\bm{\beta}} &= \{\bm{x}\in\partial Q_{\bm{\beta}} : x_k=x_{\beta_k-1/2},\,x_l\in I_{\beta_l},\,\forall\,l\neq k\}, \quad \nv_{k}^-=-\ev_{k}, \label{eq:ldg:cell:boundaries:left} \\
        \Gamma_{k,+}^{\bm{\beta}} &= \{\bm{x}\in\partial Q_{\bm{\beta}} : x_k=x_{\beta_k+1/2},\,x_l\in I_{\beta_l},\,\forall\,l\neq k\}, \quad \nv_{k}^+=\ev_{k}. \label{eq:ldg:cell:boundaries:right} 
    \end{align} 
    \end{subequations}
    These faces are disjoint up to sets of measure zero (edges) and satisfy $\partial Q_{\bm{\beta}}=\bigsqcup_{k=1}^d(\Gamma_{k,-}^{\bm{\beta}}\cup\Gamma_{k,+}^{\bm{\beta}})$.
    
    Based on this mesh structure, the tensor-product finite element space is defined as
    \begin{equation}
        \label{eq:ldg:fem-space:multi}
        \bm{V}_{\bm{\alpha}}^{\nu}:=\{v\in L^2\left((a,b)^d\right)\,:\,v|_{Q_{\bm{\beta}}}\in\mathcal{P}_{\nu}(Q_{\bm{\beta}})\} = V_{\alpha_1}^{\nu}\otimes\dots\otimes V_{\alpha_d}^{\nu},
    \end{equation}
    where $P_{\nu}(Q_{\bm{\beta}})$ denotes the space of polynomials of degree at most $\nu$ in each coordinate variable on the cell $Q_{\bm{\beta}}$.

    Having established the finite element spaces used in this work, we now introduce the weak formulation of \eqref{eq:ldg:pde}. 
    To this end, it is convenient to rewrite the convection-diffusion PDE as a degenerate, first-order differential system.
    This reformulation is necessary because LDG methods cannot be directly applied to equations containing higher-order derivatives, as the solution space \eqref{eq:ldg:fem-space:multi} lacks sufficient regularity (consisting of piecewise polynomials that are discontinuous at element interfaces), \cite{Cockburn1998,Shu2016}. 
    Consequently, we introduce the auxiliary vector variable,
    \begin{equation}
        \label{eq:ldg:definition:q}
        \bm{q} = \left[q_1,\dots,q_d\right]^T=\left[\partial_{x_1}u,\dots,\partial_{x_d}u\right]^T = \nabla_{\bm{x}} u,
    \end{equation}
    which allows us to reformulate the original equation \eqref{eq:ldg:pde} as the equivalent first-order differential system
    \begin{equation}
        \label{eq:ldg:system_pde}
        \partial_t u  = -\nabla_x \cdot\left(\bm{f}(u)-A\,\qv\right) + s(u),\quad \qv-\nabla_x u = 0,\quad (t,\xv)\in Q_T,
    \end{equation}
    with the same initial state and boundary conditions. 

    Due to the local nature of the LDG method, \joseanto{we first formulate} the variational system  restricted to an arbitrary element $Q_{\bm{\beta}}\in\Omega_{\bm{\alpha}}$.

    Let $u_h(0,\xv)$ be the $L^2$-projection of the initial condition $u_0(\xv)$ on the finite element space $\bm{V}_{\bm{\alpha}}^{\nu}$, $u_h(0,\xv)\in \bm{V}_{\bm{\alpha}}^{\nu}$. 
    We seek to find, for any $t>0$, the approximation pair $\left(u_h(t,\,\cdot\,), \qv_h(t,\,\cdot\,)\right)\in \bm{V}_{\bm{\alpha}}^{\nu}\times\left(\bm{V}_{\bm{\alpha}}^{\nu}\right)^d$, satisfying the following variational formulation
    \begin{subequations}
    \label{eq:ldg:weak_form}
    \begin{align}
        \begin{split} \label{eq:ldg:weak_form:a}
      \int_{Q_{\bm{\beta}}}\partial_t u_h v_h &= \int_{Q_{\bm{\beta}}} \bm{f}(u_h)\cdot\nabla_x v_h - \oint_{\partial Q_{\bm{\beta}}} \widehat{f}_{\bm{n}}(u_h^-,u_h^+)\, v_h \\
            &\quad - \int_{Q_{\bm{\beta}}} A\,\mathbf{q}_h \cdot \nabla_x v_h +\oint_{\partial Q_{\bm{\beta}}} A\,\widehat{\mathbf{q}}_h\cdot\mathbf{n}^+ v_h + \int_{Q_{\bm{\beta}}}s(u_h)v_h,
        \end{split} \\[1em]
        \int_{Q_{\bm{\beta}}} \mathbf{q}_h\cdot\bm{w} &= -\int_{Q_{\bm{\beta}}}u_h\,\nabla_{\xv}\cdot\bm{w} +\int_{\partial Q_{\bm{\beta}}} \widehat{u}_h \bm{w}_h\cdot\mathbf{n}^+, \label{eq:ldg:weak_form:b}
    \end{align}
    \end{subequations}
    for all $v_h \in \bm{V}_{\bm{\alpha}}^{\nu}, \bm{w}_h \in (\bm{V}_{\bm{\alpha}}^{\nu})^d$.
    In each element $Q_{\bm{\beta}}$ with $\bm{\beta}\leq\bm{\alpha}$, we have introduced the interface traces, which we denote by $u_h^-$ and $u_h^+$.
    To ensure consistency with the alternating logic and the matrix operators defined later, we adopt a global Cartesian directional convention for the trace values.

    For any interface perpendicular to the spatial direction $k$, we define the left trace, $u_h^-$, and the right trace $u_h^+$ based on the increasing direction of the coordinate axes,
    \begin{equation} \label{eq:ldg:traces:definition} 
        u_h^-(\bm{x}) := \lim_{\epsilon \to 0^+} u_h(\bm{x} - \epsilon \mathbf{e}_k), \qquad u_h^+(\bm{x}) := \lim_{\epsilon \to 0^+} u_h(\bm{x} + \epsilon \mathbf{e}_k). 
    \end{equation}
    For the assembly of the local operators, we must map these global definitions to the local perspective of an element $Q_{\bm{\beta}}$:
    \begin{itemize}
        \item On the right face ($\Gamma_{k,+}^{\bm{\beta}}$), the internal value ($u^{\text{int}}$) corresponds to the global left trace $u_h^-$, while the external neighbor value ($u^{\text{ext}}$) corresponds to the global right trace $u_h^+$.
        \item On the left face ($\Gamma_{k,-}^{\bm{\beta}}$), the internal value ($u^{\text{int}}$) corresponds to the global right trace $u_h^+$, while the external neighbor value ($u^{\text{ext}}$) corresponds to the global left trace $u_h^-$.
    \end{itemize}
    
    Moreover, the terms denoted with hats in Equation \eqref{eq:ldg:weak_form} represent numerical fluxes.
    These are designed to approximate the physical fluxes at element interfaces, ensuring information transmission between neighboring elements while guaranteeing stability and conservation. 
    Specifically, the convective numerical flux, denoted by $\widehat{f}_{\nv}(u_h^-,u_h^+)$, is typically a single-valued Lipschitz-continuous function depending on the solution values from both sides of the interface.
    It must satisfy two key properties, consistency (that is, $\widehat{f}_{\bm{n}}(u, u)=\bm{f}(u)\cdot\bm{n}$) and conservation (that is, $\hat{f}_{\nv^+}(\xi_1,\xi_2) + \hat{f}_{\nv^-}(\xi_2,\xi_1) = 0$).
    Among the standard choices available in the literature \cite{Shu2016}, in this work, we employ the local Lax-Friedrichs (or Rusanov) flux \cite{Rusanov1962}, defined as
    \begin{equation} \label{eq:ldg:num-flux:rusanov}
        \Hat{f}_{\bm{n}}(u_h^-, u_h^+) = \dfrac{1}{2}\left(\bm{f}(u_h^-)+\bm{f}(u_h^+)\right)\cdot\bm{n} - \dfrac{\alpha_R}{2}\left(u_h^+ - u_h^-\right), 
    \end{equation}
    where the viscosity parameter is given by $\alpha_R = \max_{u \in \{u_h^-, u_h^+\}}|\bm{f}'(u)\cdot\bm{n}|$.
    Here, $\alpha_R$ represents an estimate of the maximum local characteristic speed (the largest eigenvalue of the Jacobian projected onto the normal), resulting in the simplest Godunov-type upwinding flux \cite{Toro2024}.
    
    Consistent with our Cartesian directional convention, we specialize the local Lax-Friedrichs flux for each coordinate direction $k$. 
    This formulation removes the generic dependency on the normal vector, as the flux $\widehat{f}_k(u_h^-,u_h^+)$ now specifically approximates the physical flux component $f_k(u)$ flowing in the positive $k$-th direction.
    Under this setting, the scalar form of the local Lax-Friedrichs flux becomes
    \begin{equation} \label{eq:ldg:num-flux:rusanov-cartesian}
        \Hat{f}_{k}(u_h^-, u_h^+) = \dfrac{1}{2}\left(f_k(u_h^-) + f_k(u_h^+)\right) - \dfrac{\alpha_{R,k}}{2}\left(u_h^+ - u_h^-\right), 
    \end{equation}
    where $f_k$ denotes the $k$-th component of the flux vector.
    The dissipation parameter $\alpha_{R,k}$ now acts per direction, estimating the maximum local characteristic speed specifically along the $k$-th coordinate axis.
    
    For the second-order terms, namely $\widehat{u}_h$ and $\widehat{\qv}_h$, the numerical fluxes must be selected differently, as the upwinding principle used for convection \joseanto{does not apply to} diffusion. 
    While several configurations exist in the literature, see e.g., \cite{Arnold2002}, we adopt the family of alternating numerical fluxes. 
    
    Consistent with our assumptions, we define these fluxes for each spatial direction $k$. 
    The alternating logic implies that whenever $\widehat{u}$ is selected from one side of the interface, the corresponding flux component $\widehat{q}_{h,k}$ must be chosen from the opposite side.
    This relationship is formalized as
    \begin{equation}
        \label{eq:ldg:alternating_flux}
        \widehat{u}_{h} = (1-\theta) u_{h}^-+\theta u_{h}^+,\quad \widehat{q}_{h,k} = \theta q_{h_k}^-+(1-\theta)q_{h,k}^+,\quad \theta\in \{0,1\}.
    \end{equation}
    This choice guarantees the stability and local solvability of the auxiliary variables, while providing optimal convergence results \cite{Cockburn1998,Bertoluzza2008}.

    \subsection{Local basis representation and degrees of freedom} 
    \label{ldg:num:basis}
    Unlike the globally continuous bases employed in continuous FEM, DG basis functions are not constrained to be continuous across element boundaries. 
    Leveraging the tensor-product structure of the meshes considered herein, the basis for the space $\bm{V}_{\bm{\alpha}}^{\nu}$ is constructed as the tensor product of the one-dimensional bases of $V_{l}^{\nu}$. 
    In the one-dimensional setting, a generic basis function $\phi^e:[a,b]\to\R$ is supported exclusively on the cell $I_e$.
    While the specific choice of polynomials is flexible, the set must satisfy the condition
    \begin{equation}
        \label{eq:ldg:num:basis:1d}
        \mathcal{P}_{\nu}\left(I_e\right) = \text{span}\left\{\phi_i^e\right\}_{i=1}^{\nu+1},\quad\forall\,e\in\{1,\dots,l\},
    \end{equation}
    where $\nu+1$ denotes the number of local degrees of freedom, \cite{Jaust2018}.
    Although different polynomial spaces could theoretically be employed on different elements, for simplicity, we assume a uniform polynomial degree throughout the domain. 

    The choice of polynomial bases directly influences the algebraic properties of the numerical scheme, including conditioning, mass matrix structure, and dispersion or dissipation characteristics, \cite{Canuto2007}. 
    Here, we focus on families that are orthogonal with respect to the $L^2$ inner product, as orthogonality generally ensures favorable spectral properties. 
    Literature typically distinguishes between modal bases, where each basis function represents a hierarchical mode associated with a polynomial degree (and generally does not interpolate nodal values), and nodal bases, where each basis function interpolates the solution at a specific node within the reference element $I=[-1,1]$. 
    
    The mapping between the reference element and a physical element $I_e$ is given by the affine diffeomorphism $\mathcal{T}^e:I\to I_e$,
    \begin{equation}
        \label{eq:ldg:num:basis:T}
        \mathcal{T}^e(\xi) = x_{e-1/2}+\dfrac{h}{2}\left(1+\xi\right),\quad\xi\in [-1,1], 
    \end{equation}
    with its inverse $(\mathcal{T}^e)^{-1}:I_{e}\to I$ defined by
    \begin{equation}
        \label{eq:ldg:num:basis:Tm1}
        \left(\mathcal{T}^e\right)^{-1}(x) = \dfrac{2}{h}\left(x-x_{e-1/2}\right)-1,\quad x\in [x_{e-1/2}, x_{e+1/2}].
    \end{equation}
    This transformation implies that all the integrals can be efficiently evaluated on the reference element, and thus, computed only once.

    In this work, we employ the Gauss-Legendre nodal basis on the reference interval $I$.
    Specifically, the nodes correspond to the Gauss-Legendre quadrature points, and the basis consists of Lagrange polynomials constructed upon these nodes. 
    This choice is particularly advantageous as it allows the use of the same quadrature rule to obtain a diagonal mass matrix. 
    Since this quadrature rule integrates polynomials up to degree $2\nu+1$ exactly, it yields exact evaluations for the inner products of basis functions (which are products of polynomials of degree at most $\nu$, resulting in a total degree of $2\nu$).
    Nevertheless, the presentation that follows remains general and, in principle, independent of the particular orthogonal basis employed.
    Consequently, the chosen basis will be explicitly referenced only when strictly needed.

    Building upon the preceding definitions, the basis for the multidimensional space $\bm{V}_{\bm{\alpha}}^{\nu}$ is constructed via the tensor product of the one-dimensional bases established in \eqref{eq:ldg:num:basis:1d}. 
    Thus, in $d$ dimensions, a generic basis function $\phi_{\iv}^{\bm{\beta}}:\Omega\to\R$ is supported exclusively on the element $Q_{\bm{\beta}}=I_{\beta_1}\otimes I_{\beta_2}\otimes\cdots\otimes I_{\beta_d}$, where $\bm{\beta}\leq\bm{\alpha}$. 
    The corresponding local polynomial space is defined as
    \begin{equation}
        \label{eq:ldg:num:basis:multid}
        \mathcal{P}_{\nu}\left(Q_{\bm{\beta}}\right) = \text{span}\left\{\phi^{\bm{\beta}}_{\iv}\right\}_{i_1,\dots,i_d=1}^{\nu+1} = \text{span}\left\{\phi_{i_1}^{\beta_1}\otimes\phi_{i_2}^{\beta_2}\otimes\cdots\otimes\phi_{i_d}^{\beta_d}\right\}_{i_1,\dots,i_d=1}^{\nu+1}.
    \end{equation}

    Under this formulation, the number of local degrees of freedom is given by $N_{\nu} = (\nu+1)^d$, resulting in a total of $N_{\nu}\cdot N_{\bm{\alpha}}$ degrees of freedom for each function unknown.
    This property, established for the one-dimensional case, extends naturally via tensorization.
    For instance, the mapping from the reference element $Q=[-1,1]^d$ to a physical element $Q_{\bm{\beta}}$ is defined as the tensor product of the 1D mappings, i.e., 
    \begin{equation}
        \mathcal{T}^{\bm{\beta}} = \mathcal{T}^{\beta_1}\otimes\mathcal{T}^{\beta_2}\otimes\cdots\otimes\mathcal{T}^{\beta_d}:Q\to Q_{\bm{\beta}},
    \end{equation}
    where each component $\mathcal{T}^{\beta_i},\,i\in\{1,\dots,d\}$ verifies \eqref{eq:ldg:num:basis:T}. 
    An analogous tensor-product construction applies to the inverse mapping.

    We represent the local discrete solutions $u_h$ and each component of the auxiliary field $\bm{q}_h=\left[q_{h,1},\dots,q_{h,d}\right]^T$ using the local tensor-product basis.
    On a generic element $Q_{\bm{\beta}}$, these expansions take the form
    \begin{subequations}
        \label{eq:ldg:num:representation}
        \begin{align}
            u_h(t,\xv)|_{Q_{\bm{\beta}}} &= \sum_{\iv\leq\nu+1}U_{\iv}^{\bm{\beta}}(t)\phi_{\iv}^{\bm{\beta}}(\xv) = \sum_{i_1=1}^{\nu+1}\cdots\sum_{i_d=1}^{\nu+1}U_{i_1,\dots,i_d}^{\beta_1,\dots,\beta_d}(t)\prod_{j=1}^d\phi_{i_j}^{\beta_j}(x_j), \label{eq:ldg:num:representation:u}\\
            q_{h,k}(t,\xv)|_{Q_{\bm{\beta}}} &= \sum_{\iv\leq\nu+1}W_{k,\iv}^{\bm{\beta}}(t)\phi_{\iv}^{\bm{\beta}}(\xv) =  \sum_{i_1=1}^{\nu+1}\cdots\sum_{i_d=1}^{\nu+1}W_{k, i_1,\dots,i_d}^{\beta_1,\dots,\beta_d}(t)\prod_{j=1}^d\phi_{i_j}^{\beta_j}(x_j). \label{eq:ldg:num:representation:q}
        \end{align}
    \end{subequations}
    Here, $U_{\iv}^{\bm{\beta}}(t)$ and $W_{k,\iv}^{\bm{\beta}}(t)$ denote the time-dependent modes for the primal variable and the $k$-th component of the auxiliary variable, respectively. 

    To introduce a matrix-based notation (and subsequent implementation), we represent the coefficients of the discrete solution and auxiliary variables in vector form, effectively flattening their inherited tensor structure.
    To this end, we define a bijective indexing map $\psi_{loc}:\left\{1,\dots,\nu+1\right\}\to\left\{1,\dots,N_{\nu}\right\}$ that assigns an integer to each local $d$-dimensional multi-index. 
    Adopting a row-major ordering: 
    \begin{equation}
        \label{eq:ldg:num:indexing:local}
        \psi_{loc}(\mathbf{i}) = 1 + \sum_{k=1}^d (i_k - 1)(\nu+1)^{d-k}, \quad \mathbf{i} \leq \nu+1.
    \end{equation}
    
    Henceforth, whenever a local vector or matrix entry is denoted with a multi-index, it is understood that the flattening map $\psi_{loc}$ is applied.
    With this convention, the local coefficients on element $\bm{\beta}$ are allocated into $d+1$ vectors, $\bm{U}^{\bm{\beta}}(t)$ and $\bm{W}^{\bm{\beta}}_k(t)$ for $k=1,\dots,d$, defined by 
    \begin{equation}
        \label{eq:ldg:num:local_vector_coeffs}
        [\bm{U}^{\bm{\beta}}(t)]_{\iv} = U^{\bm{\beta}}_{\iv}(t)=U^{\bm{\beta}}_{i_1,\dots,i_d}(t),\qquad [\bm{W}_k^{\bm{\beta}}(t)]_{\iv} = W^{\bm{\beta}}_{k,\iv}(t) = W^{\bm{\beta}}_{k,i_1,\dots,i_d}(t).
    \end{equation}
    Analogously, we impose a global ordering on the elements via a bijection $\zeta$, mapping the multi-index $\bm{\beta}$ to a linear index in $\left\{1,\dots,N_{\bm{\alpha}}\right\}$,
    \begin{equation}
        \label{eq:ldg:num:indexing:elements}
        \zeta(\bm{\beta}) = 1 + \sum_{k=1}^d (\beta_k - 1) \left( \prod_{j=k+1}^d \alpha_j \right), \quad \bm{\beta} \leq \bm{\alpha},
    \end{equation}
    where the empty product is taken as $1$.
    The local \eqref{eq:ldg:num:indexing:local} and element \eqref{eq:ldg:num:indexing:elements} mappings induces the global indexing map
    \begin{equation}
        \label{eq:ldg:num:indexing:global}
        \psi_{glob}(\iv,\bm{\beta})=(\zeta(\bm{\beta})-1)N_{\nu}+\psi_{loc}(\iv).
    \end{equation}
        
    Accordingly, the global coefficient vectors $\bm{U}(t)$ and $\bm{W}_k(t),\,k=1,\dots,d$ are assembled by concatenating the local representation defined in \eqref{eq:ldg:num:local_vector_coeffs},
    \begin{subequations}
        \label{eq:ldg:num:global_vector_coeffs}
        \begin{align}
            \bm{U}(t) &= \left[\bm{U}^1(t),\bm{U}^2(t),\cdots,\bm{U}^{N_{\bm{\alpha}}}(t)\right]^T, \label{eq:global_vector_coeffs:u}\\
            \bm{W}_k(t) &= \left[\bm{W}_k^1(t),\bm{W}_k^2(t),\cdots,\bm{W}_k^{N_{\bm{\alpha}}}(t)\right]^T. \label{eq:global_vector_coeffs:q}
        \end{align}
    \end{subequations}
    Here, we implicitly invoke the bijective mappings defined previously, identifying the integers $1,\dots,N_{\bm{\alpha}}$ with their corresponding multi-indices via $\zeta^{-1}$.
    
    Similarly, any spatial function $g(\bm{x})$, in particular the variable PDE coefficients, can be projected onto the space $\bm{V}_{\bm{\alpha}}^{\nu}$.
    This projection admits an expansion in the corresponding basis, enabling both the local representation \eqref{eq:ldg:num:local_vector_coeffs} and the global assembly \eqref{eq:ldg:num:global_vector_coeffs}, which will be utilized in the subsequent matrix construction.

    \subsection{Matrix system of equations} \label{ldg:num:matrix-eqs}
    After the necessary notation has been established, we proceed to reformulate the semi-discrete problem \eqref{eq:ldg:weak_form} in matrix form.
    This approach enables a clear and concise presentation of the spatial discretization scheme and effectively prepares the ground for the subsequent discussion on time discretization. 

    By testing the weak formulation \eqref{eq:ldg:weak_form:a} with basis functions $v=\phi^{\bm{\beta}}_{\iv},\,\iv\leq\nu+1$, we obtain the following time-dependent system for each element $Q_{\bm{\beta}},\,\bm{\beta}\leq\bm{\alpha}$,
    \begin{equation}
        \label{eq:ldg:local_system_eqs:a}
        \begin{split}
            \sum_{\jv\leq\nu+1}\dfrac{\text{d}}{\text{d}t} U_{\jv}^{\bm{\beta}}(t)\int_{Q_{\bm{\beta}}}\phi_{\jv}^{\bm{\beta}}\phi_{\iv}^{\bm{\beta}} &=
            \sum_{k=1}^d\Biggl(\sum_{\jv\leq\nu+1}F_{k,\jv}^{\bm{\beta}}(t)\int_{Q_{\bm{\beta}}}\phi_{\jv}^{\bm{\beta}}\partial_{x_k}\phi_{\iv}^{\bm{\beta}}-\oint_{\partial Q_{\bm{\beta}}}\widehat{\mathcal{F}}_{k}(t)\phi_{\iv}^{\bm{\beta}}\Biggr)\\
            &\quad+\sum_{k,l=1}^d a_{kl}\,\sum_{\jv\leq\nu+1}\Biggl\{-W_{l,\jv}^{\bm{\beta}}(t)\int_{Q_{\bm{\beta}}}\phi_{\jv}^{\bm{\beta}}\partial_{x_k}\phi_{\iv}^{\bm{\beta}}\Biggr.\\
            &\quad\quad+(1-\theta)W_{l,\jv}^{\bm{\beta}+\ev_k} (t)\int_{\Gamma_{k,+}^{\bm{\beta}}}\phi_{\jv}^{\bm{\beta}+\ev_k}\phi_{\iv}^{\bm{\beta}}\\
            &\quad\quad+W_{l,\jv}^{\bm{\beta}}\left(\theta\int_{\Gamma_{k,+}^{\bm{\beta}}}\phi_{\jv}^{\bm{\beta}}\phi_{\iv}^{\bm{\beta}}
            -(1-\theta)\int_{\Gamma_{k,-}^{\bm{\beta}}}\phi_{\jv}^{\bm{\beta}}\phi_{\iv}^{\bm{\beta}}\right)\\
            &\quad\quad-\theta W_{l,\jv}^{\bm{\beta}-\ev_k}
            \int_{\Gamma_{k,-}^{\bm{\beta}}}\phi_{\jv}^{\bm{\beta}-\ev_k}\phi_{\iv}^{\bm{\beta}}\Biggr\}\\
            &\quad+\sum_{\jv\leq\nu+1}S_{\jv}^{\bm{\beta}}(t)\int_{Q_{\bm{\beta}}}\phi_{\jv}^{\bm{\beta}}\phi_{\iv}^{\bm{\beta}},
        \end{split}
    \end{equation}
    where the interface convective flux term combines the contributions from both faces in direction $k$,
    \begin{equation} 
        \label{eq:ldg:local_system_eqs:a:aux-conv}
        \oint_{\partial Q_{\bm{\beta}}} \widehat{\mathcal{F}}_{k}(t) \phi_{\iv}^{\bm{\beta}} = \int_{\Gamma_{k,+}^{\bm{\beta}}} \widehat{f}_{k}(U^{\bm{\beta}}(t), U^{\bm{\beta}+\ev_k}(t)) \phi_{\iv}^{\bm{\beta}} - \int_{\Gamma_{k,-}^{\bm{\beta}}} \widehat{f}_{k}(U^{\bm{\beta}-\ev_k}(t), U^{\bm{\beta}}(t)) \phi_{\iv}^{\bm{\beta}},
    \end{equation}
    which corresponds to the Rusanov flux evaluated on the traces of $u_h$ at element interfaces.
    
    For notation convenience, we introduce the local numerical flux vector $\widehat{\bm{F}}^{\bm{\beta}}_k(t)$ collecting these interface contributions,
    \begin{equation}
        \label{eq:ldg:local-num-flux-vector}
        [\widehat{\bm{F}}_k^{\bm{\beta}}(t)]_{\iv} = \oint_{\partial Q_{\bm{\beta}}} \widehat{\mathcal{F}}_{k}(t) \phi_{\iv}^{\bm{\beta}}.
    \end{equation}
    Furthermore, the physical flux and the source term are projected onto  the finite element space $\bm{V}_{\bm{\alpha}}^{\nu}$ and represented locally by the vectors $\bm{F}_k^{\bm{\beta}}(t)$ and $\bm{S}^{\bm{\beta}}(t)$. 
    Their specific entries, $[\bm{F}_k^{\bm{\beta}}(t)]_{\jv} = F_{k,\jv}^{\bm{\beta}}(t)$ and $[\bm{S}^{\bm{\beta}}(t)]_{\jv}=S^{\bm{\beta}}_{\jv}(t)$, correspond to the expansion coefficients of the $k$-th physical flux component $f_k(u_h)$, and the source term, $s(u_h)$, in the local basis, respectively. 

    Regarding the weak equation of the second equation \eqref{eq:ldg:weak_form:b}, on structured tensor-product meshes, the problem naturally decouples into $d$ independent scalar relations, each associated with a specific auxiliary component $q_k$ for $k=1,\dots,d$. 
    Testing this $k$-th relation with the basis function  $w_{h,k}=\phi^{\bm{\beta}}_{\iv},\,\iv\leq\nu+1$ yields the following coupled local system,
    \begin{equation}
        \label{eq:ldg:local_system_eqs:b}
        \begin{split}
            \sum_{\jv\leq\nu+1}W^{\bm{\beta}}_{k,\jv}(t)\int_{Q_{\bm{\beta}}}\phi_{\jv}^{\bm{\beta}}\phi_{\iv}^{\bm{\beta}} &= 
            \sum_{\jv\leq\nu+1}\Biggl\{-U_{\jv}^{\bm{\beta}}(t)\int_{Q_{\bm{\beta}}}\phi_{\jv}^{\bm{\beta}}\partial_{x_k}\phi_{\iv}^{\bm{\beta}}\\
            &\quad+\theta U_{\jv}^{\bm{\beta}+\ev_k}(t)\int_{\Gamma_{k,+}^{\bm{\beta}}}\phi_{\jv}^{\bm{\beta}+\ev_k}\phi_{\iv}^{\bm{\beta}} \Biggr.\\
            &\quad +U_{\jv}^{\bm{\beta}}(t)\left((1-\theta)\int_{\Gamma_{k,+}^{\bm{\beta}}}\phi_{\jv}^{\bm{\beta}}\phi_{\iv}^{\bm{\beta}}-\theta\int_{\Gamma_{k,-}^{\bm{\beta}}}\phi_{\jv}^{\bm{\beta}}\phi_{\iv}^{\bm{\beta}}\right)\\
            &\quad-(1-\theta) U_{\jv}^{\bm{\beta}-
            \ev_k}(t)\int_{\Gamma_{k,-}^{\bm{\beta}}}\phi_{\jv}^{\bm{\beta}-\ev_k}\phi_{\iv}^{\bm{\beta}}\Biggr\}.
        \end{split}
    \end{equation}

    For the derivation of both local formulations, we employ the alternating numerical fluxes for $u_h$ and $q_{h,k}$ defined in \eqref{eq:ldg:alternating_flux}. 
    This choice induces a specific directional exchange between a cell and its immediate neighbors. 
    Under the assumption of periodic boundary conditions, the neighbor indices naturally wrap around the domain in each coordinate direction.
    Mathematically, the index connectivity for the $k$-th component is defined as
    \begin{equation}
        \label{eq:ldg:periodic-wrap}
        (\bm{\beta}+\ev_k)_k = \begin{cases}
            \beta_k+1,&\beta_k<\alpha_k,\\
            1,&\beta_k=\alpha_k,
        \end{cases}
        \qquad (\bm{\beta}-\ev_k)_k = \begin{cases}
            \beta_k-1,&\beta_k>1,\\
            \alpha_k,&\beta_k=1.
        \end{cases}
    \end{equation}
    Conversely, for non-periodic boundary conditions, the reference to an outward neighbor is replaced by the prescribed boundary data.
    In such cases, the numerical flux must be modified to enforce the condition consistently.
    We provide more information in Section \ref{ldg:num:bcs}.

    The local formulations \eqref{eq:ldg:local_system_eqs:a} and \eqref{eq:ldg:local_system_eqs:b} can be cast into a compact matrix form. 
    This representation facilitates the assembly of the global system of differential equations from element-wise contributions. 
    We begin by defining the local mass matrix on a generic element $Q_{\bm{\beta}},\,\bm{\beta}\leq\bm{\alpha}$, as $\bm{M}^{\bm{\beta}}\in\mathcal{M}_{N_{\nu}}(\R)$, with entries given by
    \begin{equation}
        \label{eq:ldg:matrix:mass:local}
        \left[\M^{\bm{\beta}}\right]_{\iv,\jv} = \int_{Q_{\bm{\beta}}} \phi_{\jv}^{\bm{\beta}}\phi_{\iv}^{\bm{\beta}}.  
    \end{equation}
    Since the basis functions $\phi_{\iv}^{\bm{\beta}}$ are supported exclusively on $Q_{\bm{\beta}}$, the resulting global mass matrix $\bm{M}\in\mathcal{M}_{N_{\nu}N_{\bm{\alpha}}}(\R)$ exhibit a block-diagonal structure,
     \begin{equation}
        \label{eq:ldg:matrix:mass:glob}
        \M = \text{diag}\left(\M^{1}, \M^{2},\dots,\M^{N_{\bm{\alpha}}}\right), 
    \end{equation}       
    where the global element ordering follows the mappings introduced in Section \ref{ldg:num:basis}. 
    Moreover, due to the orthogonal nodal basis adopted (Gauss-Legendre), each local mass matrix, and therefore $\M$, is strictly diagonal, computationally advantageous as it enables efficient, trivial inversion.

    An analogous construction applies to the derivative matrices. 
    For each spatial direction $k=1,\dots,d$, we define the local derivative matrix with respect to the coordinate $x_k$, denoted by $\bm{C}^{\bm{\beta}}_{k}\in\mathcal{M}_{N_{\nu}}(\R)$ via the entries
    \begin{equation}
        \label{eq:ldg:matrix:der:loc}
        [\bm{C}^{\bm{\beta}}_{k}]_{\iv,\jv} = \int_{Q_{\bm{\beta}}}\phi_{\jv}^{\bm{\beta}}\partial_{x_k}\phi_{\iv}^{\bm{\beta}}.
    \end{equation}
    Due to the strictly local support of the basis functions, the corresponding global derivative matrix inherits a block-diagonal structure as the mass matrix. 
    Thus, the global derivative matrix $\bm{C}_{k}\in\mathcal{M}_{N_{\nu}N_{\bm{\alpha}}}(\R)$ is assembled as
    \begin{equation}
        \label{eq:ldg:matrix:der:glob}
        \bm{C}_k = \text{diag}\left(\bm{C}^{1}_{k}, \bm{C}^{2}_{k},\dots,\bm{C}^{N_{\bm{\alpha}}}_{k}\right),
    \end{equation}
    adhering to the global indexing convention established previously.

    Finally, we consider the matrices encoding the alternating flux logic. 
    Unlike the mass and derivative matrices, these operators require extended stencils that couple neighboring cells, thereby causing the global matrices to lose their block-diagonal form and instead exhibit a block-banded structure.
    As a preparatory step, we define the forward, $\mathcal{S}_k^+$, and backward, $\mathcal{S}_k^-$, shift operators in direction $k$.
    For a generic local vector $X^{\bm{\beta}}$ representing coefficients on $Q_{\bm{\beta}}$, we define
    \begin{equation}
        \label{eq:ldg:shif_ops}
        [\mathcal{S}_{k}^+X^{\bm{\beta}}]_{\jv} = X^{\bm{\beta}+\ev_k}_{\jv},\quad [\mathcal{S}_{k}^{-}X^{\bm{\beta}}]_{\jv} = X_{\jv}^{\bm{\beta}-\ev_k},\quad \jv\leq\nu+1, 
    \end{equation}
    assuming the periodic wrapping of indices as described in \eqref{eq:ldg:periodic-wrap}. 
    Next, we define the trace mass matrices.
    These are element-local matrices $\M_{k,\pm}^{\bm{\beta},ext},\,\M_{k,\pm}^{\bm{\beta},int}\in\mathcal{M}_{N_{\nu}}(\R)$ that compute the inner products on the boundary faces.
    Specifically,
    \begin{equation}
        \label{eq:ldg:matrix:trace:local}
        [\M_{k,\pm}^{\bm{\beta},ext}]_{\iv,\jv} = \int_{\Gamma_{k,\pm}^{\bm{\beta}}}\phi_{\jv}^{\bm{\beta\pm\ev_k}}\phi_{\iv}^{\bm{\beta}},\quad [\M_{k,\pm}^{\bm{\beta},int}]_{\iv,\jv} = \int_{\Gamma_{k,\pm}^{\bm{\beta}}}\phi_{\jv}^{\bm{\beta}}\phi_{\iv}^{\bm{\beta}}.
    \end{equation}
    Using these building blocks, we construct the alternating operator $\bm{A}^{\bm{\beta}}_{k,\theta}$, parameterized by $\theta\in\{0,1\}$.
    This operator acts on a local coefficient vector $X^{\bm{\beta}}$ by assembling the contributions from the element and its neighbors as follows,
    \begin{equation}
    \label{eq:ldg:matrix:alternating:local}
        \begin{split}
        [\bm{A}_{k,\theta}^{\bm{\beta}}X^{\bm{\beta}}]_{\iv} &= \theta\M^{\bm{\beta},ext}_{k,+}\left(\mathcal{S}^+_kX^{\bm{\beta}}\right) \\
        &\quad + \left((1-\theta)\M^{\bm{\beta},int}_{k,+} - \theta\M^{\bm{\beta},int}_{k,-}\right)X^{\bm{\beta}} \\
        &\quad - (1-\theta)\M^{\bm{\beta},ext}_{k,-}\left(\mathcal{S}_k^{-}X^{\bm{\beta}}\right).
        \end{split}
    \end{equation}
    The complementary alternating operator (required for the LDG formulation) is obtained by substituting the parameter $\theta$ with $1-\theta$ in the expression above. 

    Building upon the element-level operators introduced above, we assemble the global alternating matrix $\bm{A}_{k,\theta}$ by stacking the local unknowns according to the chosen global element ordering. 
    The resulting matrix structure is characterized as follows:
    \begin{itemize}
        \item A block-diagonal construction that consists of the interior face contributions, assembled element-wise as $(1-\theta)\M_{k,+}^{\bm{\beta},int}-\theta\M_{k,-}^{\bm{\beta},int}$.
        \item Off-diagonal bands. 
        For each spatial direction $k$, two off-diagonal block bands emerge due to neighbour coupling.
        An upper band coupling each element to its forward neighbour via $\theta\M_{k,+}^{\bm{\beta},ext}$, and a lower band coupling to its backward neighbor via $-(1-\theta)\M_{k,-}^{\bm{\beta},ext}$.
    \end{itemize}
    
    Under periodic boundary conditions, these bands wrap around to connect the first and last element blocks, yielding a block-circulant structure.
    Conversely, for non-periodic boundary conditions, the wrapping blocks are omitted.
    Instead, the boundary face contributions are incorporated into the right-hand side vector.
    Finally, aggregating the contributions for all spatial directions $k=1,\dots,d$ produces a highly sparse matrix exhibiting up to $2d$ off-diagonal block bands.

    Lastly, by combining the derivative matrix defined in \eqref{eq:ldg:matrix:der:loc} with the alternating boundary operator from \eqref{eq:ldg:matrix:alternating:local}, we construct the LDG alternating derivative operators with respect to the coordinate $x_k$.
    The local and global representations are defined, respectively, as
    \begin{equation}
        \label{eq:ldg:matrix:der-ldg}
        \bm{D}_{k, \theta}^{\bm{\beta}} = -\bm{C}_k^{\bm{\beta}}+\bm{A}_{k,\theta}^{\bm{\beta}},\quad \bm{D}_{k, \theta} = -\bm{C}_k+\bm{A}_{k,\theta}.
    \end{equation}
    Conceptually, this operator represents the discrete weak derivative obtained via Green's identity on discontinuous elements.
        
    Leveraging the definitions established above, we can reformulate the differential system of local differential equations \eqref{eq:ldg:local_system_eqs:a}-\eqref{eq:ldg:local_system_eqs:b} into the following compact matrix form,
    \begin{subequations}
    \label{eq:ldg:local_matrix_system}
        \begin{align}
            \begin{split}
                \mathbf{M}^{\bm{\beta}}\dfrac{\text{d}}{\text{d}t} \bm{U}^{\bm{\beta}}(t) &= \sum_{k=1}^d\left(\mathbf{C}^{\bm{\beta}}_k\bm{F}_k^{\bm{\beta}}(t)-\widehat{\bm{F}}_k^{\bm{\beta}}(t)\right) \\
                &\quad + \sum_{k,l=1}^d a_{kl}\mathbf{D}_{k,(1-\theta)}^{\bm{\beta}}\bm{W}_{l}^{\bm{\beta}}(t)+\mathbf{M}^{\bm{\beta}}\bm{S}^{\bm{\beta}}(t),
            \end{split} \label{eq:ldg:local_matrix_system:a} \\
            \mathbf{M}^{\bm{\beta}}\bm{W}_k^{\bm{\beta}}(t) &= \mathbf{D}_{k,\theta}^{\bm{\beta}}\bm{U}^{\bm{\beta}}(t),\quad k=1,\dots,d.\label{eq:ldg:local_matrix_system:b}
        \end{align}
    \end{subequations}
        
    The global assembly of matrix subsystems \eqref{eq:ldg:local_matrix_system} yields the ODE system
    \begin{subequations}
    \label{eq:ldg:global_matrix_system}
        \begin{align}
            \label{eq:ldg:global_matrix_system:a}
            \M\dfrac{\text{d}}{\text{d}t}\bm{U}(t)&=\sum_{k=1}^d\left(\bm{C}_k\bm{F}_k(t)-\widehat{\bm{F}}_k(t)\right)+\sum_{k,l=1}^d a_{kl}\bm{D}_{k,(1-\theta)}\bm{W}_l(t)+\M\bm{S}(t), \\
            \M \bm{W}_k(t) &= \bm{D}_{k,\theta}\bm{U}(t),\quad k=1,\dots,d. \label{eq:ldg:global_matrix_system:b} 
        \end{align}
    \end{subequations}
    This system governs the temporal evolution of the global coefficient vector $\bm{U}(t)$, which represents the expansion of the approximate solution $u_h$ in the chosen finite element basis. 
    In addition, the auxiliary vectors $\bm{W}_k(t),\,k=1,\dots,d$, which represent the expansion of $\nabla_{\bm{x}}u$ in the chosen basis, must be resolved at each time step. 
    Consequently, the total dimensionality of the problem in this formulation is $(d+1)N_{\nu}N_{\bm{\alpha}}$.
    
    Consistent with the \emph{local} nature of the LDG method, the auxiliary variables can be explicitly eliminated by static condensation, \cite{Jaust2018}.
    Specifically, since the mass matrix $\bm{M}$ is block-diagonal (and strictly diagonal with our basis choice), we can trivially invert it to solve for $\bm{W}_k(t)$ in \eqref{eq:ldg:global_matrix_system:b},
    \begin{equation}
        \bm{W}_k(t) = \bm{M}^{-1}\bm{D}_{k,\theta}\bm{U}(t),\quad k=1,\dots,d.
    \end{equation}
    Substituting this expression into \eqref{eq:ldg:global_matrix_system:a} yields the equivalent condensed ODE system
     \begin{equation}
        \label{eq:ldg:condensed_matrix_system} 
        \bm{M}\dfrac{\text{d}}{\text{d}t}\bm{U}(t) = \sum_{k=1}^d\left(\bm{C}_k\bm{F}_k(t)-\widehat{\bm{F}}_k(t)\right)+\bm{K}\bm{U}(t) +\bm{M}\bm{S}(t),
    \end{equation}
    where $\bm{K}\in\mathcal{M}_{N_{\nu}N_{\alpha}}(\R)$ represents the discrete second-order differential operator, defined as 
    \begin{equation}
        \label{eq:ldg:matrix:stiffness}
        \bm{K} = \sum_{k,l=1}^d a_{kl}\bm{D}_{k,(1-\theta)}\bm{M}^{-1}\bm{D}_{l,\theta}.
    \end{equation}
    This operator $\bm{K}$ should be viewed as the LDG analogue of the classical stiffness matrix found in continuous finite elements. 
        
    This reduction corresponds to the Schur-complement condensation of the block system involving $\left(\bm{U},\bm{W}\right)$.
    Crucially, it decreases the total number of global unknowns ($d\cdot N_{\nu}N_{\bm{\alpha}}$) to $N_{\nu}N_{\bm{\alpha}}$, retaining only those associated with the primal variable $\bm{U}$. 
    \joseanto{This is particularly} advantageous in the context of IMEX schemes for two primary reasons. 
    First, the condensation yields a significantly smaller system while preserving a sparse, block-banded structure; and second, the resulting stiffness matrix $\bm{K}$ recovers the expected properties of the second-order operator.
    This property is inherited by the linear system arising from the implicit time integration step, thereby improving the convergence rates of iterative linear solvers and enhancing preconditioner efficiency.

    \subsection{Treatment of non-periodic boundary conditions and non-constant diffusion coefficients} \label{ldg:num:bcs}
    So far, the derivation has proceeded under the simplifying assumptions of constant diffusion coefficients and periodic boundary conditions.
    However, the proposed framework extends naturally to more general settings with only minor modifications. 
    In this subsection, we outline the specific adjustments to accommodate spatially varying diffusion coefficients and boundary conditions.

    When the diffusion matrix exhibits spatial dependence, that is, $A = (a_{kl}(\xv))\in\mathcal{M}_d(\R)$, the primary modification involves revisiting the construction of the LDG derivative matrices defined in \eqref{eq:ldg:matrix:der-ldg}.
    In the general case, the necessary adaptations follow a consistent pattern.
    We compute the local matrices element-wise by explicitly incorporating the spatially varying diffusion coefficient into the quadrature evaluation at the local integration points. 
    Accordingly, rather than factoring out a constant coefficient, we assemble specific derivative matrices for each diffusion tensor component and spatial direction.
    This procedure reflects the matrix-vector product structure of the divergence operator applied to the variable diffusive flux.
    Thus, the weighted local derivative matrix components are defined as 
    \begin{equation}
        \label{eq:ldg:matrix:der:non_constant}
        [\bm{C}_{k,l}^{\bm{\beta}}]_{\iv,\jv} = \int_{Q_{\bm{\beta}}} a_{k,l}(\xv)\phi_{\bm{j}}^{\bm{\beta}}(\xv)\partial_{x_k}\phi_{\bm{i}}^{\bm{\beta}}(\xv)\text{d}\xv,\quad k,l=1,2,\dots,d,
    \end{equation}
    where the dependency on $\xv$ is made explicit to highlight the spatial variability of the diffusion coefficient within the integral. 
    
    Analogously, we must redefine the trace operators to incorporate the diffusion matrix components at the element interfaces.
    Thus, the weighted trace mass matrices are given by 
    \begin{equation}
    \label{eq:ldg:matrix:trace:local:non_constant}
        \begin{split}
        [\M_{k, l,\pm}^{\bm{\beta},ext}]_{\iv,\jv} &= \int_{\Gamma_{k,\pm}^{\bm{\beta}}}a_{k,l}(\xv)\phi_{\jv}^{\bm{\beta}\pm\ev_k}(\xv)\phi_{\iv}^{\bm{\beta}}(\xv)\,\text{d}\Gamma, \\
        [\M_{k,l,\pm}^{\bm{\beta},int}]_{\iv,\jv} &= \int_{\Gamma_{k,\pm}^{\bm{\beta}}}a_{k,l}(\xv)\phi_{\jv}^{\bm{\beta}}(\xv)\phi_{\iv}^{\bm{\beta}}(\xv)\,\text{d}\Gamma,
        \end{split} \quad k,l=1,\dots,d.
    \end{equation}
    These definitions allow us to construct the alternating matrices $\bm{A}_{k,l,\theta}^{\bm{\beta}}$ by simply substituting these weighted trace operators into the structure defined in \eqref{eq:ldg:matrix:alternating:local}.
    
    Upon global assembly of the weighted derivative matrices $\bm{C}_{k,l}$ and alternating matrices $\bm{A}_{k,l,\theta}$, we define the weighted LDG alternating derivative operator associated with coordinate $x_k$ and the diffusion coefficient $a_{k,l}$ as 
    \begin{equation}
        \label{eq:ldg:matrix:der-ldg:non-constant}
        \bm{D}_{k, l, \theta} = -\bm{C}_{k,l}+\bm{A}_{k,l,\theta},\quad k,l=1,2,\dots,d,
    \end{equation}
    which yields the generalized stiffness matrix $\bm{K}$ formulated as
    \begin{equation}
        \label{eq:ldg:matrix:stiffness:non_constant}
        \bm{K} = \sum_{k,l=1}^d\bm{D}_{k,l,(1-\theta)}\bm{M}^{-1}\bm{D}_{l,\theta},
    \end{equation}
    where $\bm{D}_{l,\theta}$ represents the standard derivative operator defined in \eqref{eq:ldg:matrix:der-ldg}. 
    We retain the notation $\bm{K}$ for this operator to emphasize the algebraic analogy with the constant coefficient case.
    To address non-periodic boundary conditions, we must first categorize the enforced constraints, as each type requires a distinct numerical treatment.
    Assuming the \joseanto{previously defined} lexicographical \joseanto{convention}, we adopt the fixed left-alternating flux configuration, $\theta=0$.
    This choice entails no loss of generality, as the opposite strategy, $\theta=1$, follows an analogous \joseanto{derivation}.
    
    Under this convention, the stability of the LDG method \joel{relies on a careful algebraic distribution of the \joseanto{boundary constraints}.} 
    \joel{In the numerical experiments presented in this work, \joseanto{we primarily} employ Neumann and purely outflow boundary conditions.}
    \joel{Nevertheless, to provide a comprehensive framework capable of accommodating alternative setups, such as those discussed in \cite{Wang2018}, we also formalize the derivation of the Dirichlet load vectors.}

    The enforcement of these non-periodic boundary conditions modifies the system through the introduction of  specific right-hand side load vectors.
    Because of the asymmetry in the alternating flux, the algebraic treatment differs between the lower and upper boundaries.
    
    First, regarding standard Dirichlet conditions on the lower boundaries (where $\beta_k=1$), the exterior trace evaluates the prescribed condition, which we denote by $g_{k}^u(t,\xv)$.
    This generates a local source term $\bm{b}_k^u(t)$, defined as
    \begin{equation}
        \label{eq:ldg:vector:dirchlet:local}
       \left[\bm{b}_{k}^u(t)\right]_{\iv} = \begin{cases}\displaystyle
           -\int_{\Gamma_{k,-}^{\bm{\beta}}}g_{k}^u(t,\xv)\phi_{\iv}^{\bm{\beta}}(\xv)\,\text{d}\Gamma,& \beta_k=1,\\
           0,&\text{otherwise}.
       \end{cases} 
    \end{equation}
    Since these terms appear in the auxiliary equation for $\bm{W}_k(t)$, they enter the final condensed system through the diffusion operator, contributing to the global vector $\bm{B}_u(t)$. 

    \joel{In the present work, however, the physical \joseanto{boundary constraints} are imposed \joseanto{primarily through} the gradient variables. 
    Consequently, we prescribe an outflow condition for the primal variable $u$ across the lower boundaries.
    Under the left-alternating convection, the numerical flux at the lower boundaries would evaluate the exterior trace.
    To enforce a transmissive behavior, we set:
    \begin{equation*}
        \widehat{u}_h = u_h^+\quad \text{on } \Gamma_{k,-}^{\bm{\beta}} \text{ with } \beta_k = 1,\, k\in\{1,\dots,d\}.
    \end{equation*}
    Because \joseanto{these} transmissive conditions \joseanto{rely} on the internal state and \joseanto{are absorbed} implicitly into the discrete differential operators, no external data is introduced.
    Therefore, the associated explicit load vector vanishes, yielding $\bm{B}_u(t) = \bm{0}$.
    Algebraically, this implies that the backward dependence, given by the shift operator $\mathcal{S}_k^-$ vanishes.
    }

    \joel{Second, regarding the Neumann conditions imposed on the gradient auxiliary variables, the numerical treatment must distinguish between the lower and upper boundaries of the domain.
    On the upper boundaries (where $\beta_k=\alpha_k$), the left-alternating choice naturally evaluates the auxiliary numerical flux using the exterior trace.
    Thus, the numerical flux selects the prescribed exterior gradient, $g_{k,+}^q(t,\xv)$.
    Conversely, to enforce a Neumann condition on the lower boundaries (where $\beta_k=1$), the standard alternating structure in the $k$-th direction must be explicitly broken.
    \joseanto{Instead of} relying on a ghost cell, the prescribed physical datum, $g_{k,-}^q(t,\xv)$, is injected to evaluate the normal flux.
    }
    This generates a local term $\bm{b}_{k,l}^q(t,\xv)$, which can be generalized for both boundaries as:
    \begin{equation}
        \label{eq:ldg:vector:neumann:local}
         \left[\bm{b}_{k,l}^q(t)\right]_{\iv} = \begin{cases}\displaystyle
            -\int_{\Gamma_{k,-}^{\bm{\beta}}} a_{k,l}(\xv)g_{k,-}^q(t,\xv)\phi_{\iv}^{\bm{\beta}}(\xv),\text{d}\Gamma,& \beta_k=1,\\
           \displaystyle\int_{\Gamma_{k,+}^{\bm{\beta}}}a_{k,l}(\xv)g_{k,+}^q(t,\xv)\phi_{\iv}^{\bm{\beta}}(\xv)\,\text{d}\Gamma,& \beta_k=\alpha_k,\\
           0,&\text{otherwise}.
       \end{cases} 
    \end{equation} 
    Unlike the Dirichlet terms, these contributions are directly incorporated into the primal equation, yielding the global Neumann load vector $\bm{B}_q(t)$.

    \joel{
    Finally, for boundaries where the outflow is required, we enforce the behavior for the auxiliary variables by assuming continuity across the domain. 
    Under the left alternating, evaluating the numerical flux at $\beta_k=\alpha_k$ relies on the exterior trace, so we follow the same approach given for the outflow of the primal variable.
    We set this exterior trace \joseanto{equal to} the interior trace:
    }  
    \begin{equation*}
        q_{h,k}^+=q_{h,k}^-\quad \text{on } \Gamma_{k,+}^{\bm{\beta}} \text{ with } \beta_k=\alpha_k.
    \end{equation*}
    Consequently, since no external boundary data is introduced, the corresponding Neumann load vector is zero, $\bm{B}_q(t)=0$.
    Algebraically, for those boundary elements the forward neighbor dependence, represented by the shift operator $\mathcal{S}_k^+$, vanishes and is then resolved by local terms.
    
    Consequently, by incorporating the aforementioned boundary effects, the fully discrete condensed ODE system is generalized to:
    \begin{equation}
        \label{eq:ldg:condensed_matrix_system:general}
        \bm{M}\dfrac{\text{d}}{\text{d}t}\bm{U}(t) = \sum_{k=1}^d\left(\bm{C}_k\bm{F}_k(t)-\widehat{\bm{F}}_k(t)\right)+\bm{K}\bm{U}(t)+ \underbrace{\bm{B}_u(t) + \bm{B}_q(t)}_{\text{Boundary loads}} +\bm{M}\bm{S}(t),
    \end{equation}
    where the global boundary vectors are assembled as
    \begin{equation}
        \label{eq:ldg:vectors:global}
        \bm{B}_u(t)=\sum_{k,l=1}^d\bm{D}_{k,l,(1-\theta)}\bm{M}^{-1}\bm{b}_{l}^u(t),\qquad \bm{B}_q(t)=\sum_{k,l=1}^d\bm{b}_{k,l}^q(t).
    \end{equation}

    \joel{As previously established, under the specific boundary conditions prescribed for the primal variable in the numerical experiments of Section \ref{ldg:experiments}, the Dirichlet load vector vanishes.
    However, we deliberately retain the fully generalized structure in \eqref{eq:ldg:condensed_matrix_system:general}.
    This ensures that the IMEX time integration scheme, detailed in the following section, is formulated in \joseanto{a universal} manner.
}
    
     \section{IMEX-RK time integration and fully discrete scheme} \label{ldg:imex}
    Following spatial discretization, the PDE is reduced to a high-dimensional ODE system. 
    The size of this system scales exponentially with the number of spatial dimensions, a phenomenon termed the curse of dimensionality, as both the number of elements and the degrees of freedom per element increase. 
    For the convection-diffusion problem addressed herein, standard explicit time integrators, e.g., explicit RK schemes, become inefficient on refined meshes.
    This is due to the semi-discrete system becoming increasingly stiff, \cite{Calvo2001}.
    In particular, the stability constraint imposed by the diffusion term dictates that the time step must scale quadratically with the mesh size, that is, $\Delta\text{t} \propto \min_{k}\{h_k^2\}$. 
    This restriction makes explicit integration computationally prohibitive for high-resolution simulations.
    
    Conversely, fully implicit methods effectively remove this severe time-step restriction. However, in the presence of nonlinear convective fluxes \joseanto{or source terms}, they require, at least, the solution of a coupled nonlinear system at every time step. 
    For large-scale problems, particularly those employing spectral discretizations, this approach often becomes computationally prohibitive.
    Nonlinear solvers are significantly less efficient than linear solvers, rendering fully implicit schemes impractical for this class of problems, \cite{Calvo2001}.

    To circumvent the limitations inherent in fully implicit or fully explicit approaches for \joseanto{the} semi-discrete nonlinear LDG system \eqref{eq:ldg:condensed_matrix_system:general}, we adopt IMEX-RK methods. 
    In this framework, the stiff parabolic/diffusive terms are treated implicitly, while the hyperbolic/convective terms are handled explicitly. 
    This strategy provides a dual advantage. 
    On the one hand, it relaxes the restrictive diffusive stability constraint, as the implicit treatment of diffusion yields linear elliptic algebraic systems that are, in principle, amenable to efficient iterative solvers.
    On the other hand, it avoids the computational cost of implicitly solving nonlinear convective terms by advancing them explicitly \cite{Wang2015}. 
        
    Let $\{t^n=n\Delta t\}_{n=0}^{N_T}$ denote a uniform partition of the time interval $[0,T]$, with time step $\Delta t$ such that $N_T\Delta t = T$. 
    Although the time step could actually change from step to step, we take it to be constant here \joseanto{to simplify the presentation}. 
    Given $\bm{U}^n\approx\bm{U}(t_n)$, we seek to compute the numerical solution of \eqref{eq:ldg:condensed_matrix_system:general} at the subsequent time level $t^{n+1}$.
    To illustrate the fundamental algorithmic challenges, we begin our discussion with the simplest first-order IMEX-RK scheme, the forward-backward Euler method, which involves no intermediate stages. 
    This baseline scheme provides a clear template for coupling explicit and implicit contributions, as higher-order IMEX-RK methods can be conceptualized as a stage-wise linear combination of such forward-backward updates. 
    
    Assuming consistent notation for all time-dependent vectors, the \joseanto{forward-backward   Euler IMEX}  method applies a first-order approximation to the time derivative, evaluating the diffusive contribution implicitly at time $t^{n+1}$ and the convective and source term contributions at time $t^n$,
        \begin{equation}
            \bm{M}\dfrac{1}{\Delta t}\left(\bm{U}^{n+1}-\bm{U}^n\right) = \sum_{k=1}^d\left(\bm{C}_k\bm{F}_k^n-\widehat{\bm{F}}_k^n\right)+\bm{K}\bm{U}^{n+1} +\bm{M}\bm{S}^n+\bm{B}_{bcs}^{n+1},
        \end{equation}
        where $\bm{B}_{bcs}^{n+1}=\bm{B}_u(t^{n+1})+\bm{B}_q(t^{n+1})$ represents the total boundary vector evaluated at the new time level.
        Rearranging the terms, we obtain the linear system to be solved at each time step,
        \begin{equation}
            \label{eq:ldg:imex:forward_backward_euler}
            \left(\bm{M}-\Delta t\bm{K}\right)\bm{U}^{n+1} = \bm{M}\bm{U}^n +\Delta t\,\left(\bm{R}_{conv}^n+\bm{B}_{bcs}^{n+1}\right),
        \end{equation}
        where $\bm{R}_{conv}^n$ aggregates the explicit convective and source contributions computed from the previous state.

        Consequently, advancing the solution to time $t^{n+1}$ requires solving the linear system \eqref{eq:ldg:imex:forward_backward_euler} governed by the coefficient matrix $\bm{M}-\Delta t\bm{K}$.
        It is worth noting that this system matrix is algebraically analogous to those arising in elliptic problems, specifically Helmholtz-type operators. 
        \joel{Importantly, the condensed formulation inherits key properties of the second-order operator, such as its coercive nature, retaining the favorable behavior to use efficient iterative solvers.}

        \joel{Given the substantial size of the resulting linear systems, particularly in high-dimensional settings, the use of iterative solvers becomes indispensable.
        In this work, we employ Krylov subspace methods \cite{Golub2013}, which are highly effective for large, sparse algebraic systems.
        This choice is further justified as they seamlessly accommodate matrix-free implementations, enabling computations in higher dimensions by mitigating memory constrains.
        Specific details regarding the implementation of these solvers, including preconditioning strategies for assembled matrices and their unconditioned matrix-free counterparts, are discussed alongside the numerical experiments in Section \ref{ldg:experiments}.
        }

        Having addressed the linear algebra requirements of the time-marching scheme, we now introduce the general IMEX-RK formulation. 
        As in standard RK methods, the coefficients and stages of the method are encoded via Butcher tableaux. 
        The scheme is constructed by coupling two distinct methods:
        \begin{itemize}
            \item For the diffusive component, we employ an $s$-stage diagonally implicit RK (DIRK) scheme characterized by coefficients $G\in\mathcal{M}_s(\R)$, and vectors $c, b\,\in\R^s$.
            \item For the convective and source term, we utilize an $(s+1)$-stage explicit RK scheme defined by the extended abscissae $\hat{c}=(0,c)^T$ and coefficients $\hat{G}\in\mathcal{M}_{s+1}(\R),\hat{b}\in\R^{s+1}$
        \end{itemize}
        To unify these methods into a single framework, the $s$-stage DIRK tableau is embedded into an $(s+1)$-stage structure. 
        This is achieved by padding the implicit matrix with a zero first row and column, thereby yielding compatible implicit-explicit tableaux suitable for the IMEX coupling, \cite{Ascher1997},

        Applying the general $s$-stage IMEX-RK discretization, the solution of the system \eqref{eq:ldg:condensed_matrix_system:general} advances from time $t^n$ to $t^{n+1}=t^n+\Delta t$ via the following multi-stage procedure
         \begin{subequations}
        \label{eq:ldg:imex}
        \begin{align}
            \bm{U}^{n,1} &= \bm{U}^n, \\
            \begin{split}
                \left(\bm{M}-\gamma\Delta t\bm{K}\right)\bm{U}^{n,i} &= \bm{M}\bm{U}^n + \gamma\Delta t \bm{B}_{bcs}^{n,i} \\
                &\quad+ \Delta t\Biggl[ \sum_{j=2}^{i-1} g_{ij}\left(\bm{K}\bm{U}^{n,j}+\bm{B}_{bcs}^{n,j}\right) + \sum_{j=1}^{i-1}\hat{g}_{ij}\bm{R}_{conv}^{n,j} \Biggr],
            \end{split} \quad 2\leq i\leq s+1,\\
            \bm{U}^{n+1} &= \bm{U}^n + \Delta t \bm{M}^{-1}\Biggl(\sum_{i=2}^{s+1}b_i\left(\bm{K}\bm{U}^{n,i}+\bm{B}_{bcs}^{n,i}\right)+ \sum_{i=1}^{s+1}\hat{b}_i\bm{R}_{conv}^{n,i}\Biggr),
            \end{align}
        \end{subequations}
        where $\bm{U}^{n,i}$ denote the intermediate stage solution, whilst $\bm{R}_{conv}^{n,i}$ and $\bm{B}_{bcs}^{n,i}$ represent the explicit convective/source and boundary terms evaluated at stage $i$, respectively. 
        Here, $\gamma$ denotes the standard diagonal entry of the DIRK tableau, and it is worth recalling that, in the spectral case employed here, the mass matrix $\bm{M}$ is strictly diagonal. The application of its inverse, required in the update step, is computationally trivial, preserving the overall efficiency of the scheme.

        Two specific IMEX schemes are considered in this work. 
        First, the L-stable, two-stage, second-order DIRK$(2,2,2)$ method proposed in \cite{Ascher1997} with the coefficients defined by the following Butcher tableau, where the implicit scheme is shown on the left and the explicit one on the right,
        \begin{equation}
            \label{eq:ldg:imex:tableau:dirk2}\renewcommand{\arraystretch}{1.3}
            \begin{array}{c|ccc|ccc}
            0        & 0 &   0       &    0     &  0   &      0      & 0              \\
            \gamma   & 0 & \gamma    &    0     &  \gamma   &       0  &  0          \\
            1        & 0 & 1-\gamma  & \gamma   &    \delta &      1-\delta    & 0   \\
            \hline
                     & 0 & 1-\gamma  & \gamma   & \delta  & 1-\delta       & 0 
            \end{array}
        \end{equation}
        with parameters $\gamma=1-\dfrac{\sqrt{2}}{2}$ and $\delta=1-\dfrac{1}{2\gamma}$. 
        Under this setting, the general IMEX scheme \eqref{eq:ldg:imex} reads
        \begin{subequations}
            \label{eq:ldg:imex:dirk2}
            \begin{align}
                \bm{U}^{n,1} &= \bm{U}^n, \\
                \left(\bm{M}-\gamma\Delta t\bm{K}\right)\bm{U}^{n,2} &= \bm{M}\bm{U}^n + \Delta t\gamma\bm{R}_{conv}^{n,1}+\gamma\Delta t\bm{B}_{bcs}^{n,2},\\
                \begin{split}
                    \left(\bm{M}-\gamma\Delta t\bm{K}\right)\bm{U}^{n+1} &= \bm{M}\bm{U}^n +\Delta t\gamma\bm{B}_{bcs}^{n+1} + \Delta t(1-\gamma)\left(\bm{K}\bm{U}^{n,2}+\bm{B}_{bcs}^{n,2}\right)\\
                    &\quad+\Delta t\delta\bm{R}_{conv}^{n,1}+\Delta t(1-\delta)\bm{R}_{conv}^{n,2}.
                \end{split}
            \end{align}
        \end{subequations}
        Note that we exploit the stiffly accurate property of the scheme (the last step equals the new solution) to avoid an explicit reconstruction step. 
        
        Second, we employ the third-order, four-stage DIRK scheme described in \cite{Calvo2001}.
        The coefficients are defined by the following Butcher tableau,
         \begin{equation}
            \label{eq:ldg:imex:tableau:dirk3}\renewcommand{\arraystretch}{1.3}
            \begin{array}{c|cccc|cccc}
            0        & 0 &   0       &    0     &  0   &      0      & 0       &0     & 0      \\
            \gamma   & 0 & \gamma    &    0     &  0   &\gamma   &       0  &  0      & 0   \\
            \frac{1+\gamma}{2}  & 0 & \frac{1-\gamma}{2}  & \gamma  & 0 &   \dfrac{1+\gamma}{2}-\iota_1 &   \iota_1    & 0   &0 \\
            1 & 0 & \zeta_1 & \zeta_2 & \gamma & 0 & 1-\iota_2 & \iota_2 & 0 \\
            \hline
                     & 0 & \zeta_1  & \zeta_2   & \gamma  & 0  & \zeta_1   & \zeta_2     & \gamma 
            \end{array}
        \end{equation}  
        with parameters $\gamma \approx 0.435866521508459$, $\zeta_1=-\frac{3}{2}\gamma^2+4\gamma+\frac{1}{4}$, $\zeta_2=\frac{3}{2}\gamma^2-5\gamma+\frac{5}{4}$, $\iota_1=-0.35$ and $\iota_2=\frac{\frac{1}{3}-2\gamma^2-2\zeta_2\iota_1\gamma}{\gamma(1-\gamma)}$. 
        
        Expanding the stages according to \eqref{eq:ldg:imex} yields the following sequential partial updates,
        \begin{subequations}
            \label{eq:ldg:imex:dirk3}
            \begin{align}
                \bm{U}^{n,1} &= \bm{U}^n, \\
                \left(\bm{M}-\gamma\Delta t\bm{K}\right)\bm{U}^{n,2} &= \bm{M}\bm{U}^n + \Delta t\gamma\bm{B}_{bcs}^{n,2}+\Delta t \gamma \bm{R}_{conv}^{n,1},\\
                \begin{split}
                \left(\bm{M}-\gamma\Delta t\bm{K}\right)\bm{U}^{n,3} &= \bm{M}\bm{U}^n + \Delta t\gamma\bm{B}_{bcs}^{n,3} + \Delta t\dfrac{1-\gamma}{2}\left(\bm{K}\bm{U}^{n,2}+\bm{B}_{bcs}^{n,2}\right)\\
                &\quad+\Delta t\left(\dfrac{1+\gamma}{2}-\iota_1\right)\bm{R}_{conv}^{n,1}+\Delta t\iota_1\bm{R}_{conv}^{n,2},
                \end{split} \\
                \begin{split}
                \left(\bm{M}-\gamma\Delta t\bm{K}\right)\bm{U}^{n,4} &= \bm{M}\bm{U}^n + \Delta t\gamma\bm{B}_{bcs}^{n,4}+\Delta t\zeta_1\left(\bm{K}\bm{U}^{n,2}+\bm{B}_{bcs}^{n,2}\right)\\
                &\quad+\Delta t\zeta_2\left(\bm{K}\bm{U}^{n,3}+\bm{B}_{bcs}^{n,3}\right)+\Delta t\left(1-\iota_2\right)\bm{R}_{conv}^{n,2}+\Delta t\,\iota_2\,\bm{R}_{conv}^{n,3},
                \end{split}
            \end{align}
        \end{subequations}
        with the time-marching update,
        \begin{equation*}
         \begin{split}
                \bm{U}^{n+1} &= \bm{U}^{n}+\Delta t\bm{M}^{-1}\biggl(
                \zeta_1\left(\bm{K}\bm{U}^{n,2}+\bm{B}_{bcs}^{n,2}+\bm{R}_{conv}^{n,2}\right)\\
                &\quad+\zeta_2\left(\bm{K}\bm{U}^{n,3}+\bm{B}_{bcs}^{n,3}+\bm{R}_{conv}^{n,3}\right)+
                \gamma\left(\bm{K}\bm{U}^{n,4}+\bm{B}_{bcs}^{n,4}+\bm{R}_{conv}^{n,4}\right)
                \biggr).
        \end{split}
        \end{equation*}

        These specific choices are motivated by the stability and accuracy properties established for LDG discretizations. 
        Under these IMEX time-marching schemes, the presented LDG formulation is unconditionally stable in the sense that the stable time step is upper-bounded by a positive constant independent of the spatial mesh size. 
        Furthermore, the fully discrete schemes exhibit an optimal convergence rate of order $\mathcal{O}\left(\Bar{h}^{\nu+1}+\Delta t^p\right)$, where $\Bar{h}$ denotes the element diameter, and $p$ is the temporal order of the IMEX scheme \cite{Wang2015,Shu2016}.

        From a practical point of view, the introduction of the IMEX time-marching schemes allows us to work with the time step restriction dictated by the standard Courant-Friedrich-Lewy (CFL) condition.
        In our case, we apply the constraint imposed by the explicit convective part unless stated otherwise.
        This directly translates to the CFL relationship,
        \begin{equation} \label{eq:ldg:imex:cfl}
            \Delta t\leq C\dfrac{\Bar{h}}{v_{\max}(2\nu+1)},\quad 0<C<1,
        \end{equation} 
        where $v_{\max}$ is the maximum characteristic speed, and we have taken into account the DG restriction arising from the distance between the interpolating nodes inside the element.

    \section{Sparse Grids Combination Technique} \label{ldg:combination}
    
    \joel{
        Although the proposed IMEX-LDG scheme is a highly suitable candidate for solving the aforementioned convection-diffusion problems, it presents several drawbacks related to its associated computational cost. 
        Specifically, the number of degrees of freedom becomes prohibitively large when utilizing highly refined meshes or working in higher dimensions, as both the number of elements and the nodes within each element suffer from the curse of dimensionality.
        This directly impacts the resulting sparse linear system.
        To overcome this challenge and to enhance the efficiency of high-order IMEX-LDG \joseanto{schemes}, we seek to couple such a methodology with sparse-grid implementations, specifically with the so-called combination technique, \cite{Griebel1990}.
    }
    
    The sparse-grid methodology includes a family of discretization techniques designed to mitigate the curse of dimensionality to some extent by reducing the number of degrees of freedom while retaining the prescribed accuracy.
    These approaches rely on two key ingredients: a hierarchical basis that provides an equivalent representation of the discrete function space to conventional nodal or modal bases; and a sparse tensor-product construction that systematically prunes anisotropic high-level refinements while preserving the dominant mixed-resolution interactions, \cite{Garcke2013}.

    \joel{
        Within this context, the authors of \cite{Bungartz2004} demonstrate that, when working with continuous FEM spaces build upon structured grids and polynomials of degree $\nu$, sparse-grid spaces reduce the number of degrees of freedom from $\mathcal{O}(\Bar{h}^{-d})$ to $\mathcal{O}(\Bar{h}^{-1}\,|\log_2 \Bar{h}|^{d-1})$ while achieving an interpolation error of the order $\mathcal{O}\left(\Bar{h}^{\nu+1}\,|\log_2\Bar{h}|^{d-1}\right)$ with respect to both the $L_{\infty}$- and the $L_{2}$-norm.
        In addition, efforts have been made to extend the sparse-grid methodology to the discontinuous Galerkin framework. 
        Several studies \cite{Cheng2016,Cheng2019,Tao2019,Guo2014} prove that for various types of PDEs, such as elliptic, transport, and systems of hyperbolic equations, operating within sparse-grid spaces reduces the degrees of freedom following the same asymptotic scaling as the continuous case. 
        This reduction is achieved while largely preserving the theoretical interpolation properties (typically bounded by $\mathcal{O}(\Bar{h}^{\nu+1}\,|\log_2 \Bar{h}|^d)$ or $\mathcal{O}(\Bar{h}^{\nu+1}|\log_2 \Bar{h}|^{d-1})$) in the $L^2$-norm, although reduced accuracy is often observed in the $L_{\infty}$-norm compared to the continuous case. 
    }

    In practice, the deployment of sparse grids requires specialized data structures and algorithms, particularly to support efficient matrix-vector products within Galerkin formulations.
    Representative developments in this direction include \cite{Cheng2016,Cheng2019}, where computations are performed directly in hierarchical spaces using basis functions that are mutually orthogonal across levels. 
    However, this approach imposes significant computational demands. 
    The overall performance of the method is critically dependent on the efficiency of hierarchical basis evaluation, basis transformation, and operator application during residual assembly.

    Among the various methodologies available, the presented work focuses on the combination technique, \cite{Griebel1990}.
    This approach constructs a sparse-grid representation by solving the PDE on a specific set of semi-coarsened, anisotropic grids of varying resolutions.
    The final solution is obtained by linearly combining these sub-solutions to approximate the result on a target full grid \cite{Garcke2013}.
    Theoretically, this method aims to achieve an accuracy comparable to that obtained by solving the PDE directly on a single, highly refined full mesh, or, equivalently, to achieve an order of convergence similar to \joseanto{that} obtained with the pure sparse-grid approach. 
    
    In this framework, the computational domain is partitioned into a family of component grids.
    Rather than solving on the prohibitively expensive single full grid, the PDE is solved independently on these selected semi-coarsened grids.
    Consequently, the total number of degrees of freedom, calculated as the sum over all component grids, is substantially smaller than that required for the equivalent full grid, thereby yielding significant reductions in computational cost.
    
    To formalize the combination technique, we first establish the necessary notation and standing assumptions.
    Let $\Omega_{\bm{\alpha}}$ denote the target full grid, characterized by a uniform refinement $\bm{\alpha} = (N_b 2^{N_L}, N_b 2^{N_L},\dots,N_b 2^{N_L})\in\N^d$.
    Here, $N_b>0$ represents the base number of elements in the coarsest mesh, and $N_L>0$ represents the maximum level of dyadic refinement.
    A semi-coarsened subgrid is denoted by $\Omega_{\bm{\alpha}}^{\bm{\kappa}}$, where the multi-index $\bm{\kappa}$ encodes the specific dyadic refinement level for each coordinate direction. This index must satisfy $0\leq\bm{\kappa}\leq N_L$.
    Consistent with the notation established in Section \ref{sec:ldg:num}, the structure of such a grid is defined via the tensor product of one-dimensional meshes,
    \begin{equation}
        \label{eq:ldg:ct:subgrid}
        \Omega_{\bm{\alpha}}^{\bm{\kappa}}([a,b]^d)=\bigotimes_{j=1}^d\Omega_{N_b 2^{\kappa_j}}([a,b]).
    \end{equation}

    The coarsest mesh within this hierarchy is termed the base grid (or base refinement).
    It corresponds to the zero-level index $\bm{\kappa}=0$ and is denoted by $\Omega_{\bm{\alpha}}^{\bm{0}}$. 
    We define the characteristic mesh size of this base grid as $H_b=(b-a)/N_b$.
    Consequently, for a general component grid $\Omega_{\bm{\alpha}}^{\bm{\kappa}}$, the grid spacing in the $k$-th coordinate direction is explicitly given by $h_k=H_b2^{-\kappa_k}$.

    Under these conventions, we consider the family of semi-coarsened sparse grids $\{\Omega_{\bm{\alpha}}^{\bm{\kappa}}\}_{0\leq\kappa\leq N_L}$.
    We apply the sparse-grid combination technique to the IMEX-LDG scheme described in the previous section. 
    Let $\bm{U}^{\bm{\kappa}}$ denote the numerical solution of Equation \eqref{eq:ldg:pde} computed via IMEX-LDG discretization on the specific component grid $\Omega_{\bm{\alpha}}^{\bm{\kappa}}$.
    
\begin{figure}
        \centering
        \includegraphics[width=0.5\linewidth]{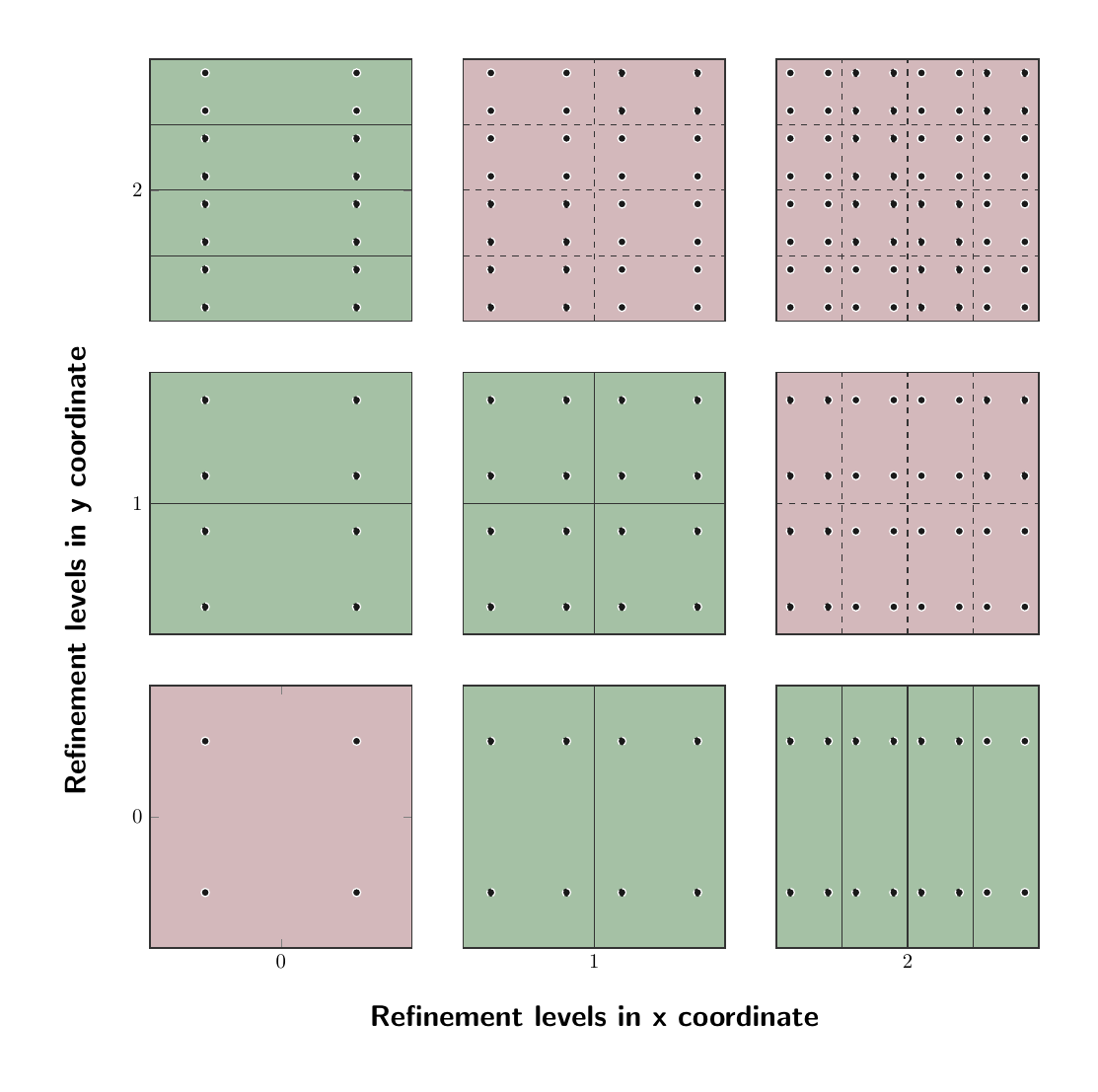}
        \caption[Sparse grid hierarchy and components]{Illustration of the sparse grid hierarchy for the case $\Omega_{4,4}$, with base refinement $N_b=1$ and maximum refinement level $N_L=2$. The grids marked with solid lines represent the subfamily of semi-coarsened meshes where the PDE is explicitly solved, whilst the grids marked with dashed lines represent the subgrids where the PDE is not solved. Furthermore, Gauss-Legendre nodes are depicted within each element to visualize the distribution of degrees of freedom.}
        \label{fig:ldg:ct:family-grids}
    \end{figure}

    Rather than solving the PDE on the single full grid $\Omega_{\bm{\alpha}}$, the problem is solved on a specific subfamily of grids $\{\Omega_{\bm{\alpha}}^{\bm{\kappa}}\}_{\bm{\kappa}\in\Lambda}$.
    The relevant index set $\Lambda$ is defined based on the $l_1$-norm of the refinement levels,
    \begin{equation}
        \label{eq:ldg:ct:index-set}
        \Lambda = \left\{0\leq\bm{\kappa}\leq N_L\,:\, |\bm{\kappa}|_1=\kappa_1+\cdots+\kappa_d=N_L-q,\,\, q=0,\dots,d-1\right\},
    \end{equation}
    which contains
    \begin{equation}
        \label{eq:ldg:ct:index-set:cardinal}
        \left|\Lambda\right| = \sum_{q=0}^{d-1}\binom{N_L-q+d-1}{d-1}
    \end{equation}
    component grids, \cite{Griebel1990}. 
    Crucially, the full grid $\Omega_{\bm{\alpha}}$ is not included in this subfamily, so the PDE is never solved at the finest resolution directly, yielding substantial savings in computational cost.
    This construction is illustrated in Figure \ref{fig:ldg:ct:family-grids}.

    Once the PDE has been solved on the entire set of subgrids $\{\Omega_{\bm{\alpha}}^{\bm{\kappa}}\}_{\bm{\kappa}\in\Lambda}$, the final fine-grid approximation $\bm{U}$ on $\Omega_{\bm{\alpha}}$ is obtained in two steps.
    First, each subgrid solution $\bm{U}^{\bm{\kappa}}$ is projected to the fine-grid space $\bm{V}_{\bm{\alpha}}^{\nu}$ via a prolongation operator.
    Second, the $\left|\Lambda\right|$ prolonged solutions are linearly combined to synthesize the final solution $\bm{U}$.

    \joel{
        The combination technique does not directly solve the PDE within the sparse-grid space. 
        Instead, by solving the problem on a defined set of subgrids, it yields a representation of the solution in the sparse space that, under suitable conditions, satisfies the same theoretical error bounds. 
        First, the condition of bounded mixed derivatives remains a prerequisite for ensuring optimal approximation properties.
        Furthermore, a critical requirement for these theoretical bounds to hold is the existence of an error splitting, \cite{bungartz1994}.
        That is, the functional error can be decomposable into a direct sum of one-dimensional error contributions along each coordinate direction, plus higher-order mixed error terms that are asymptotically negligible.
        While the underlying IMEX-LDG scheme possesses robust stability and consistency properties, and this error splitting is generally well established for linear PDEs, the primary challenge arises in nonlinear problems, where such an error splitting cannot always be guaranteed.
        Nevertheless, numerical experiments, \joseanto{such as those} presented in \cite{Griebel1990}, demonstrate that the combination technique performs reliably well provided the solution is sufficiently smooth.
    }

    \subsection{The prolongation operator} \label{ldg:ct:prolongation}
        A fundamental component of the combination technique within the presented IMEX-LDG framework is the transfer operator defined to map the semi-coarsened solutions, $\bm{U}^{\bm{\kappa}}$, onto the full grid space, thereby enabling the final recombination.
        In the literature, such transfer mappings are called prolongation operators.
        While  prolongation is a classical component of geometric multigrid methods, where nested mesh hierarchy and intergrid transfers are exploited to accelerate the solution of large linear systems, see e.g. \cite{Briggs2000,Oosterlee2000}, its role here is distinct. 
        In this context, the operator serves to lift the sparse-grid solutions to the finest-grid approximation space prior to their linear combination, ensuring that all contributions are represented within a unified function space.

        More specifically, let $\bm{V}_{\bm{\lambda}}^{\nu}$ denote the LDG vector space associated with the semi-coarsened grid $\Omega_{\bm{\alpha}}^{\bm{\kappa}}$, with $\lambda_k=N_b 2^{\kappa_k}$ for $\kappa_k<N_L,\,k=1,\dots,d$, and let $\bm{V}_{\alpha}^{\nu}$ denote the LDG space on the target full grid.
        The prolongation operator $\Pi^{\bm{\kappa}}:\bm{V}^{\nu}_{\bm{\lambda}}\to\bm{V}_{\bm{\alpha}}^{\nu}$ is defined as the natural injection induced by mesh refinement. 
        It can be characterized as the $L^2$-orthogonal projection onto the fine space $\bm{V}_{\bm{\alpha}}^{\nu}$: find $\Pi^{\bm{\kappa}}v_{\bm{\lambda}}\in\bm{V}_{\bm{\alpha}}^{\nu}$ such that
        \begin{equation}
            \label{eq:ldg:ct:prolongation}
            \int_{\Omega} \left(\Pi^{\bm{\kappa}}v_{\bm{\lambda}}\right)w_{\bm{\alpha}} = \int_{\Omega}v_{\bm{\lambda}}w_{\bm{\alpha}},\quad\forall  w_{\bm{\alpha}}\in\bm{V}^{\nu}_{\bm{\alpha}}.
        \end{equation}
        Since the semi-coarsened partitions are nested within the full refinement by construction, $\bm{V}_{\bm{\lambda}}^{\nu}\subset\bm{V}_{\bm{\alpha}}^{\nu}$, the function $v_{\bm{\lambda}}$ is already exactly representable in the fine space.
        Consequently, $\Pi^{\bm{\kappa}}$ acts as the natural inclusion, meaning the projection introduces no approximation error, \cite{Cheng2019}. 

        In a basis representation, the prolongation matrix is constructed by applying the operator $\Pi^{\bm{\kappa}}$ to each basis function of the coarse space $\bm{V}_{\bm{\lambda}}^{\nu}$ and expanding the resulting image in the basis of the fine space $\bm{V}^{\nu}_{\bm{\alpha}}$.
        Because the LDG polynomial basis is built as a tensor product of one-dimensional bases, the multidimensional transfer operator is assembled from the tensor product of one-dimensional operators.
            
        To derive the explicit matrix form, it suffices to define the action element-wise.
        Consider a fine element $Q_{\bm{\beta}}\in\Omega_{\bm{\alpha}}$, where $\bm{\beta}\leq\bm{\alpha}$, that is geometrically contained within a coarse element $Q_{\bm{o}}$ of the sparse grid $\Omega_{\bm{\alpha}}^{\bm{\kappa}}$, which we denote by $\Omega_{\bm{\lambda}}$, where $\bm{o}\leq\bm{\lambda}$. 
        By construction, each coarse cell $Q_{\bm{o}}$ decomposes into a disjoint union of fine cells.
        Specifically, the number of fine cells contained within a single coarse cell is given by $2^{\sum_{e=1}^d(N_L-\kappa_e)}$, which corresponds to the dyadic refinement ratio accumulated across all spatial directions. 

        Let $\{\phi_{\bm{i}}^{\bm{\beta}}\}_{\bm{i}\leq\nu+1}$ and  $\{\phi_{\bm{j}}^{\bm{o}}\}_{\bm{j}\leq\nu+1}$ denote the local polynomial bases on the fine element $Q_{\bm{\beta}}$ and the coarse element $Q_{\bm{o}}$, respectively. 
        We define the local prolongation operator $\Pi_{\beta}^{\bm{o}}:\bm{V}^{\nu}_{\bm{\lambda}}|_{Q_{\bm{o}}}\to\bm{V}^{\nu}_{\bm{\alpha}}|_{Q_{\bm{\beta}}}$ as the $L^2$-projection restricted to the domain of the fine element $Q_{\bm{\beta}}$.
        Thus, for any test $w_{\bm{\beta}}=\phi_{\bm{i}}^{\bm{\beta}}\in\bm{V}_{\bm{\alpha}}^{\nu}|_{Q_{\bm{\beta}}}$, the projection satisfies the variational equality
        \begin{equation}
            \label{eq:ldg:ct:local_projection:variational}
            \int_{Q_{\bm{\beta}}} \left(\Pi_{\bm{\beta}}^{\bm{o}}v_{\bm{o}}\right)\phi_{\bm{i}}^{\bm{\beta}} = \int_{Q_{\bm{\beta}}}v_{\bm{o}}\phi_{\bm{i}}^{\bm{\beta}}.
        \end{equation}
        
        On the left-hand side, expanding $\Pi_{\bm{\beta}}^{\bm{o}}v_{\bm{o}}$ in the basis on $Q_{\bm{\beta}}$ yields the local fine-grid mass matrix $\bm{M}^{\bm{\beta}}$ acting on the unknown coefficient vector, which we denote by $\bm{v}^{\bm{\beta}}$, whose coefficients are represented by $[\bm{v}^{\bm{\beta}}]_{\bm{i}}$. 
        
        On the right-hand side, we expand the coarse function on $Q_{\bm{o}}$ as $v_{\bm{o}}=\sum_{\bm{j}\leq\nu+1}v_{\bm{j}}^{\bm{o}}\phi_{\bm{j}}^{\bm{o}}$. 
        Thus, restricting the integration to the fine element $Q_{\bm{\beta}}\subset Q_{\bm{o}}$, we obtain
        \begin{equation}
            \label{eq:ldg:ct:aux-prolong:1}
             \int_{Q_{\bm{\beta}}}v_{\bm{o}}\phi_{\bm{i}}^{\bm{\beta}} = \sum_{\bm{j}\leq\nu+1} v_{\bm{j}}^{\bm{o}}\int_{Q_{\bm{\beta}}}\phi_{\bm{j}}^{\bm{o}}\phi_{\bm{i}}^{\bm{\beta}}. 
        \end{equation}
        Since the support of each coarse basis function $\phi_{\bm{j}}^{\bm{o}}$ extends over the entire coarse element $Q_{\bm{o}}$, while the integration domain is restricted to the fine element $Q_{\bm{\beta}}$, these coupling integrals involve evaluating the coarse basis on a restricted domain.
        Consequently, the local orthogonality properties of the basis do not hold locally on $Q_{\bm{\beta}}$, rendering the resulting interaction matrix generally non-diagonal.

        Let $\bm{v}^{\bm{o}}$ be the local vector of coarse coefficients.
        We define the fine-coarse coupling matrix $\bm{\Xi}_{\bm{\beta}}^{\bm{o}}\in\mathcal{M}_{N_{\nu}}(\R)$ with entries given by 
        \begin{equation} \label{eq:ldg:ct:matrix:fine-coarse-coupling}
        [\bm{\Xi}_{\bm{\beta}}^{\bm{o}}]_{\bm{i}\bm{j}}=\int_{Q_{\bm{\beta}}}\phi_{\bm{j}}^{\bm{o}}\phi_{\bm{i}}^{\bm{\beta}}.
        \end{equation} 
        The local prolongation relationship between the coarse element $Q_{\bm{o}}\in\Omega_{\bm{\lambda}}$ and the fine element $Q_{\bm{\beta}}\in\Omega_{\bm{\alpha}}$ is governed by the linear system,
        \begin{equation}
            \label{eq:ldg:ct:local_prolongation:eq}
            \bm{M}^{\bm{\beta}}\bm{v}^{\bm{\beta}} = \bm{\Xi}_{\bm{\beta}}^{\bm{o}}\bm{v}^{\bm{o}},\quad\implies\quad\bm{v}^{\bm{\beta}} = \left(\bm{M}^{\bm{\beta}}\right)^{-1}\bm{\Xi}_{\bm{\beta}}^{\bm{o}}\bm{v}^{\bm{o}}.
        \end{equation}
        This yields the explicit definition of the local prolongation matrix $\bm{\Pi}_{\bm{\beta}}^{\bm{o}}$, which maps coefficients from the coarse to the fine,
        \begin{equation}
            \label{eq:ldg:ct:local_prolongation:matrix}
            \bm{\Pi}_{\bm{\beta}}^{\bm{o}}=\left(\bm{M}^{\bm{\beta}}\right)^{-1}\bm{\Xi}_{\bm{\beta}}^{\bm{o}}.
        \end{equation}
        
        The global assembly proceeds by placing each local block $\bm{\Pi}^{\bm{o}}_{\bm{\beta}}$ into the global matrix structure. 
        Rows are indexed by the fine element $Q_{\bm{\beta}},\,\bm{\beta}\leq\bm{\alpha}$ and columns by the coarse cell $Q_{\bm{o}},\,\bm{o}\leq\bm{\lambda}$, strictly  following to the ordering convention established in Section \ref{ldg:num:basis}.
        This process yields a rectangular global matrix $\bm{\Pi}^{\bm{\kappa}}\in\mathcal{M}_{N_{\bm{\alpha}}N_{\nu}, N_{\lambda}\times N_{\nu}}(\R)$ with more rows than columns.
        This dimensionality reflects the one-to-many mapping inherent in the mesh refinement, where a single coarse element distributes its information to $2^{\sum_{k=1}^d(N_L-\kappa_k)}$ distinct fine cells.

        \begin{remark}
            When employing spectral nodal bases defined on Gauss-Legendre nodes, computation of the fine-coarse coupling matrix entries \eqref{eq:ldg:ct:matrix:fine-coarse-coupling} becomes highly efficient. 
            By using the Gaussian quadrature rule associated with the fine element nodes, the fine-grid basis functions satisfy $\phi_{\bm{i}}^{\bm{\beta}}(\xi_{\bm{j}})=\delta_{\bm{i}}^{\bm{j}}$. 
            Consequently, the integral reduces to a simple weighted evaluation of the coarse basis function $\phi_{\mathbf{j}}^{\mathbf{o}}$ at the mapped fine-grid quadrature nodes, performed dimension-by-dimension.
        \end{remark}
            
        \begin{remark}
            Although the global prolongation matrix mapping from each sparse subgrid to the finest grid is well-defined, its explicit assembly is rarely practical. 
            Storing transfer matrices for every sparse-grid component would require memory resources proportional to the full grid, reducing the storage benefits inherent to the sparse-grid methodology. 
            In practice, we rely on the element-wise prolongation operator $\bm{\Pi}_{\bm{\beta}}^{\bm{o}}=\left(\bm{M}^{\bm{\beta}}\right)^{-1}\bm{\Xi}_{\bm{\beta}}^{\bm{o}}$.
            This local operator enables a matrix-free implementation, where prolongation is computed on the fly during the combination process.
        \end{remark}

        \subsection{Recombination and practical considerations}
            With the prolongation operator established, we can now state the final recombination formula within the extrapolation-based sparse grid framework.
            As detailed in \cite{Bungartz2004,Griebel1990}, a specific linear combination of the prolonged subgrid solutions yields an approximation equivalent to that of the sparse-grid space.
            
            This methodology is based on the seminal work of Smolyak, who introduced the hierarchical combination of anisotropic tensor-product contributions to approximate multivariate operators and functionals in an abstract setting \cite{Smolyak1962}.
            
            Formally, let $\bm{U}^{\bm{\kappa}}$ denote the IMEX-LDG solution computed on the semi-coarsened grid $\Omega_{\bm{\alpha}}^{\bm{\kappa}}$, and let $\bm{\Pi}^{\bm{\kappa}}$ represent the global prolongation operator mapping this solution to the finest grid space. 
            The sparse-grid combination technique formula is then explicitly given by
            \begin{equation}
                \label{eq:ldg:ct:recombination}
                \bm{U} = \sum_{q=0}^{d-1}(-1)^q\binom{d-1}{q}\sum_{|\bm{\kappa}|_1=N_L-q}\bm{\Pi}^{\bm{\kappa}}\bm{U}^{\bm{\kappa}}.
            \end{equation}

            While the full recombination expression \eqref{eq:ldg:ct:recombination} is valuable for theoretical analysis and convergence verification, it is not the most practical form for large-scale computations. 
            Solving the PDE \eqref{eq:ldg:pde} on the subfamily indexed by \eqref{eq:ldg:ct:index-set} reduces simulation time for two primary reasons. 
            The total number of degrees of freedom is significantly lower than that of the full grid, and coarser refinements in the component grids typically allow for larger time steps.
            However, explicitly projecting all subgrid solutions onto the full grid to assemble $\bm{U}$ would negate memory savings.
            To avoid this bottleneck, the recombination formula \eqref{eq:ldg:ct:recombination} is interpreted as a sparse functional representation.
            Instead of assembling a global vector, we store the solutions as a collection of coefficient vectors for each component grid.
            The final approximation is then evaluated on demand at specific query points by summing the prolonged contributions locally, without ever materializing the full-grid coefficient vector.

            We define the sparse-grid interpolant $u_{\bm{\alpha}}^{SG}$ as the final approximation to the solution of \eqref{eq:ldg:pde} evaluated at an arbitrary spatial point $\bm{x}^*\in\Omega$. 
            The evaluation proceeds locally and on demand via the following steps:
            \begin{itemize}
                \item Identify the unique fine element $Q_{\bm{\beta}}\subset\Omega_{\bm{\alpha}}$ that contains the query point $\xv^*$.
                \item For each component grid $\Omega_{\bm{\alpha}}^{\bm{\kappa}}$ with $\kappa\in\Lambda$, determine the unique coarse element $Q_{\bm{o}(\bm{\kappa})}\subset\Omega_{\bm{\alpha}}^{\bm{\kappa}}$ such that $\xv^*\in Q_{\bm{o}(\bm{\kappa})}$.
                \item Retrieve the coarse coefficient vector $\bm{v}^{\bm{o}}$ associated with the solution on that grid.
                Apply the local prolongation matrix $\bm{\Pi}^{\bm{o}}_{\bm{\beta}}$ to obtain the equivalent fine-element coefficients $\bm{v}^{\bm{\beta}}=\bm{\Pi}_{\bm{\beta}}^{\bm{o}}\bm{v}^{\bm{o}}$. 
                \item Evaluate the fine-element basis functions at $\xv^*$ and compute the inner product with $\bm{v}^{\bm{\beta}}$, weighted by the corresponding combination coefficient.
            \end{itemize}
            
\begin{figure}
                \centering
                \includegraphics[width=0.75\linewidth]{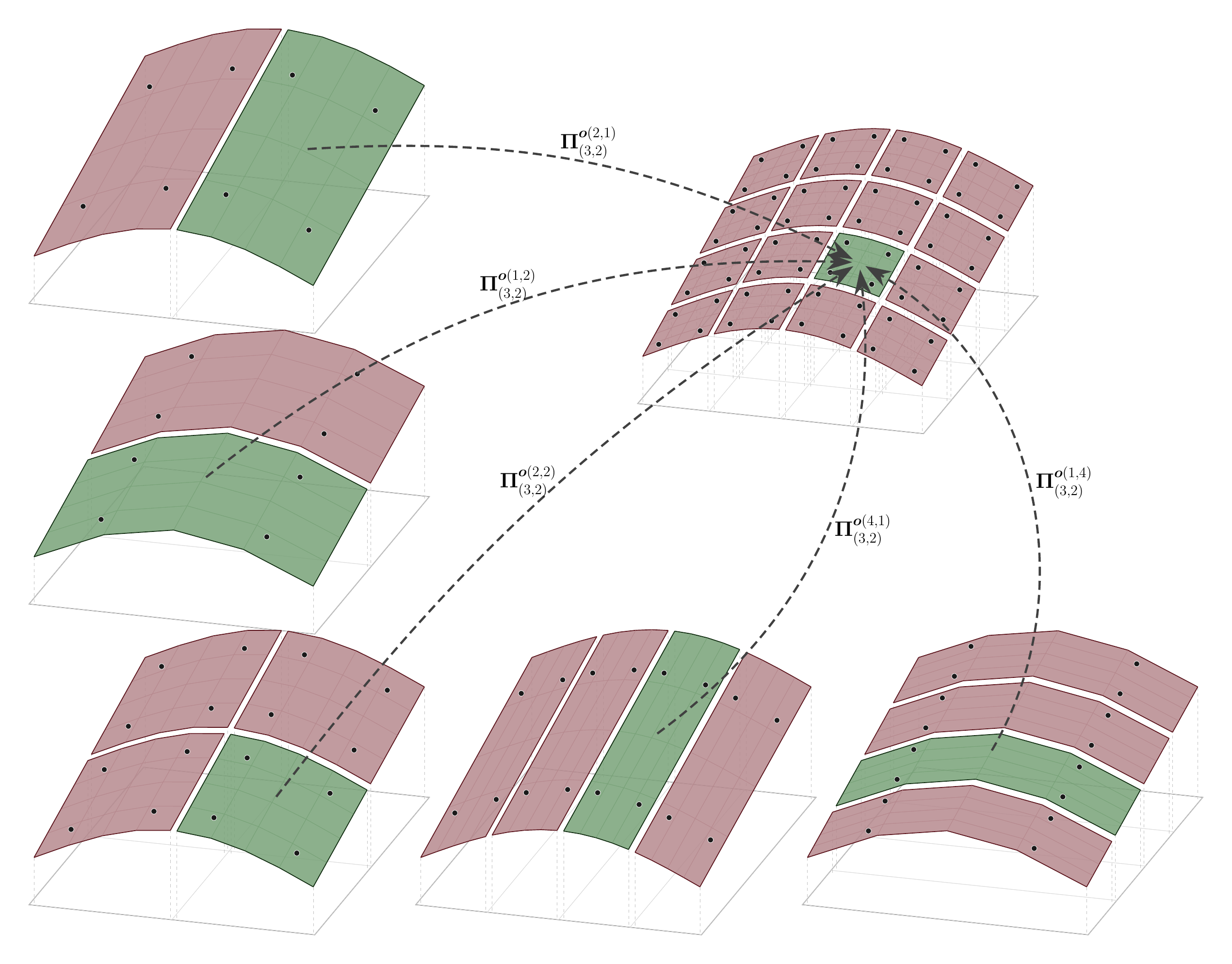}
                \caption{Schematic representation of the evaluation process using the combination technique, highlighting the projection from the component grids defined in Figure \ref{fig:ldg:ct:family-grids}.}
                \label{fig:ldg:ct:projection-ct}
            \end{figure}

            Mathematically, the point-wise value is assembled as 
            \begin{equation}
                \label{eq:ldg:ct:recombination:local}
                u_{\bm{\alpha}}^{SG}(\xv^*)=\sum_{q=0}^{d-1}(-1)^q\binom{d-1}{q}\sum_{|\bm{\kappa}|_1=N_L-q}\left(\sum_{\bm{i}\leq\nu+1}\left[\bm{\Pi}_{\bm{\beta}}^{\bm{o}(\bm{\kappa})}\bm{v}^{\bm{o}(\bm{\kappa})}\right]_{\bm{i}}\phi_{\bm{i}}^{\bm{\beta}}(\xv^*)\right), 
            \end{equation}
            where $Q_{\bm{o}(\bm{\kappa})}$ denotes the unique coarse element from grid $\bm{\kappa}$ containing $\xv^*$.
            A schematic representation of this projection and recombination is provided in Figure \ref{fig:ldg:ct:projection-ct}. 

\section{Numerical experiments} \label{ldg:experiments}

    \joseanto{
    This section evaluates the performance, accuracy, and computational efficiency of the proposed high-order IMEX-LDG scheme and its sparse-grid combination technique variant. The numerical experiments are organized into two main blocks, followed by two application-oriented tests. The first block, presented in Sections~\ref{ldg:experiments:linear} and~\ref{ldg:experiments:vburgers}, assesses the methodologies using classical multidimensional convection--diffusion benchmarks: a linear Gaussian pulse governed by a constant-coefficient convection--diffusion equation, and a nonlinear viscous Burgers' equation with constant viscosity and an adjusted source term. These examples are designed to thoroughly test convergence rates and scalability, while also exposing the inherent limitations of sparse-grid extrapolation when dealing with non-aligned nonlinear fluxes and high-frequency aliasing effects. Regarding the second block, Section~\ref{ldg:experiments:ac}, focuses on the advective Allen-Cahn equation to demonstrate the robustness of both schemes in handling highly stiff reaction sources. In this context, preserving the geometric integrity and physical metrics of moving phase boundaries is crucial, providing a rigorous test for the combination technique's ability to filter out numerical dissipation. Finally, Section~\ref{ldg:experiments:bs} applies the methodology to a multidimensional Black-Scholes model under credit risk adjustments, illustrating its behavior in a practical high-dimensional financial setting.
    }
    
    To quantify the accuracy of our schemes, all numerical errors are computed with respect to an exact reference solution, denoted by $u^*$, in the norms of the associated functional spaces. Specifically, we utilize the $L^p$-norm, $p \in \mathbb{N}$, 
    \begin{equation}
                \label{eq:ldg:experiments:lp-error}
                ||e||_p = \left(\int_{\Omega} \left|u^*(T,\xv)-u(T,\xv)\right|^{p}\text{d}\xv\right)^{1/p}.
            \end{equation}
            In practice, we apply the Gauss-Legendre quadrature rule to approximate the integral. 
            Consequently, the computation relies on the solution values evaluated specifically at the Gauss-Legendre nodes.
            Additionally, we consider the $L^{\infty}$-norm, defined formally by the limit $p\to\infty$ as
            \begin{equation}
                \label{eq:ldg:experiments:linfty-error}
                ||e||_{\infty}=\sup_{\xv\in\Omega}\left|u^*(T,\xv)-u(T,\xv)\right|.
            \end{equation}
            
            Correspondingly, the experimental order of convergence is the $L^p$-norm calculated using the errors obtained from two consecutive mesh refinements, characterized by element diameters $\bar{h}$ and $\bar{h}/2$, respectively, i.e.,
            \begin{equation}
                \label{eq:ldg:experiments:order}
                \text{Order}_p = \log_2\left(\dfrac{||e(\bar{h})||_p}{||e(\bar{h}/2)||_p}\right). 
            \end{equation}
    
    A fundamental computational challenge in high-order, multidimensional settings is the massive memory footprint and computational cost associated with solving the linear systems at each IMEX-RK sub-step. To mitigate these memory bottlenecks, we primarily employ the Generalized Minimum Residual (GMRES) iterative solver. GMRES facilitates a matrix-free framework by strictly requiring matrix-vector products, thereby circumventing the need to explicitly assemble the global system matrix. Instead, the solver iteratively evaluates the explicit diffusive operator on a given solution vector. While this matrix-free approach precludes the use of matrix-dependent algebraic preconditioners and increases the floating-point operations per iteration, it yields critical memory savings that make highly refined and high-dimensional simulations---such as those utilizing quadratic polynomials---feasible. For lower-dimensional or highly coarsened scenarios where memory permits, we revert to explicit matrix assembly to accelerate convergence.
    
    Finally, both the full-grid and sparse-grid methodologies share a unified, dimension-agnostic core implementation written in C++ using template metaprogramming. This architecture guarantees that comparisons across different spatial dimensions and numerical tests remain consistent. To handle the immense computational workloads, the software is fully parallelized using the Message Passing Interface (MPI) standard. All simulations were executed on the Picasso supercomputer cluster at the University of M\'alaga. Although the heterogeneous nature of this computing environment may introduce minor variations in raw execution times, the reported performance metrics remain highly robust for demonstrating the relative efficiency gains achieved by the evaluated methodologies.

        \subsection{Linear problem} \label{ldg:experiments:linear}
        The first benchmark addresses the classical linear problem with constant coefficients, governed by 
        \begin{equation}
            \label{eq:ldg:experiments:linear-pde}
            \partial_t u + c_v\sum_{k=1}^d\partial_{x_k}u = a\sum_{k=1}^d\partial_{x_kx_k}u. 
        \end{equation}
        This equation aligns with the general conservative formulation \eqref{eq:ldg:pde} by setting the flux vector to $\bm{f}(u)=\left[c_vu,\dots,c_vu\right]^T$, the source term to zero, and the diffusion matrix to $A=a\bm{I}_d$, where $\bm{I}_d$ denotes the identity matrix.
        In this context, we present results concerning the well-known $d$-dimensional Gaussian pulse problem, which is defined by taking a point mass $m$ as an initial condition.
        This is mathematically modeled as a Dirac delta function, $\delta(\xv-\xv_0)$, where $\xv_0$ is \joseanto{the position}of the mass. 
        Hereafter, we assume that the mass is centered at $\xv_0 = \bm{0}$.
        The analytical solution of the problem is given by:
        \begin{equation} 
            \label{eq:ldg:experiments:linear-pde:sol}
            u^*(t,\xv) = \dfrac{m}{(4\pi at)^{d/2}}\exp\left\{-\dfrac{\sum_{k=1}^d\left(x_k-c_v t\right)^2}{4at}\right\},\quad t>0.
        \end{equation}

        First, we study the experimental order of convergence and the errors incurred by both the full-grid IMEX-LDG scheme and the sparse-grid combination technique in the $2$- and $3$-dimensional cases, utilizing polynomial degrees $\nu=1$ and $\nu=2$. 
        To start the simulation from a smooth initial condition, we use a time-advanced state of the exact solution.
        That is, for a fixed time $t^*>0$, we choose $u_0(\xv) = u^*(t^*, \xv)$, which yields the solution \eqref{eq:ldg:experiments:linear-pde:sol} with the time shift $t \to t^*+t$.

        Due to the exponential decay of the solution and its symmetry, we impose periodic boundary conditions on a sufficiently large domain. 
        This allows us to isolate the performance of the numerical scheme from boundary effects.
        For these experiments, we set a unit mass $m=1$ and define the initial condition as the exact solution at time $t^*=0.4$.
        The model parameters are chosen as follows: wave speed $c_v=1$, diffusion coefficient $a=0.05$, and final simulation time $T=0.5$ (corresponding to the physical time $t=0.9$).
        Under this setting, the spatial domain is defined as $\Omega=(-2.35, 3.65)^d$.
        Finally, the time step is selected in accordance with the CFL condition \eqref{eq:ldg:imex:cfl} by setting the constant to $C=0.9$.

        To study the order of convergence of our schemes in the $2$D case, we design an experimental setup that allows for a direct error comparison between the full-grid and combination technique approaches at equivalent effective resolutions. 
        To this end, we establish a set of reference full-grid resolutions, defined by $\bm{\alpha}=\{80, 160, 320, 640\}$.
        Based on these, we conduct the combination technique experiments using two distinct strategies.

        \begin{table}[ht]
            \centering
            \caption{Convergence analysis for the $2$D linear convection-diffusion problem: Full Grid vs. Combination Technique}
            \setlength{\tabcolsep}{6pt} 
            \begin{tabular}{ll S[table-format=1.2e-1,scientific-notation=true] S[table-format=1.2] S[table-format=1.2e-1, scientific-notation=true] S[table-format=1.2] S[table-format=1.2e-1,scientific-notation=true]}
                \toprule
                Case & Config. & {$||e||_1$} & {EOC} & {$||e||_2$} & {EOC} & {Time (ms)} \\
                \midrule
                \multicolumn{7}{l}{\textbf{Polynomial degree $\nu=1$}}\\
                \textit{Full Grid} & $\bm{\alpha}=80$  & 4.69e-3 & {--} & 4.04e-3 & {--} & 3.43e2 \\
                                   & $\bm{\alpha}=160$ & 1.13e-3 & 2.05 & 9.69e-4 & 2.15 & 1.57e3 \\
                                   & $\bm{\alpha}=320$ & 2.79e-4 & 2.02 & 2.38e-4 & 2.02 & 1.07e4 \\
                                   & $\bm{\alpha}=640$ & 6.93e-5 & 2.01 & 2.01e-4 & 2.05 & 1.21e5 \\
                \noalign{\smallskip}
                \textit{CT}        & $N_B=10$ & 1.07e-2 & {--} & 8.40e-3 & {--} & 4.24e2 \\
                Max. level $N_L=3$ & $N_B=20$ & 1.22e-3 & 3.14 & 1.06e-3 & 2.99 & 1.10e3 \\
                                   & $N_B=40$ & 2.79e-4 & 2.13 & 2.40e-4 & 2.14 & 4.66e3 \\
                                   & $N_B=80$ & 6.93e-5 & 2.01 & 5.92e-5 & 2.02 & 4.90e4 \\
                \noalign{\smallskip}
                \textit{CT}        & $N_L=3$  & 1.07e-2 & {--} & 8.40e-3 & {--} & 2.95e2 \\
                Base ref. $N_B=10$ & $N_L=4$  & 2.73e-3 & 1.98 & 2.16e-3 & 1.96 & 9.60e2 \\
                                   & $N_L=5$  & 6.85e-4 & 2.00 & 5.50e-4 & 1.97 & 2.94e3 \\
                                   & $N_L=6$  & 1.73e-4 & 1.99 & 1.41e-4 & 1.96 & 7.92e3 \\
                \midrule
                \multicolumn{7}{l}{\textbf{Polynomial degree $\nu=2$}}\\
                \textit{Full Grid} & $\bm{\alpha}=80$  & 1.22e-4 & {--} & 1.06e-4 & {--} & 1.91e3 \\
                                   & $\bm{\alpha}=160$ & 1.53e-5 & 3.00 & 1.33e-5 & 3.00 & 1.76e4 \\
                                   & $\bm{\alpha}=320$ & 1.91e-6 & 3.00 & 1.66e-6 & 3.00 & 1.93e5 \\
                                   & $\bm{\alpha}=640$ & 2.39e-7 & 3.00 & 2.07e-7 & 3.00 & 2.07e6 \\
                \noalign{\smallskip}
                \textit{CT}        & $N_B=10$ & 1.51e-4 & {--} & 1.25e-4 & {--} & 1.68e3 \\
                Max. level $N_L=3$ & $N_B=20$ & 1.53e-5 & 3.30 & 1.33e-5 & 3.23 & 8.52e3 \\
                                   & $N_B=40$ & 1.91e-6 & 3.01 & 1.66e-6 & 3.01 & 6.79e4 \\
                                   & $N_B=80$ & 2.39e-7 & 3.00 & 2.07e-7 & 3.00 & 7.21e5 \\
                \noalign{\smallskip}
                \textit{CT}        & $N_L=3$  & 1.51e-4 & {--} & 1.25e-4 & {--} & 1.71e3 \\
                Base ref. $N_B=10$ & $N_L=4$  & 1.97e-5 & 2.94 & 1.65e-5 & 2.92 & 5.74e3 \\
                                   & $N_L=5$  & 2.55e-6 & 2.94 & 2.16e-6 & 2.93 & 2.20e4 \\
                                   & $N_L=6$  & 3.31e-7 & 2.95 & 2.82e-7 & 2.94 & 8.70e4 \\
                \bottomrule
            \end{tabular}
            \label{tab:ldg:experiments:linear-pde:2d}
        \end{table}
 
        The first strategy \joseanto{consists of} fixing the maximum refinement level at $N_L=3$, while successively doubling the base refinement from $N_b=10$ to $N_b=80$. 
        Under this approach, the convergence rate should closely match that of the full-grid case, as we are explicitly enforcing a reduction in the element diameter. 
        This strategy serves to verify the proper behavior of the interpolation procedure.
        
        Alternatively, in the second strategy, we fix the base refinement at $N_b=10$ and perform simulations across varying maximum refinement levels, $N_L=\{3,4,5,6\}$, which correspond to the effective resolutions of the full-grid sequence. 
        This scenario is where the combination technique is expected to prove most advantageous. 
        Because we are analyzing a linear problem, the resulting convergence order should closely align with the theoretical estimates established in Section \ref{ldg:combination}.
        
        Regarding the technical configuration of the linear solver, for both the full-grid and the combination technique approaches, we explicitly assemble the global matrix system, as the memory requirements remain computationally manageable. 
        This also allows us to employ an Algebraic Multigrid (AMG) preconditioner \cite{Yang2006}, which is widely recognized as a gold standard for solving elliptic problems.
        Finally, to ensure the fairest possible comparison of computational times, the workload for all tested meshes was distributed across $50$ MPI processes, even though this may not represent the optimal parallelization strategy for every individual resolution.
        
        Table \ref{tab:ldg:experiments:linear-pde:2d} presents the error, convergence rates, and simulation times required for each considered scenario.
        It is worth noting that the errors are computed using the $L^1$ and $L^2$ functional space norms.
        To facilitate this, the combination technique solution is projected onto the effective fine-grid resolution.

        For both full-grid cases, the expected convergence rates are successfully achieved, demonstrating the capabilities of the method. 
        In addition, we compute the $L^{\infty}$-norm using \eqref{eq:ldg:experiments:linfty-error}. 
        Both polynomial configurations keep the expected convergence rates, achieving maximum absolute errors at the finest resolution of $2.01\times 10^{-4}$ and $4.80\times 10^{-7}$ for $\nu=1$ and $\nu=2$, respectively.
        As expected, the convergence rate for the scenarios with a fixed maximum level and varying base refinement \joseanto{matches} the theoretical convergence rate of the full-grid schemes.
        Furthermore, they converge at the expected rate in the $L^{\infty}$-norm, exhibiting maximum absolute errors in the finest configurations of $2.03\times 10^{-4}$ and $4.79\times 10^{-7}$.
        These represent negligible differences compared to the maximum absolute errors observed in the full-grid case.

        \begin{figure}[htbp]
            \centering
            \begin{subfigure}[b]{0.49\textwidth}
                \centering
                \includegraphics[width=\linewidth]{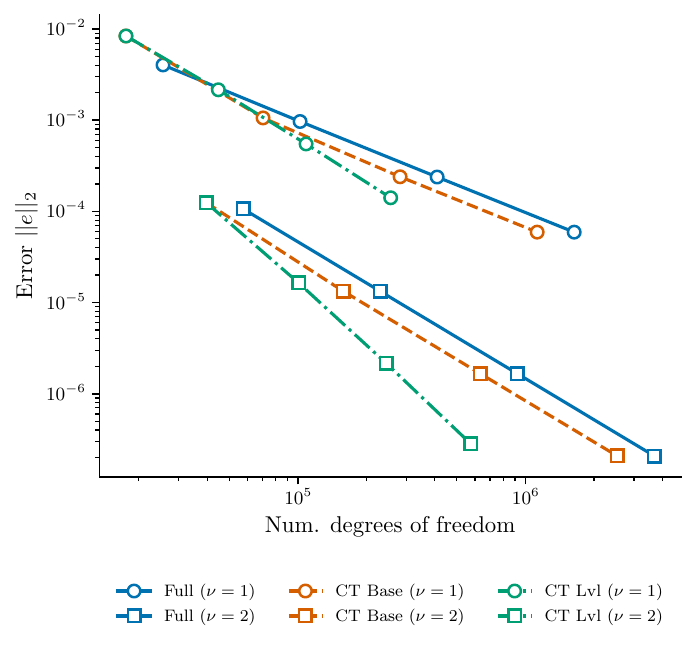}
                \caption{Efficiency in terms of DoFs.}
                \label{fig:ldg:experiments:linear-pde:2d:efficiency:DoF}
            \end{subfigure}
            \hfill 
            \begin{subfigure}[b]{0.49\textwidth}
            \centering
                \includegraphics[width=\linewidth]{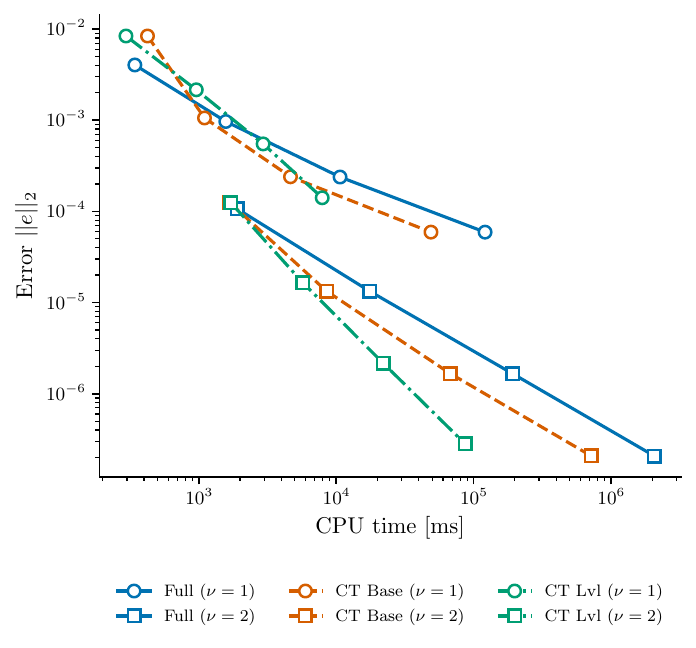}
            \caption{Efficiency in terms of CPU time.}
            \label{fig:ldg:experiments:linear-pde:2d:efficiency:Time}
            \end{subfigure}
        \caption{Efficiency comparison between the full-grid IMEX-LDG scheme and the sparse-grid combination technique for the 2D linear convection-diffusion problem. Data derived from Table \ref{tab:ldg:experiments:linear-pde:2d}.}
        \label{fig:ldg:experiments:linear-pde:2d:efficiency}
        \end{figure}

        Conversely, for the cases with a fixed base refinement and varying maximum levels, we observe \joseanto{a slight} decrease in the convergence rates for both the $L^1$ and $L^2$ norms predicted by the sparse-grid theory discussed above.
        We also confirm that the convergence rate in the $L^{\infty}$ norm fails to keep this theoretical scaling, dropping to $1.48$ for $\nu=1$ and $2.22$ for $\nu=2$.
        The finest levels exhibit maximum absolute errors of $1.48\times 10^{-3}$ and $3.40\times 10^{-6}$, which are nearly an order of magnitude larger than those previously presented.

        However, the situation is entirely different when analyzing the efficiency achieved across the three scenarios.
        Figure \ref{fig:ldg:experiments:linear-pde:2d:efficiency} illustrates the efficiency of these methods by evaluating both the total number of degrees of freedom and the total computational time against the resulting $L^2$-norm error.
        In these plots, optimal efficiency is represented by curves positioned closer to the bottom-left corner, indicating configurations that achieve lower errors for a given computational cost.
        
        From both plots, we first observe that all configurations utilizing a polynomial basis of degree $\nu=2$ significantly outperform those with $\nu=1$, highlighting the inherent advantage of high-order methods. 
        Second, because these \joseanto{third-order} methods ($\nu=2$) demand a greater number of degrees of freedom and computational time, the superior performance of the pure combination technique (where the maximum level is varied) becomes much more pronounced.
        Overall, the combination technique with varying base refinement consistently outperforms the full-grid scheme; however, the combination technique with varying maximum levels yields the highest efficiency among all tested approaches.
        This performance hierarchy is less distinct for the $\nu=1$ cases, as the memory and operational requirements for a $2$D domain remain sufficiently low that the relative saving provided by the sparse-grid approach \joseanto{is less} pronounced.

        The experimental setup for the $3$D cases closely follows that of the $2$D scenarios, with a few critical modifications required to address \joseanto{computational constraints}.
        The most notable difference is that it is no longer feasible to compute the full-grid scheme for the finest resolution considered before ($\bm{\alpha}=640$) within a reasonable time-frame or resource budget, particularly for polynomials of degree $\nu=2$.
        In such a configuration, the system would reach $7,077,888,000$ degrees of freedom, and under the hardware resources detailed below, the simulation time would be on the order of months.
        Consequently, we evaluate a shifted, coarser sequence of full-grid resolutions $\bm{\alpha}=\{40,80,160,320\}$.
        However, we maintain the higher effective resolutions for the combination technique to explicitly demonstrate its advantages in reducing both simulation time and memory footprint.

        \begin{table}[ht]
            \centering
            \caption{Convergence analysis for the $3$D linear convection-diffusion problem: Full Grid vs. Combination Technique}
            \setlength{\tabcolsep}{6pt} 
            \begin{tabular}{ll S[table-format=1.2e-1,scientific-notation=true] S[table-format=1.2] S[table-format=1.2e-1,scientific-notation=true] S[table-format=1.2] S[table-format=1.2e-1,scientific-notation=true]}
                \toprule
                Case & Config. & {$||e||_1$} & {EOC} & {$||e||_2$} & {EOC} & {Time (ms)} \\
                \midrule
                \multicolumn{7}{l}{\textbf{Polynomial degree $\nu=1$}}\\
                \textit{Full Grid} & $\bm{\alpha}=40$  & 2.576e-2 & {--} & 2.13e-2 & {--} & 2.31e3 \\
                                   & $\bm{\alpha}=80$ & 5.671e-3 & 2.18 & 4.61e-3 & 2.20 & 4.65e4 \\
                                   & $\bm{\alpha}=160$ & 1.372e-3 & 2.05 & 1.10e-3 & 2.07 & 7.74e5 \\
                                   & $\bm{\alpha}=320$ & 3.404e-4 & 2.01 & 2.70e-4 & 2.02 & 4.97e6 \\
                \noalign{\smallskip}
                \textit{CT}        & $N_B=10$ & 2.788e-2 & {--} & 1.95e-2 & {--} & 3.87e3 \\
                Max. level $N_L=3$ & $N_B=20$ & 3.537e-3 & 2.98 & 2.93e-3 & 2.73 & 1.78e5 \\
                                   & $N_B=40$ & 8.519e-4 & 2.05 & 7.17e-4 & 2.03 & 1.34e6 \\
                                   & $N_B=80$ & 2.119e-4 & 2.01 & 1.78e-4 & 2.01 & 1.80e7 \\
                \noalign{\smallskip}
                \textit{CT}        & $N_L=3$  & 2.788e-2 & {--} & 1.95e-2 & {--} & 4.23e3 \\
                Base ref. $N_B=10$ & $N_L=4$  & 7.737e-3 & 1.85 & 5.48e-3 & 1.83 & 5.14e4 \\
                                   & $N_L=5$  & 2.096e-3 & 1.88 & 1.49e-3 & 1.88 & 3.28e5 \\
                                   & $N_L=6$  & 5.577e-4 & 1.91 & 3.95e-4 & 1.91 & 1.74e6 \\
                \midrule
                \multicolumn{7}{l}{\textbf{Polynomial degree $\nu=2$}}\\
                \textit{Full Grid} & $\bm{\alpha}=40$  & 1.240e-3 & {--} & 1.01e-3 & {--} & 2.84e4 \\
                                   & $\bm{\alpha}=80$ & 1.549e-4 & 3.00 & 1.26e-4 & 3.00 & 5.24e5 \\
                                   & $\bm{\alpha}=160$ & 1.934e-5 & 3.00 & 1.58e-5 & 3.00 & 4.41e6 \\
                                   & $\bm{\alpha}=320$ & 2.429e-6 & 2.99 & 1.97e-6 & 3.00 & 2.08e8 \\
                \noalign{\smallskip}
                \textit{CT}        & $N_B=10$ & 2.555e-4 & {--} & 1.88e-4 & {--} & 2.25e5 \\
                Max. level $N_L=3$ & $N_B=20$ & 2.369e-5 & 3.43 & 1.96e-5 & 3.26 & 3.10e6 \\
                                   & $N_B=40$ & 2.981e-6 & 2.99 & 2.47e-6 & 2.99 & 3.71e7 \\
                                   & $N_B=80$ & 3.841e-7 & 2.96 & 3.11e-7 & 2.99 & 1.69e8 \\
                \noalign{\smallskip}
                \textit{CT}        & $N_L=3$  & 2.555e-4 & {--} & 1.88e-4 & {--} & 2.22e5 \\
                Base ref. $N_B=10$ & $N_L=4$  & 3.441e-5 & 2.89 & 2.58e-5 & 2.87 & 1.21e6 \\
                                   & $N_L=5$  & 4.595e-6 & 2.90 & 3.47e-6 & 2.89 & 4.92e6 \\
                                   & $N_L=6$  & 6.089e-7 & 2.92 & 4.62e-7 & 2.91 & 2.58e7 \\
                \bottomrule
            \end{tabular}
            \label{tab:ldg:experiments:linear-pde:3d}
        \end{table}

        Furthermore, explicit matrix assembly is no longer viable across all cases due to excessive RAM requirements.
        Therefore, simulations for the finest levels of the full-grid approach, as well as for the combination technique with varying base refinement, are executed using a matrix-free framework without preconditioning.
        For the remaining, less demanding scenarios, we retain explicit matrix assembly but can no longer rely \joseanto{on an AMG} preconditioner\joseanto{, which  demands substantial memory.}
        Instead, we employ an Additive Schwarz preconditioner, which scales efficiently in parallel environments, coupled with an incomplete LU (ILU) local sub-solver.
        Finally, the workload for all $3$D meshes was distributed across $100$ MPI processes, even though, as in the $2$D case, this may not represent the optimal parallelization strategy for every individual resolution.
         
        Table \ref{tab:ldg:experiments:linear-pde:3d} presents the error, convergence rates, and simulation times for each considered $3$D scenario.
        The errors are computed following the same methodology established for the $2$D case. 
        Overall, the error patterns observed in the $2$D scenarios are preserved. 
        For the full-grid scenario, the expected convergence rates are successfully achieved for both polynomial degrees across all norms, including the $L^{\infty}$ norm.
        At the finest resolution, the maximum absolute errors are $1.55\times10^{-3}$ for $\nu=1$ and $6.70\times 10^{-6}$ for $\nu=2$.
        Similarly, the combination technique with varying base refinement successfully recovers the optimal convergence rates, preserving the expected theoretical scaling even in the $L^{\infty}$ norm, with maximum absolute errors of $7.47\times 10^{-4}$ and $1.18\times 10^{-6}$, respectively.

        Regarding the combination technique with varying maximum levels of refinement, we observe the expected degradation in the convergence rates.
        Notably, this degradation is more pronounced than in the $2$D case. 
        This behavior aligns perfectly with the theoretical estimates, as the logarithmic penalty term worsens with increasing dimensionality.
        \joseanto{The reduction in the convergence rate} is particularly severe in the $L^{\infty}$ norm, where it drops to $1.42$ for $\nu=1$ and $2.65$ for $\nu=2$ at the finest level.
        Correspondingly, the maximum absolute errors obtained in this configuration are $4.83\times 10^{-3}$ and $7.98\times 10^{-6}$.

        \begin{figure}[htbp]
            \centering
            \begin{subfigure}[b]{0.49\textwidth}
                \centering
                \includegraphics[width=\linewidth]{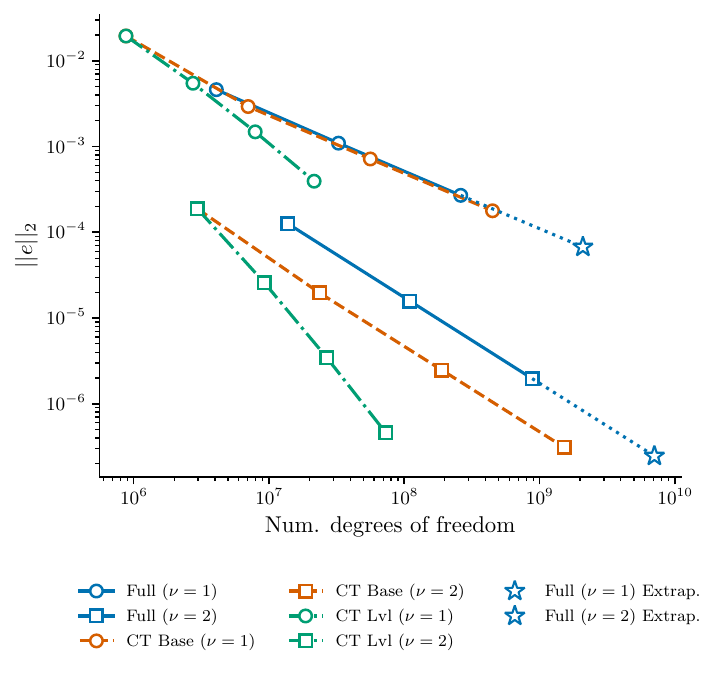}
                \caption{Efficiency in terms of DoFs.}
                \label{fig:ldg:experiments:linear-pde:3d:efficiency:DoF}
            \end{subfigure}
            \hfill 
            \begin{subfigure}[b]{0.49\textwidth}
            \centering
                \includegraphics[width=\linewidth]{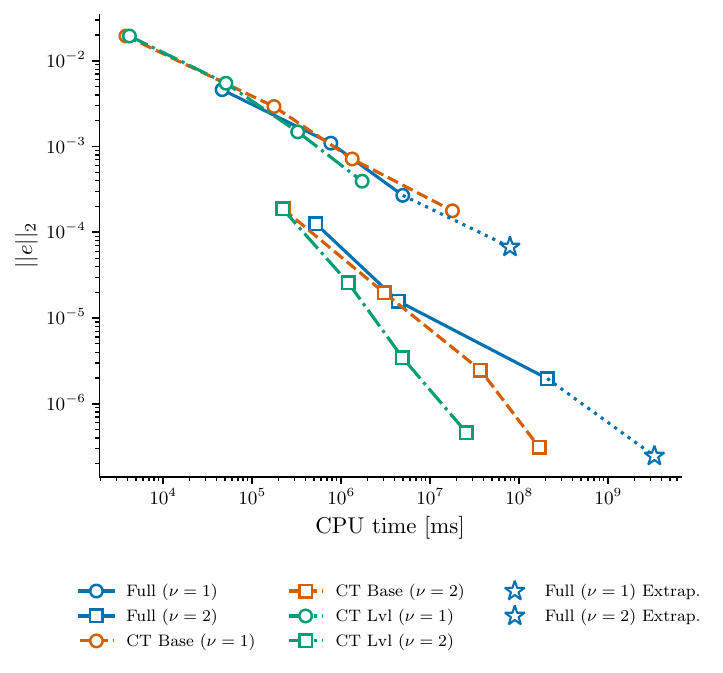}
            \caption{Efficiency in terms of CPU time.}
            \label{fig:ldg:experiments:linear-pde:3d:efficiency:Time}
            \end{subfigure}
        \caption{Efficiency comparison between the full-grid IMEX-LDG scheme and the sparse-grid combination technique for the 3D linear convection-diffusion problem. Data derived from Table \ref{tab:ldg:experiments:linear-pde:3d}.}
        \label{fig:ldg:experiments:linear-pde:3d:efficiency}
        \end{figure}

        The efficiency of each scenario is depicted in Figure \ref{fig:ldg:experiments:linear-pde:3d:efficiency}, following the same considerations established for the $2$D case.
        To ensure a fair comparison and facilitate a clearer understanding of the efficiency gains, we exclude the lowest-resolution simulation entry from the full-grid table.
        Instead, we provide an optimistic estimation of the error, simulation time, and degrees of freedom for the unresolved finest full-grid case ($\bm{\alpha}=640$).
        This estimation is based on the theoretical scaling of a second- or third- order method. 
        Refining the mesh by a factor of two in $3$D yields eight times more elements, requires twice the number of time steps due to the CFL condition, and roughly doubles the computational effort required by the linear solver.
        Furthermore, this is an optimistic projection, as it ignores the parallel communication overhead that would be incurred at that scale.

        The general conclusions drawn from the $2$D efficiency analysis hold true for the $3$D case, but the benefits of the sparse-grid combination technique are now significantly more pronounced.
        Specifically, when varying the maximum refinement level, the combination technique vastly outperforms the full-grid approach for the $\nu=2$ case, achieving comparable accuracy with an approximate $100$ times reduction in both the number of degrees of freedom and simulation time.

        However, the efficiency plots also reveal that the combination technique, in all its considered variants, is no longer computationally advantageous for the linear polynomial case ($\nu=1$).
        We emphasize that this limitation does not stem from a degradation in the convergence rate, which, as demonstrated in the experiment, remains robust.
        Rather, the inefficiency arises from the interplay between the initial error and the parallel communication overhead.

        From a numerical perspective, while the asymptotic scaling is preserved, approximating the multidimensional curvature of a Gaussian pulse on the highly stretched elements of an anisotropic $3$D mesh yields a severe baseline error for purely linear basis functions.
        Consequently, achieving a competitive target accuracy with $\nu=1$ requires pushing the refinement to much deeper levels.
        In contrast, quadratic polynomials ($\nu=2$) possess the intrinsic curvature necessary to keep a remarkably low error constant even on these highly semi-coarsened grids.
        Furthermore, the efficiency gap is increased by the parallel execution setup. 
        To ensure a uniform comparison, all meshes were distributed across $100$ MPI processes. 
        While this scale is mandatory for the finest full-grid resolutions, such a parallelization \joseanto{of the significantly} smaller subgrids of the combination technique degrades the ratio between computation and communication.

        Lastly, we apply the same example to a $5$-dimensional case to demonstrate the capabilities of our proposal. 
        Since the number of degrees of freedom becomes prohibitively large for the full IMEX-LDG scheme, we exclusively adopt the sparse-grid combination technique. 
        All problem parameters are kept constant, with the exception of the domain. 
        Because the primary goal here is not to achieve the extremely low errors required for studying convergence orders, we can utilize a smaller domain without a significant loss of accuracy due to the periodic wrap. 
        Thus, we have $\Omega=(-1.15, 2.45)^5$.
        
        Since analyzing the error throughout the entire $5$D domain is computationally challenging, we evaluate the simulation results on the $(x_1,x_2)$-hyperplane, while fixing the non-varying dimensions at the center of the domain, i.e., $x_3=x_4=x_5=0.65$.
        For this experiment, we select a low-resolution framework with $N_b=3$ to avoid the fill-in effect produced by base refinement, and $N_L=7$, which is equivalent to a full mesh of $\bm{\alpha}=384$ elements per dimension.
        To ensure robust extrapolation properties given the information obtained in the $3$D case, we perform the simulation using polynomials of degree $\nu=2$ to better capture the pulse's geometry.
        This configuration requires solving a total of $771$ subgrids, representing a combined $3,608,602,488$ degrees of freedom.
        Due to the high number of degrees of freedom per subgrid, we employ a matrix-free context without preconditioning, using a relative error tolerance of $10^{-7}$ as the stopping criterion for the linear system. 
        Under this configuration, the total simulation time was approximately $45.5$ hours, with the workload distributed across $416$ processes.

        \begin{figure}[htbp]
         \centering
         \begin{subfigure}[b]{0.48\textwidth}
            \centering
             \includegraphics[width=\textwidth]{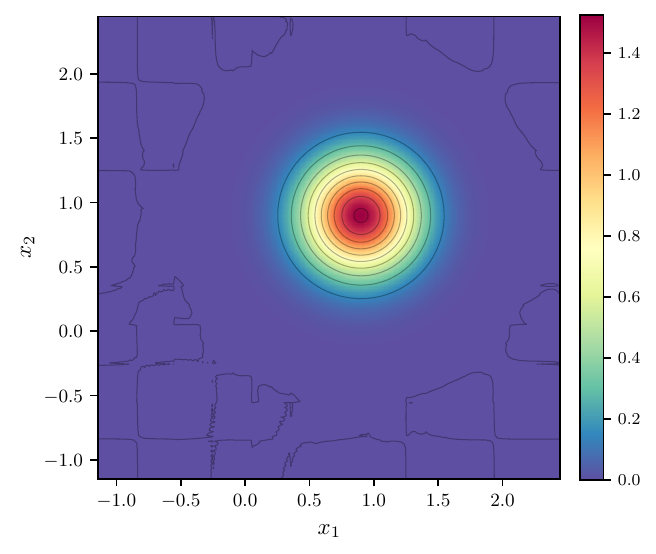}
             \caption{Numerical solution.}
              \label{fig:ldg:experiments:linear-pde:5d:sol} 
         \end{subfigure}
         \hfill
         \begin{subfigure}{0.48\textwidth}
            \centering
            \includegraphics[width=\textwidth]{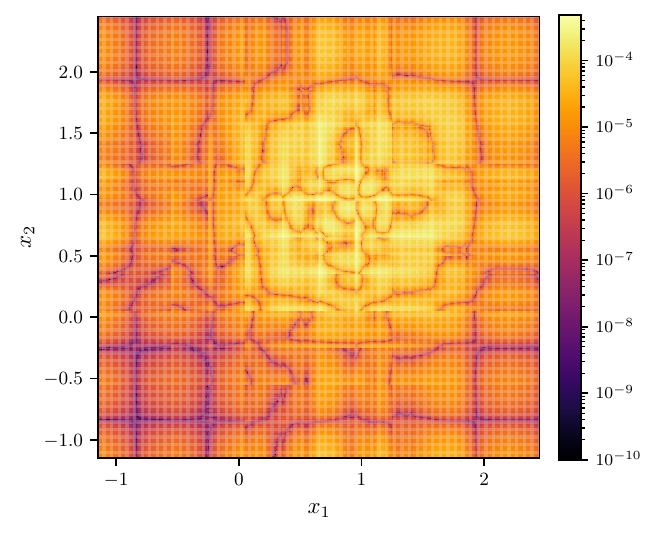}
             \caption{Absolute error.}
              \label{fig:ldg:experiments:linear-pde:5d:err} 
         \end{subfigure}
         \caption{Two-dimensional slice of the numerical solution and absolute error for the $5$D linear convection-diffusion example. The multidimensional domain is intersected at the geometric center of the non-varying dimensions, $x_3=x_4=x_5=0.65$.}
         \label{fig:ldg:experiments:linear-pde:5d}
        \end{figure}

        As we are no longer dealing with the full-grid case, we approximate the functional error \eqref{eq:ldg:experiments:lp-error} on the hyperplane by means of its discretized version.
        The error is evaluated over a discrete set of points $\mathcal{X}\subset \Omega$.
        Thus, \joseanto{the normalized} discrete $L^p$ error norm is then defined as:
        \begin{equation}
            \label{eq:ldg:experiments:discrete_error}
            \text{E}_p=\left(\dfrac{1}{|\mathcal{X}|}\sum_{\xv\in\mathcal{X}}|u^*(T,\xv) - u_h(T,\xv)|^p\right)^{\frac{1}{p}}
        \end{equation} 
        which, for $p=1$ and $p=2$, corresponds to the metrics known as mean absolute error (MAE) and root mean squared error (RMSE), respectively. 
        The discrete $L^{\infty}$-norm, $\text{E}_{\infty}$ is computed using expression \eqref{eq:ldg:experiments:linfty-error} applied to the corresponding discrete subdomain.
        
        In this instance, the evaluation points on the hyperplane are chosen to be the cell centers of the equivalent full-resolution case. 
        For this example, we obtain the following errors: $\text{E}_1=2.219\times 10^{-5}$, $\text{E}_2=4.159\times 10^{-5}$, and $\text{E}_{\infty}=4.851\times 10^{-4}$, which indicates strong approximation properties for the linear case.
        Additionally, Figures \ref{fig:ldg:experiments:linear-pde:5d:sol} and \ref{fig:ldg:experiments:linear-pde:5d:err} depict the resulting numerical solution and the absolute error with respect to the exact solution on the hyperplane.
        As expected, the largest errors are located in the vicinity of the bell, and the subgrid pattern is clearly identifiable on the absolute error construction.
        
    \subsection{Viscous Burgers' problem} \label{ldg:experiments:vburgers}
        The second benchmark assesses the performance of the methodology using the viscous Burgers' equation, which introduces a nonlinear advection term characterized by the flux vector $\bm{f}(u)=[\frac{1}{2}u^2,\dots,\frac{1}{2}u^2]^T$.
        To ensure the availability of an analytic solution for arbitrary dimension $d$, we augment the standard equation with the forcing term $s(t,\xv)$.
        Thus, the problem is governed by:
        \begin{equation}
            \label{eq:ldg:experiments:burgers-pde}
            \partial_t u + \sum_{k=1}^d \partial_{x_k}\left(\dfrac{u^2}{2}\right)= a \sum_{k=1}^d \partial_{x_k x_k} u + s(t,\xv),\quad \text{with } s(t,\xv) = \dfrac{d}{2}e^{-2adt}\sin\left(2\sum_{k=1}^dx_k\right),
        \end{equation}
        yielding the analytical solution:
        \begin{equation}
            \label{eq:ldg:experiments:burgers-pde:sol}
            u^*(t,\xv)= e^{-adt}\sin\left(\sum_{k=1}^dx_k\right).
        \end{equation}

        This problem naturally admits periodic boundary conditions on a domain of length $2\pi$ in each coordinate direction. 
        Consequently, we define the computational domain as $\Omega=(-\pi,\pi)^d$.
        For these experiments, we set the diffusion coefficient $a=0.1$ and the final simulation time $T=0.5$.
        
        In this Burgers' case, the maximum characteristic speed is given by $c_v=\max|\bm{f}'(u)|$, a quantity that generally needs to be updated at each time step.
        However, by construction of our test case, the sine wave experiences a continuous decay in amplitude throughout the simulation.
        This allows us to safely bound the global time step restriction using the characteristic speed of the initial condition.
        Thus, the time step is once again selected to satisfy the CFL condition \eqref{eq:ldg:imex:cfl} by setting the constant to $C=0.9$. 
        
        Furthermore, the general experimental setup, including the grid resolutions, polynomial degrees, linear solvers, and parallelization strategy, remains identical to the configuration detailed in Section \ref{ldg:experiments:linear} for the linear convection-diffusion problem.

        \begin{table}[ht]
            \centering
            \caption{Convergence analysis for the $2$D viscous Burgers problem: Full Grid vs. Combination Technique}
            \setlength{\tabcolsep}{6pt} 
            \begin{tabular}{ll S[table-format=1.2e-1,scientific-notation=true] S[table-format=1.2] S[table-format=1.2e-1,scientific-notation=true] S[table-format=1.2] S[table-format=1.2e-1,scientific-notation=true]}
                \toprule
                Case & Config. & {$||e||_1$} & {EOC} & {$||e||_2$} & {EOC} & {Time (ms)} \\
                \midrule
                \multicolumn{7}{l}{\textbf{Polynomial degree $\nu=1$}}\\
                \textit{Full Grid} & $\bm{\alpha}=80$  & 6.268e-3 & {--} & 1.54e-3 & {--} & 3.17e2 \\
                                   & $\bm{\alpha}=160$ & 1.624e-3 & 1.95 & 4.02e-4 & 1.94 & 1.29e3 \\
                                   & $\bm{\alpha}=320$ & 4.136e-4 & 1.97 & 1.03e-4 & 1.97 & 1.17e4 \\
                                   & $\bm{\alpha}=640$ & 1.044e-4 & 1.99 & 2.60e-5 & 1.98 & 8.61e4 \\
                \noalign{\smallskip}
                \textit{CT}        & $N_B=10$ & 6.903e-3 & {--} & 1.58e-3 & {--} & 3.11e2 \\
                Max. level $N_L=3$ & $N_B=20$ & 1.648e-3 & 2.07 & 4.03e-4 & 1.97 & 1.13e3 \\
                                   & $N_B=40$ & 4.145e-4 & 1.99 & 1.03e-4 & 1.97 & 6.20e3 \\
                                   & $N_B=80$ & 1.044e-4 & 1.99 & 2.60e-5 & 1.98 & 7.12e4 \\
                \noalign{\smallskip}
                \textit{CT}        & $N_L=3$  & 6.903e-3 & {--} & 1.58e-3 & {--} & 3.55e2 \\
                Base ref. $N_B=10$ & $N_L=4$  & 1.819e-3 & 1.92 & 4.12e-4 & 1.93 & 9.77e2 \\
                                   & $N_L=5$  & 4.693e-4 & 1.95 & 1.06e-4 & 1.96 & 2.66e3 \\
                                   & $N_L=6$  & 1.198e-4 & 1.97 & 2.69e-5 & 1.98 & 8.75e3 \\
                \midrule
                \multicolumn{7}{l}{\textbf{Polynomial degree $\nu=2$}}\\
                \textit{Full Grid} & $\bm{\alpha}=80$  & 4.373e-5 & {--} & 9.98e-6 & {--} & 1.70e3 \\
                                   & $\bm{\alpha}=160$ & 5.565e-6 & 2.97 & 1.26e-6 & 2.98 & 1.45e4 \\
                                   & $\bm{\alpha}=320$ & 7.021e-7 & 2.99 & 1.59e-7 & 2.99 & 1.28e5 \\
                                   & $\bm{\alpha}=640$ & 8.879e-8 & 2.99 & 2.00e-8 & 2.99 & 2.02e6 \\
                \noalign{\smallskip}
                \textit{CT}        & $N_B=10$ & 4.486e-5 & {--} & 1.02e-5 & {--} & 2.20e3 \\
                Max. level $N_L=3$ & $N_B=20$ & 5.568e-6 & 3.01 & 1.27e-6 & 3.01 & 8.84e3 \\
                                   & $N_B=40$ & 7.021e-7 & 2.99 & 1.59e-7 & 2.99 & 4.72e4 \\
                                   & $N_B=80$ & 8.879e-8 & 2.98 & 2.00e-8 & 2.99 & 6.77e5 \\
                \noalign{\smallskip}
                \textit{CT}        & $N_L=3$  & 4.486e-5 & {--} & 1.02e-5 & {--} & 1.37e3 \\
                Base ref. $N_B=10$ & $N_L=4$  & 5.774e-6 & 2.96 & 1.32e-6 & 2.95 & 6.15e3 \\
                                   & $N_L=5$  & 7.389e-7 & 2.97 & 1.71e-7 & 2.95 & 2.46e4 \\
                                   & $N_L=6$  & 1.074e-7 & 2.78 & 2.78e-8 & 2.62 & 5.54e4 \\
                \bottomrule
            \end{tabular}
            \label{tab:ldg:experiments:burgers:2d}
        \end{table}

        Table \ref{tab:ldg:experiments:burgers:2d} presents the error, convergence rates, and simulation times for the two-dimensional Burgers' case.
        The full-grid approach achieves the expected performance across all scenarios, delivering the optimal theoretical convergence rates for both polynomial degrees.
        In addition, the full-grid scheme successfully maintains the expected convergence rate in the $L^{\infty}$ norm, yielding maximum absolute errors of $8.41\times 10^{-6}$ and $7.14\times 10^{-9}$ for $\nu=1$ and $\nu=2$, respectively.
        For the combination technique with varying base refinement, the performance metrics generally align with theoretical expectations, successfully recovering the optimal convergence rates for the $L^1$ and $L^2$ norms.
        However, a slight degradation in the $L^{\infty}$ convergence rate is observed, particularly for the quadratic case, where the rate drops to $2.25$. 
        Despite this, the maximum absolute errors remain well-controlled at $8.42\times 10^{-6}$ and $1.24\times 10^{-8}$.

        \begin{figure}[htbp]
            \centering
            \begin{subfigure}[b]{0.49\textwidth}
                \centering
                \includegraphics[width=\linewidth]{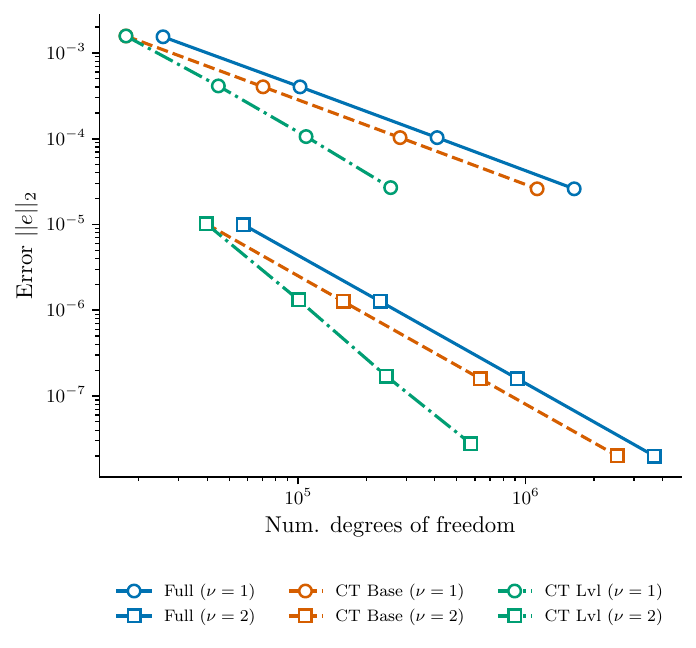}
                \caption{Efficiency in terms of DoFs.}
                \label{fig:ldg:experiments:burgers:2d:efficiency:DoF}
            \end{subfigure}
            \hfill 
            \begin{subfigure}[b]{0.49\textwidth}
            \centering
                \includegraphics[width=\linewidth]{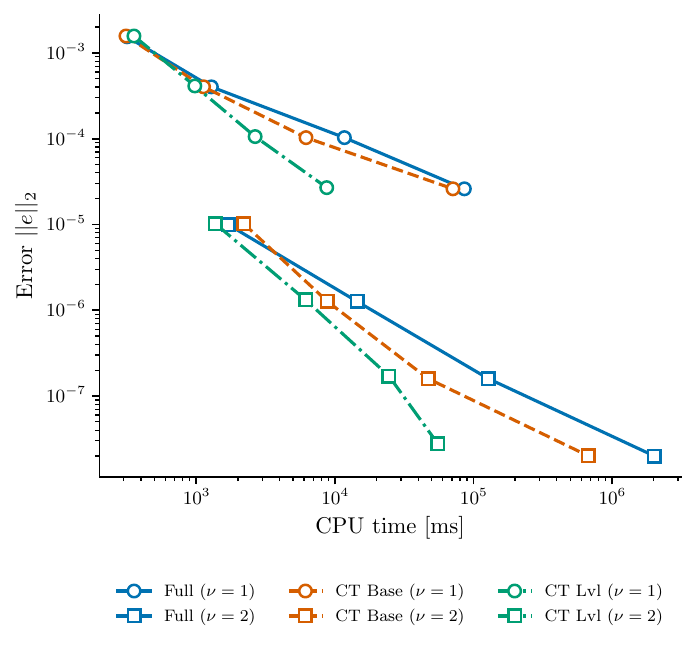}
            \caption{Efficiency in terms of CPU time.}
            \label{fig:ldg:experiments:burgers:2d:efficiency:Time}
            \end{subfigure}
        \caption{Efficiency comparison between the full-grid IMEX-LDG scheme and the sparse-grid combination technique for the 2D viscous Burgers problem. Data derived from Table \ref{tab:ldg:experiments:burgers:2d}.}
        \label{fig:ldg:experiments:burgers:2d:efficiency}
        \end{figure}

        The most significant deviation compared to the $2$D linear benchmark is observed in the combination technique utilizing a varying maximum refinement level. 
        While the $\nu=1$ case behaves similarly to its linear PDE counterpart, the configuration with $\nu=2$ experiences a more pronounced drop in its convergence rates.
        This degradation is particularly evident when we account for the $L^{\infty}$ norm, where the expected convergence rate drops below $2$, settling specifically as $1.71$.
        Nevertheless, the maximum absolute errors remain highly acceptable, reaching $1.72\times 10^{-5}$ for $\nu=1$ and $1.74\times 10^{-7}$ for $\nu=2$.

        Regarding the computational efficiency, the combination technique with varying maximum refinement levels clearly outperforms the alternative approaches with respect to the $L^2$ error norm, highlighting the substantial advantage of the sparse-grid methodology in reducing simulation time.
        These results also demonstrate that the introduction of a nonlinear advection flux does not severely impact the global approximation quality across the domain, at least in the two-dimensional context.
        The corresponding efficiency curves are depicted in Figure \ref{fig:ldg:experiments:burgers:2d:efficiency}, which corroborate the general conclusions drawn previously for the linear case.
        
        In the $3$D Burgers' case, the performance landscape shifts significantly.
        This is particularly true for the combination technique with varying \joseanto{refinement levels}, where severe catastrophic cancellation emerges when employing quadratic elements.
        For the full-grid cases, the theoretical convergence rates are generally preserved, although the $\nu=1$ configuration has not yet reached full asymptotic convergence in the $L^1$ and $L^2$ norms.
        Nevertheless, the optimal $L^{\infty}$ convergence rates are successfully attained for both polynomial degrees, yielding maximum absolute errors of $4.84\times 10^{-5}$ and $8.19\times 10^{-8}$ at the finest resolutions.
        When applying the combination technique with a varying base refinement, the case of linear elements continues to perform robustly.
        However, the $\nu=2$ case exhibits a pronounced degradation in convergence.
        An analysis of the $L^{\infty}$ errors reveals a loss of monotonicity.
        While the $\nu=2$ scheme \joseanto{initially demonstrates} rapid error reduction (with convergence rates of $5.16$ and $4.09$ in the early stages), the final refinement step actually amplifies the maximum error, resulting in a negative convergence rate of $-0.73$.
        The maximum absolute errors for these configurations are $1.41\times 10^{-4}$ and $2.20\times 10^{-7}$.

        This performance degradation becomes even more severe when analyzing the combination technique with varying maximum refinement levels for $\nu=2$.
        Here, the convergence rate drops to approximately $1$ across all norms, falling drastically short of theoretical expectations.
        However, this behavior is not entirely unexpected.
        Previous studies, such as \cite{Lastdrager2001}, have empirically demonstrated that the sparse-grid combination technique struggles with the Burgers' equation when the advection features are not strictly aligned with the computational grid.

         \begin{table}[ht]
            \centering
            \caption{Convergence analysis for the $3$D viscous Burgers problem: Full Grid vs. Combination Technique}
            \setlength{\tabcolsep}{6pt} 
            \begin{tabular}{ll S[table-format=1.2e-1,scientific-notation=true] S[table-format=1.2] S[table-format=1.2e-1,scientific-notation=true] S[table-format=1.2] S[table-format=1.2e-1,scientific-notation=true]}
                \toprule
                Case & Config. & {$||e||_1$} & {EOC} & {$||e||_2$} & {EOC} & {Time (ms)} \\
                \midrule
                \multicolumn{7}{l}{\textbf{Polynomial degree $\nu=1$}}\\
                \textit{Full Grid} & $\bm{\alpha}=40$  & 2.062e-1 & {--} & 1.69e-2 & {--} & 2.18e3 \\
                                   & $\bm{\alpha}=80$ & 5.563e-2 & 1.89 & 4.52e-3 & 1.90 & 4.60e4 \\
                                   & $\bm{\alpha}=160$ & 1.449e-2 & 1.94 & 1.18e-3 & 1.94 & 5.18e5 \\
                                   & $\bm{\alpha}=320$ & 3.700e-3 & 1.97 & 3.01e-4 & 1.97 & 3.93e6 \\
                \noalign{\smallskip}
                \textit{CT}        & $N_B=10$ & 6.188e-2 & {--} & 4.89e-3 & {--} & 3.80e3 \\
                Max. level $N_L=3$ & $N_B=20$ & 1.453e-2 & 2.09 & 1.18e-3 & 2.05 & 6.45e4 \\
                                   & $N_B=40$ & 3.700e-3 & 1.97 & 3.01e-4 & 1.97 & 1.50e6 \\
                                   & $N_B=80$ & 9.351e-4 & 1.98 & 7.61e-5 & 1.98 & 1.48e7 \\
                \noalign{\smallskip}
                \textit{CT}        & $N_L=3$  & 6.189e-2 & {--} & 4.89e-3 & {--} & 4.37e3 \\
                Base ref. $N_B=10$ & $N_L=4$  & 1.85e-2 & 1.74 & 1.42e-3 & 1.78 & 5.07e4 \\
                                   & $N_L=5$  & 5.394e-3 & 1.78 & 4.19e-4 & 1.76 & 3.34e5 \\
                                   & $N_L=6$  & 1.816e-3 & 1.57 & 1.50e-4 & 1.48 & 1.52e6 \\
                \midrule
                \multicolumn{7}{l}{\textbf{Polynomial degree $\nu=2$}}\\
                \textit{Full Grid} & $\bm{\alpha}=40$  & 2.698e-3 & {--} & 2.28e-4 & {--} & 2.65e4 \\
                                   & $\bm{\alpha}=80$ & 3.488e-4 & 3.00 & 2.92e-5 & 2.97 & 4.68e5 \\
                                   & $\bm{\alpha}=160$ & 4.439e-5 & 3.00 & 3.70e-6 & 2.98 & 6.48e6 \\
                                   & $\bm{\alpha}=320$ & 5.607e-6 & 2.99 & 4.66e-7 & 2.99 & 1.15e8 \\
                \noalign{\smallskip}
                \textit{CT}        & $N_B=10$ & 7.542e-4 & {--} & 8.98e-5 & {--} & 2.19e5 \\
                Max. level $N_L=3$ & $N_B=20$ & 4.768e-5 & 3.98 & 4.26e-6 & 4.40 & 3.08e6 \\
                                   & $N_B=40$ & 5.613e-6 & 3.09 & 4.68e-7 & 3.19 & 3.76e7 \\
                                   & $N_B=80$ & 8.558e-7 & 2.71 & 7.13e-8 & 2.71 & 1.34e8 \\
                \noalign{\smallskip}
                \textit{CT}        & $N_L=3$  & 7.542e-4 & {--} & 8.98e-5 & {--} & 2.22e5 \\
                Base ref. $N_B=10$ & $N_L=4$  & 3.375e-4 & 1.16 & 4.42e-5 & 1.02 & 1.22e6 \\
                                   & $N_L=5$  & 1.320e-4 & 1.35 & 1.75e-5 & 1.34 & 6.16e6 \\
                                   & $N_L=6$  & 7.315e-5 & 0.85 & 9.73e-6 & 0.85 & 1.60e7 \\
                \bottomrule
            \end{tabular}
            \label{tab:ldg:experiments:burgers:3d}
        \end{table}

        The failure of the third-order scheme ($\nu=2$) can be explained through sparse-grid approximation theory \cite{Bungartz_Griebel_2004}.
        Elevating the polynomial degree from $\nu=1$ to $\nu=2$ imposes stricter regularity requirements, specifically demanding higher-order bounded \joseanto{mixed derivatives}.
        In the context of the nonlinear Burgers' flux, the integrals in the weak formulation generate cross-polynomial terms up to degree $4$. 
        Although in practice we employ an interpolated flux, this projection still contributes to the aliasing effect.
        These under-resolved high-frequency components suffer from aliasing \cite{Hesthaven2007}, a phenomenon exacerbated by the fact that the highly anisotropic meshes force the computations into a marginally resolved regime. 
        This aliasing contaminates the low-frequency modes in the discrete representation (an effect that is more disruptive for $\nu=2$ than for $\nu=1$).
        Furthermore, because the advective flux is not grid-aligned, characteristic information must propagate diagonally across highly anisotropic cells of the sparse subgrids.
        For $\nu=1$, the inherent numerical dissipation introduced by the LDG numerical fluxes \cite{Cockburn1998} effectively smears the numerical error, acting as a stabilizing diffusive mechanism.
        Conversely, the high-precision $\nu=2$ scheme is inherently less dissipative, causing its leading truncation error to be dominated by dispersion.
        This dispersion generates spurious high-frequency numerical oscillations. 
        When the sparse-grid formula linearly combines these subgrid solutions, it fails to coherently cancel these terms, leading instead to a catastrophic cancellation of the physical solution.
        Ultimately, the loss of the linear superposition property in this nonlinear regime causes the leading directional errors to inextricably mix, violating the fundamental error-splitting assumption upon which the combination technique relies.
        
\begin{figure}[htbp]
            \centering
            \begin{subfigure}[b]{0.49\textwidth}
                \centering
                \includegraphics[width=\linewidth]{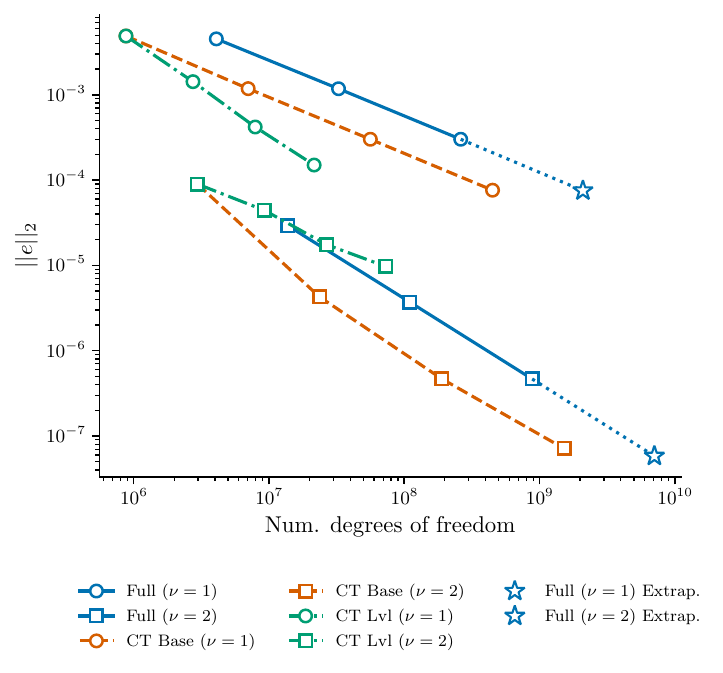}
                \caption{Efficiency in terms of DoFs.}
                \label{fig:ldg:experiments:burgers:3d:efficiency:DoF}
            \end{subfigure}
            \hfill 
            \begin{subfigure}[b]{0.49\textwidth}
            \centering
                \includegraphics[width=\linewidth]{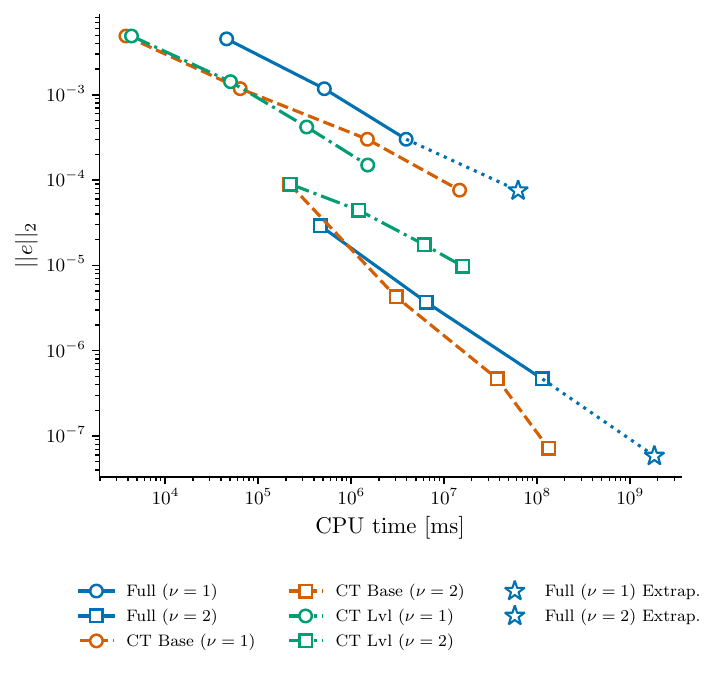}
            \caption{Efficiency in terms of CPU time.}
            \label{fig:ldg:experiments:burgers:3d:efficiency:Time}
            \end{subfigure}
        \caption{Efficiency comparison between the full-grid IMEX-LDG scheme and the sparse-grid combination technique for the 3D viscous Burgers problem. Data derived from Table \ref{tab:ldg:experiments:burgers:3d}.}
        \label{fig:ldg:experiments:burgers:3d:efficiency}
        \end{figure}

        The preceding analysis provides insight for interpreting the efficiency curves depicted in Figure \ref{fig:ldg:experiments:burgers:3d:efficiency}.
        When applying the combination technique with a varying base refinement, the initial grid resolution proves sufficient to accurately capture the underlying physics before the sparse-grid extrapolation is applied.
        In this scenario, the methodology avoids severe aliasing and successfully demonstrates superior computational efficiency over the full-grid scheme for both linear and quadratic elements.
        Conversely, the efficiency plots reflect the failure of the combination technique when varying the maximum refinement level with quadratic polynomials.
        As detailed above, the inability to coherently resolve the nonlinearity completely nullifies any theoretical reduction in computational cost, \joseanto{rendering it} the least efficient option among all tested quadratic configurations.
        
        In order to present the $5$D Burgers' case, the configuration remains consistent with the previously described $5$D linear setup, though the domain size is preserved to ensure periodicity.
        We evaluate the numerical results of the IMEX-LDG sparse-grid combination technique on the $(x_1,x_2)$-hyperplane, while the remaining variables are fixed at the geometric center of the domain, i.e., $x_3=x_4=x_5=0$.
        We maintain the low-resolution framework with $N_b=3$ and $N_L=7$, yielding an equivalent full-grid resolution as established in the linear benchmark.
        For this case, we use quadratic polynomials.
        Although linear elements provided convergence rates closer to theoretical expectations, higher-order elements deliver significantly lower absolute errors with fewer degrees of freedom, a trend consistently observed in the efficiency plots.
        Under this configuration, the total simulation time was approximately $23.45$ hours using the same workload distribution as the linear example.
        
\begin{figure}[htbp]
             \centering
             \begin{subfigure}[b]{0.48\textwidth}
                \centering
                 \includegraphics[width=\textwidth]{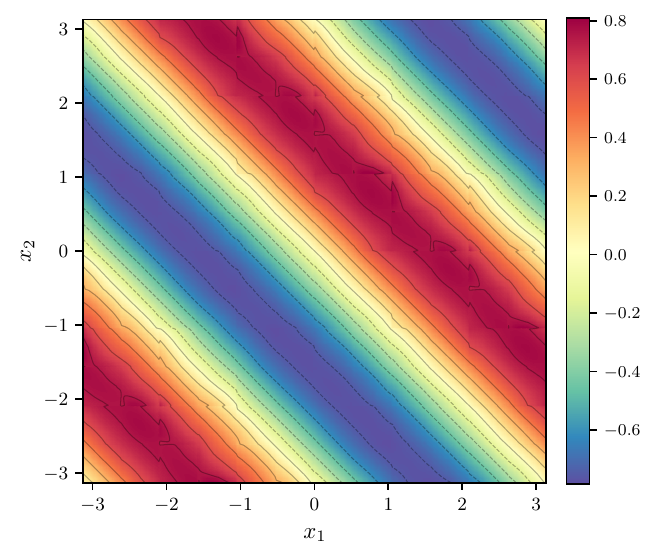}
                 \caption{Numerical solution.}
                  \label{fig:ldg:experiments:burgers:5d:sol} 
             \end{subfigure}
             \hfill
             \begin{subfigure}{0.48\textwidth}
                \centering
                \includegraphics[width=\textwidth]{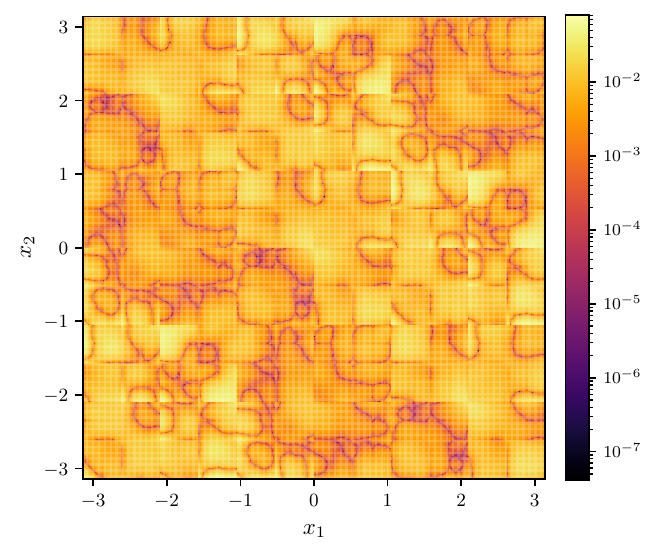}
                 \caption{Absolute error.}
                  \label{fig:ldg:experiments:burgers:5d:err} 
             \end{subfigure}
             \caption{\joseanto{Two-dimensional} slice of the numerical solution and absolute error for the $5$D Burgers' example. The multidimensional domain is intersected at the geometric center of the non-varying dimensions, $x_3=x_4=x_5=0$.}
             \label{fig:ldg:experiments:burgers:5d}
        \end{figure}

        For this experiment, we obtain the following errors against the analytical solution: $\text{E}_1=7.127\times 10^{-3}$, $\text{E}_2=9.89\times 10^{-3}$, and the maximum absolute error of $\text{E}_{\infty}=8.07\times 10^{-2}$.
        These values are substantially higher than those recorded for the analogous linear case, further demonstrating the increasing complexity of the Burgers' dynamics.
        This is better illustrated by the numerical solution and point-wise absolute error depicted in Figure \ref{fig:ldg:experiments:burgers:5d}.
        In the solution heatmap, the $10$ superimposed contour levels provide clear evidence of the aliasing effect. 
        The positive crests of the sine wave exhibit irregular contours, and the error map confirms that the maximum absolute errors are localized precisely at these positions.
        Interestingly, this phenomenon is not observed in the negative valleys.
        
        This spatial asymmetry directly builds upon the degradation mechanism established in the $3$D analysis.
        While the generation of under-resolved high-frequency harmonics and subsequent aliasing explains the global loss of convergence, the localized manifestation of these spurious oscillations is driven by the relative phase error intrinsic to high-order DG methods \cite{Hesthaven2007}.
        As the unresolved modes fold back into the lower-degree polynomial space, they undergo a phase shift.
        The superposition of the fundamental physical wave with these phase-shifted, aliased spurious modes yields a spatially asymmetric error distribution.
        Depending on the local phase relationship, this superposition leads to a localized amplification of the error, manifesting as sharp oscillations on the positive crests, or to an attenuation and artificial smoothing in the negative valleys.
        In our context, extracting the slice $x_3=x_4=x_5=0$ places the evaluation plane exactly at a sinusoidal transition for the transverse dimensions, maximizing this phase discrepancy at the crests. 
        Similar localized residual fluctuations are documented in \cite{Wang2016}, confirming that the sparse-grid representation framework is highly susceptible to strong localized, phase-driven residual errors.

        In order to corroborate the theoretical effect discussed above, we simulated the same $5$D Burgers' test case but increasing the base refinement to $N_b=4$.
        From a computational standpoint, this represents a massive increase in the number of degrees of freedom solved during the simulation, reaching   a total of $15,206,621,184$, with an overall execution time of $79.5$ hours.
        As anticipated, we observe a substantial improvement in the error metrics, which now stand at $\text{E}_1=4.92\times 10^{-3}$, $\text{E}_2=7.34\times 10^{-3}$, and a maximum absolute error of $\text{E}_{\infty}=4.88\times 10^{-2}$.
        Notably, this represents an approximate $40\%$ reduction in the maximum absolute error compared to the $N_b=3$ configuration.
        
\begin{figure}[htbp]
             \centering
             \begin{subfigure}[b]{0.48\textwidth}
                \centering
                 \includegraphics[width=\textwidth]{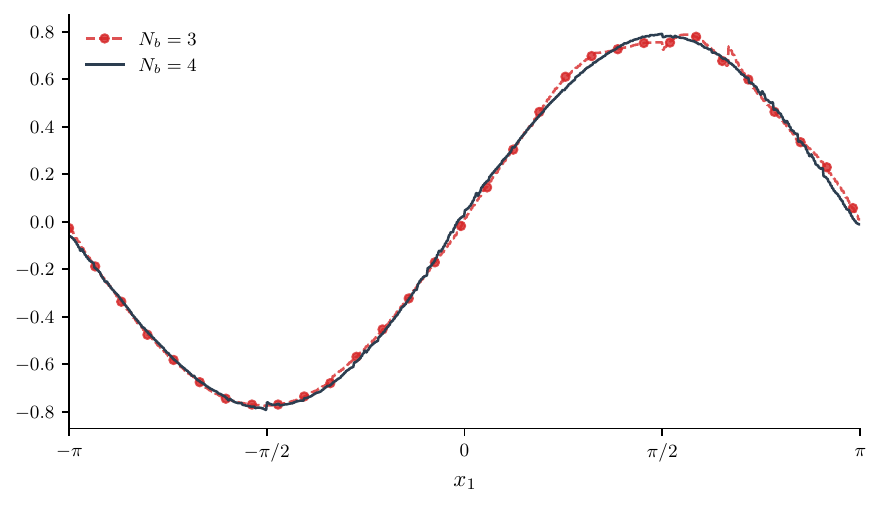}
                 \caption{Cross-section of the numerical approximations.}
                  \label{fig:ldg:experiments:burgers:5d:aliasing:sol} 
             \end{subfigure}
             \hfill
             \begin{subfigure}{0.48\textwidth}
                \centering
                \includegraphics[width=\textwidth]{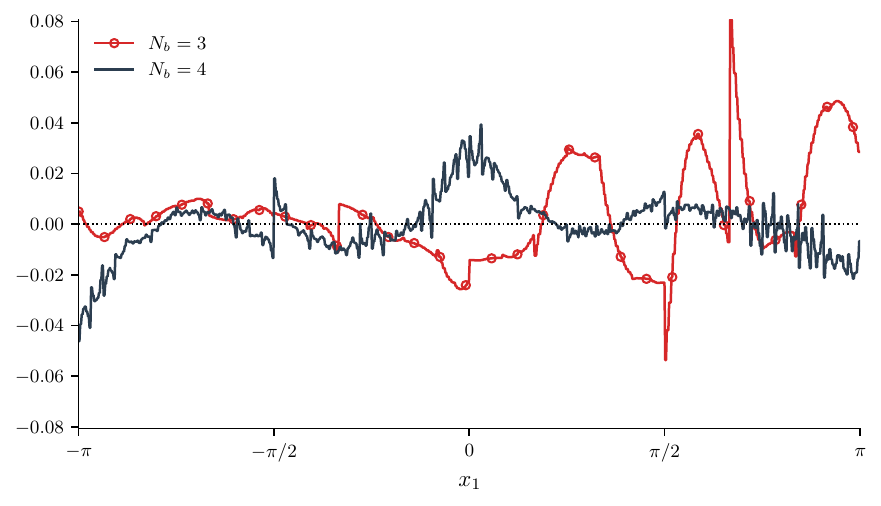}
                 \caption{Pointwise difference profiles.}
                  \label{fig:ldg:experiments:burgers:5d:aliasing:err} 
             \end{subfigure}
             \caption{One-dimensional cross-section of the numerical solution and its corresponding differences for the $5$D Burgers' equation. The multidimensional domain is intersected along the $x_1$-axis at the geometric center $x_2=x_3=x_4=x_5=0$.}
             \label{fig:ldg:experiments:burgers:5d:aliasing}
        \end{figure}

        Increasing the base level imposes a strong damping effect on the nonlinear aliasing, as the higher isotropic base resolution is capable of more accurately capturing the underlying geometry of the physical solution before the sparse-grid extrapolation is applied.
        This mitigating effect is clearly visible in Figure \ref{fig:ldg:experiments:burgers:5d:aliasing}, which plots the one-dimensional cross-section of the numerical approximations and their absolute differences with respect to the exact solution, intersected along the $x_1$-axis at the geometric center $x_2=x_3=x_4=x_5=0$.
        The profiles illustrate the severe oscillatory artifacts present at the lower base resolution, and their effective suppression upon increasing the base grid level.
        
        Furthermore, the difference profiles reveal another interesting phenomenon.
        The spatial locations of the error peaks do not align between the two discretizations.
        This spatial shift occurs because incrementing the base refinement fundamentally alters the nodal distribution of the sparse grid.
        Consequently, this changes the spatial support of the basis functions, modifies the local interpolation nodes, and alters the relative alignment between the structured grid and the nonlinear advective flux.

        The primary drawback of this approach is the drastic increase in computational complexity.
        Operating at a larger base level is exponentially more expensive than incrementing the maximum refinement level.
        However, the fundamental issue is that simply increasing the maximum refinement level $N_L$ fails to resolve the aliasing effect.
        Because a higher $N_L$ merely introduces additional highly anisotropic subgrids without improving the coarse cross-directional resolution,
        it cannot prevent the high-frequency harmonics from folding back into the lower-degree polynomial space.

    \subsection{Advective Allen-Cahn problem with rotational flow} \label{ldg:experiments:ac}   
        We now consider the Allen-Cahn equation with a normalized diffusion coefficient, subject to a prescribed rotational velocity field in two- and three-dimensional spatial domains.
        As a prototypical nonlinear reaction-diffusion PDE, the Allen-Cahn model describes the process of phase separation across various scientific disciplines, including chemistry, population dynamics, and materials science.
        Consequently, the phase variable $u$ is typically bounded within the range $[-1,1]$, where the extrema represent the two distinct physical states.

        Specifically, defining the computational domain as $\Omega=[0,1]^d$ for $d\in\{2,3\}$, the advective mean curvature flow problem simulated in this work is governed by the following initial-boundary value problem: 

        \begin{equation} \label{eq:ldg:experiments:ac:pde}
            \begin{cases}
                \partial_t u + \nabla_{\bm{x}}\cdot \bm{f}(u,\xv) = \Delta u - \dfrac{1}{\epsilon^2}\left(u^3-u\right),& \text{in } (0,T)\times\Omega, \\[7pt]
                u(0,\xv) = \tanh\left(\dfrac{R_0 - ||\xv-\bm{c}||}{\sqrt{2}\epsilon}\right), & \text{in } \Omega, \\[7pt] 
                \nabla_{\xv}u\cdot \bm{n} = 0,& \text{in } (0,T)\times \partial\Omega
            \end{cases}
        \end{equation}
        
        where $\epsilon>0$ dictates the thickness of the diffuse transition region between the pure phases, and $\bm{c}$ denotes the geometric center of the initial spherical interface of radius $R_0=0.25$.
        Specifically, we fix the initial center at $\bm{c}=[0.65, 0.5]^T$ for the $2$D case, and $\bm{c}=[0.65,0.5,0.5]^T$ for the $3$D case.
        The advective flux is driven by the rotational velocity field, defined in $2$D as $\bm{f}(u,\xv) = \left[-\omega(x_2-1/2),\,\omega(x_1-1/2)\right]^Tu$, with $\omega$ being the angular velocity.
        For the $3$D scenarios, this planar rotation is naturally extended by setting the third velocity component to zero.

        Physically, this problem describes the continuous shrinking of a circular ($2$D) or spherical ($3D$) interface driven by mean curvature, while simultaneously rotating around the center of the domain.
        Let $R(t)$ and $V(t)$ denote the effective radius and volume of the phase region at time $t$, respectively.
        As demonstrated in \cite{Ju2015} and references therein, in the sharp-interface limit $\epsilon\to 0$, the theoretical limit radius $R_{\mathrm{ana}}(t)$ is given by $R_{\mathrm{ana}}(t)=\sqrt{R_0^2+2(1-d)t}$.
        This evolution can be deduced directly from the analytical volume $V_{\mathrm{ana}}(t) = \pi \left(R_0^2-2t\right)$ for $d=2$, and $V_{\mathrm{ana}}^{2/3}(t)=(4\pi/3)^{2/3}(R_0^2-4t)$ for $d=3$.  
        Thus, we validate our numerical results by tracking the evolution of the effective simulated radius, which is extracted from the mass reduction.
        This is achieved by computing the discrete volume:
        \begin{equation*}
            V(t)=\int_{\Omega}\dfrac{1}{2}(1+u)\,\text{d}\xv,
        \end{equation*}
        where the integrand shifts and scales the phase variable $u$ to ensure \joseanto{it takes} a value of $1$ inside the evolving region and $0$ outside.
        The simulated radius $R(t)$ is then computed by inverting the aforementioned volume formulas.
        
        This benchmark couples the stiff, nonlinear phase separation dynamics characteristic of \joseanto{the Allen-Cahn} model with linear solid-body rotational transport.
        The configuration poses a challenge for the standard full-grid IMEX-LDG scheme, and even more severe test for the sparse-grid combination technique.
        Specifically, the numerical method must accurately resolve the steep moving interfaces driven by the stiff nonlinear reaction term, while simultaneously preserving the geometric integrity of the phase boundary as it is advected across the highly anisotropic, \joseanto{structure-sparse} subgrids.

        \begin{figure}[htbp]
            \centering
            \begin{minipage}{0.88\textwidth}
                \centering
                \begin{subfigure}{0.24\textwidth}
                    \centering
                    \includegraphics[width=\textwidth]{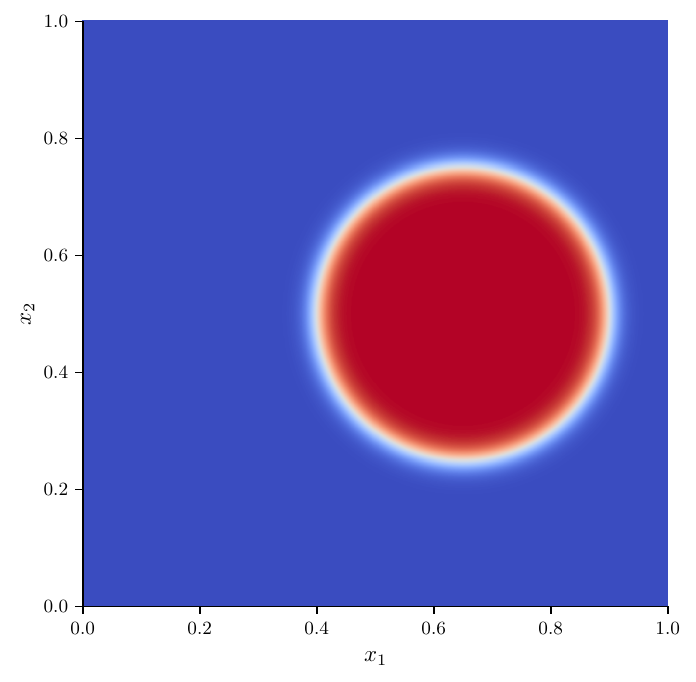}
                \end{subfigure}
                \hfill
                \begin{subfigure}{0.24\textwidth}
                    \centering
                    \includegraphics[width=\textwidth]{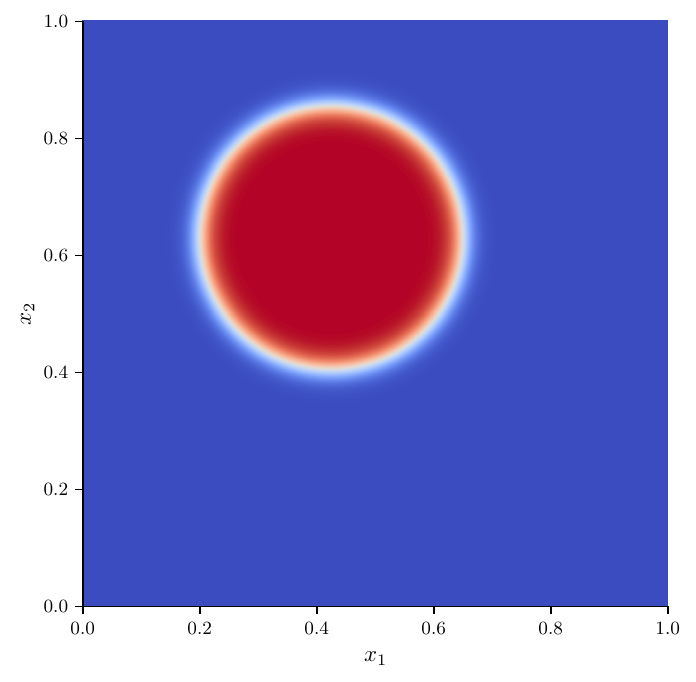}
                \end{subfigure}
                \hfill
                \begin{subfigure}{0.24\textwidth}
                    \centering
                    \includegraphics[width=\textwidth]{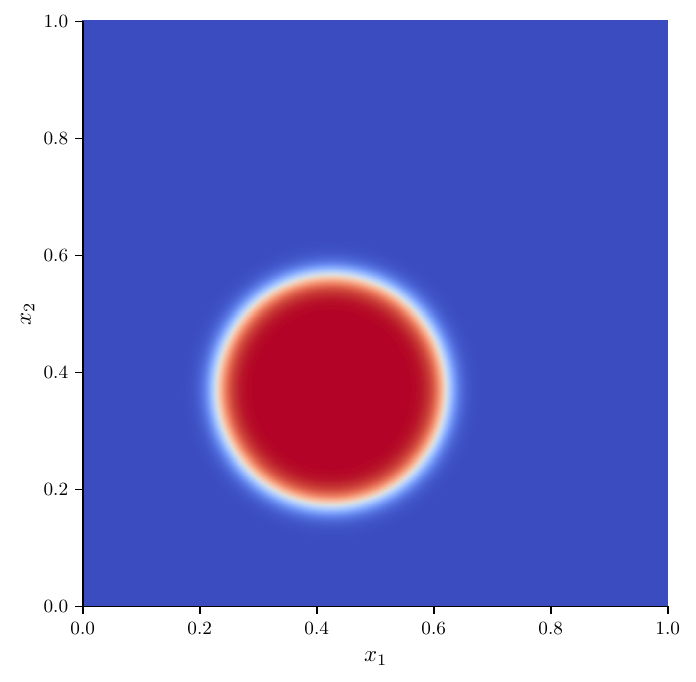}
                \end{subfigure}
                \hfill
                \begin{subfigure}{0.24\textwidth}
                    \centering
                    \includegraphics[width=\textwidth]{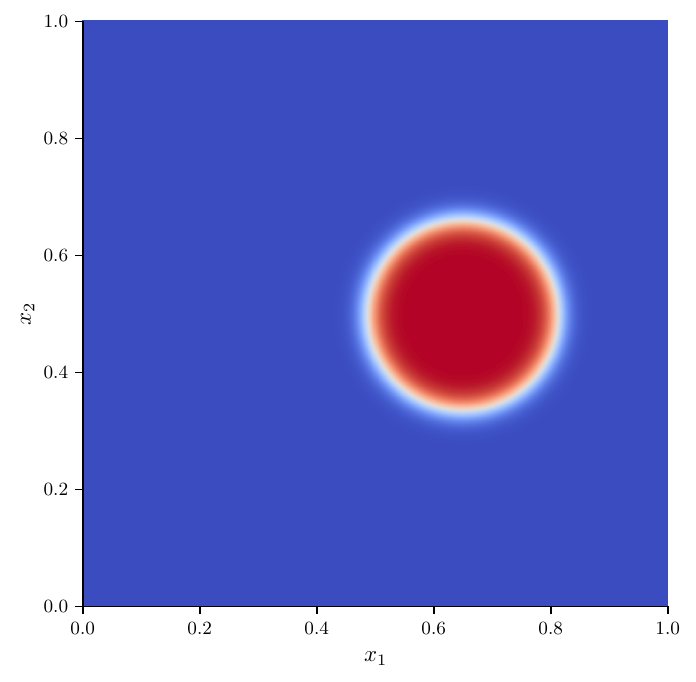}
                \end{subfigure}
        
                \vspace{1mm}
        
                \begin{subfigure}{0.24\textwidth}
                    \centering
                    \includegraphics[width=\textwidth]{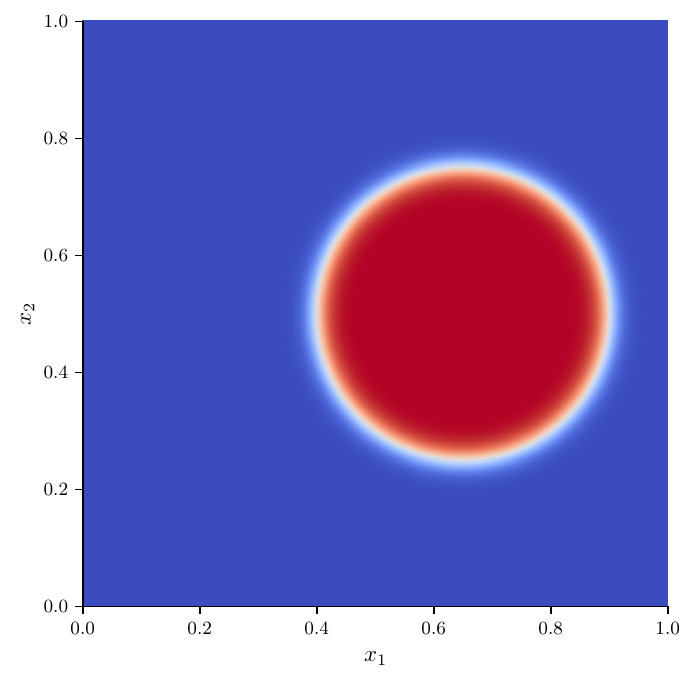}
                    \caption{$t = 0$}
                \end{subfigure}
                \hfill
                \begin{subfigure}{0.24\textwidth}
                    \centering
                    \includegraphics[width=\textwidth]{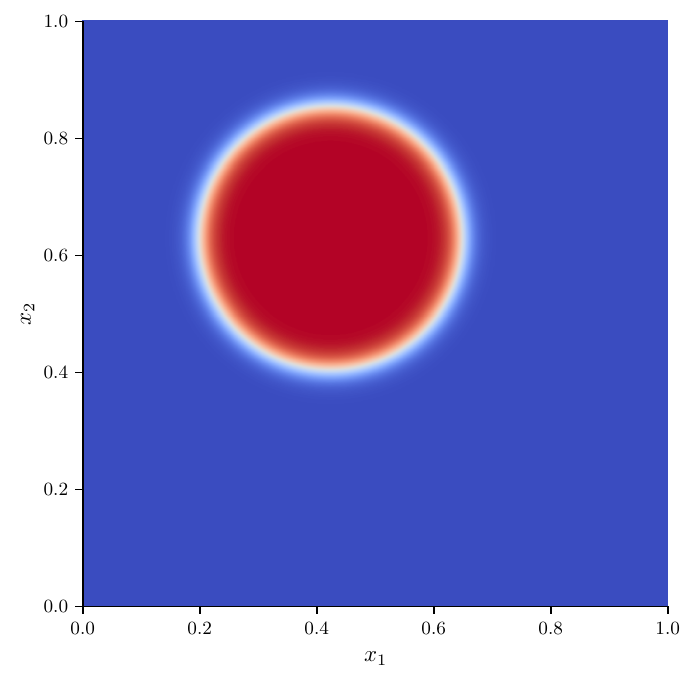}
                    \caption{$t = 0.006$}
                \end{subfigure}
                \hfill
                \begin{subfigure}{0.24\textwidth}
                    \centering
                    \includegraphics[width=\textwidth]{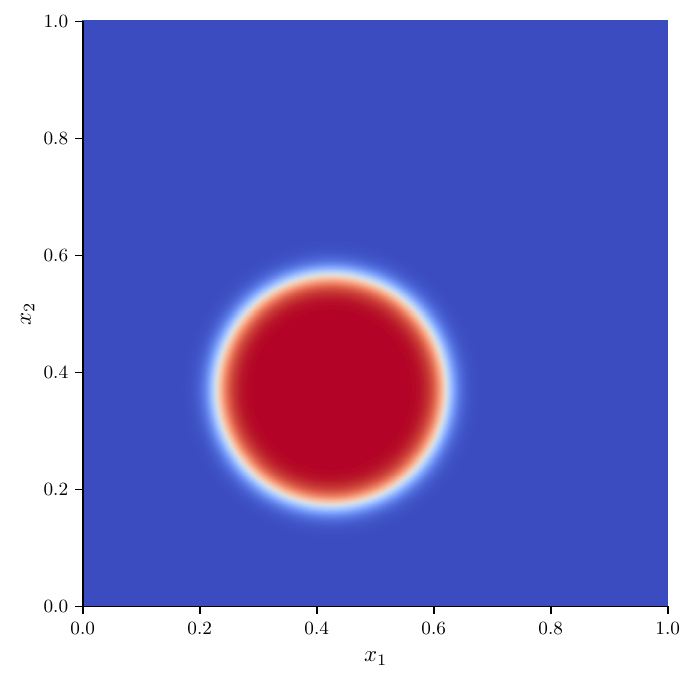}
                    \caption{$t = 0.012$}
                \end{subfigure}
                \hfill
                \begin{subfigure}{0.24\textwidth}
                    \centering
                    \includegraphics[width=\textwidth]{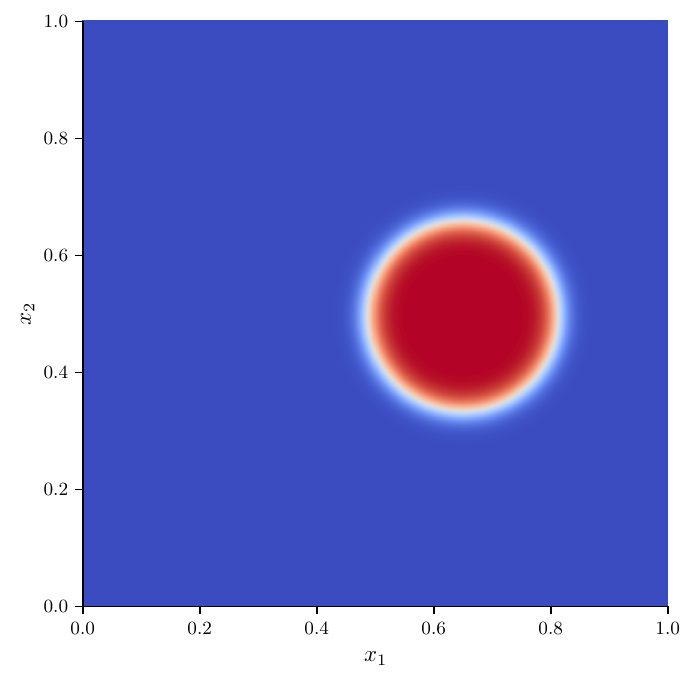}
                    \caption{$t = 0.018$}
                \end{subfigure}
            \end{minipage}
            \hfill
            \begin{minipage}{0.08\textwidth}
                \centering
                \includegraphics[width=\linewidth,height=0.22\textheight,keepaspectratio]{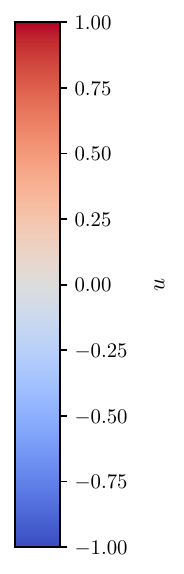}
                \vspace{6mm} 
            \end{minipage}
        
            \caption{Evolution snapshots of the $2$D advective Allen-Cahn simulation at selected physical times $t\in\{0,\,0.006,\,0.012,\,0.018\}$. The top row displays the phase evolution computed using the full-grid IMEX-LDG scheme, while the bottom row illustrates the corresponding results obtained via the sparse-grid combination technique.}
            \label{fig:ldg:experiments:ac:2d:dynamics}
        \end{figure}

        We first evaluate the circular shrinking process in the $2$D domain by setting the final simulation time to $T=0.018$ and the diffuse interface width parameter to $\epsilon=0.015$.
        To ensure that the phase region completes exactly one full revolution over the simulation period, the angular velocity is fixed at $\omega=2\pi/T$.
        Under this specific configuration, the theoretical limit for the final radius is $R_{\mathrm{ana}}(T)=0.1627$.
        We compare the performance of the full-grid IMEX-LDG scheme against its sparse-grid combination technique counterpart, employing quadratic polynomial elements ($\nu=2$) for both approaches.
        The full-grid resolution is set to $\bm{\alpha}=512$, whereas the combination technique parameters are chosen as $N_B=16$ and $N_L=5$, ensuring an equivalent effective resolution across both frameworks. 
        It is important for the combination technique to succeed that the chosen base refinement level must provide a sufficiently dense local grid to adequately resolve the steep gradients of the transition region.
        If the baseline interpolation nodes fail to capture the underlying physics of the diffuse interface, the sparse-grid extrapolation will be unable to preserve the structural integrity of the moving boundary.

        We maintain the linear solver configuration employed for the previous periodic benchmarks.
        Specifically, we explicitly assemble the global system matrix and utilize an Additive Schwarz preconditioner coupled with an ILU local sub-solver.
        For both the full-grid and the combination technique scenarios, the computational workload is distributed across $54$ MPI processes.
        It is important to note that the Allen-Cahn experiments impose stringent restrictions on the selection of the time step.
        In addition to the standard CFL condition \eqref{eq:ldg:imex:cfl}, the explicit temporal treatment of the stiff nonlinear reaction term introduces another stability constraint.
        As established in \cite{Shen2010}, this requires the time step to satisfy the bound $\Delta t \leq \mathcal{O}(\epsilon^2)$.
        However, in our specific experimental configuration, the spatial grid resolution is sufficiently fine that the standard CFL restriction is more restrictive and dominates stability criteria, thereby naturally satisfying the reaction-driven bound.
        Consequently, establishing a CFL number of $C=0.2$ yields a highly constrained time step of $\Delta t=3.165\times 10^{-7}$.
        
        \begin{figure}[htbp]
            \centering
            \begin{subfigure}{0.48\textwidth}
                \centering
                \includegraphics[width=\textwidth]{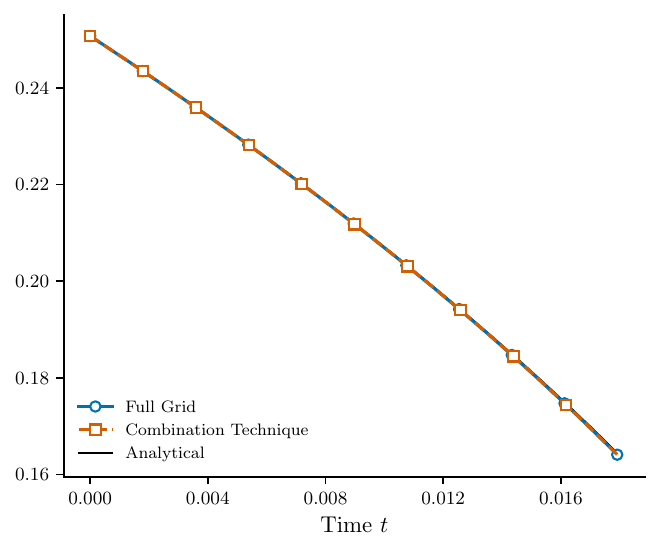}
                \caption{Radius evolution $R(t)$.}
                \label{fig:ldg:experiments:ac:2d:radius:val}
            \end{subfigure}
            \hfill 
            \begin{subfigure}{0.48\textwidth}
                \centering
                \includegraphics[width=\textwidth]{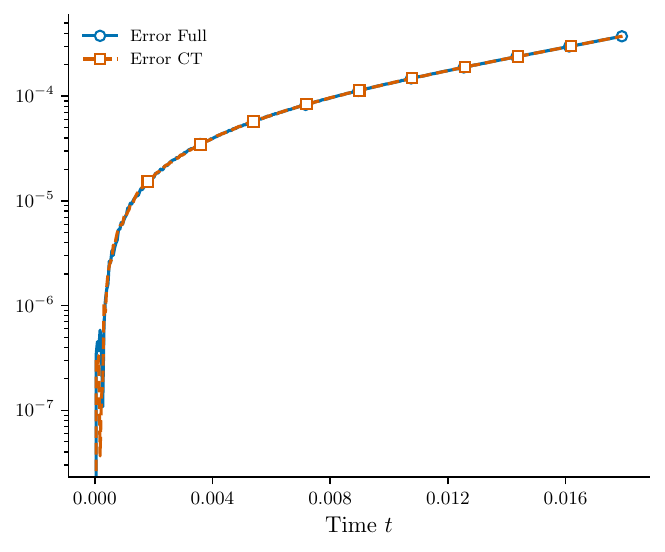}
                \caption{Absolute error $|R_{\mathrm{sim}} - R_{\mathrm{ana}}|$.}
                \label{fig:ldg:experiments:ac:2d:radius:error}
            \end{subfigure}
        
            \caption{Temporal evolution of the effective phase radius (left) and the corresponding absolute error relative to the theoretical sharp-interface limit (right) for the $2$D advective Allen-Cahn problem.}
            \label{fig:ldg:experiments:ac:2d:radius}
        \end{figure}
        
        The numerical results obtained at the final simulation time are virtually identical for both methodologies.
        We observe a decay in the total phase mass from an initial value of $0.196350$ to a final value of $0.084026$, representing a substantial mass reduction of $57.206\%$.
        Correspondingly, both experiments yield a final effective radius of $0.163543$.
        When compared against the theoretical limit radius, this translates to a highly accurate relative error of approximately $0.49\%$, demonstrating that both schemes successfully capture the sharp-interface dynamics.

        The critical distinction between the two approaches lies in their computational efficiency.
        While the full-grid simulation required $7.01$ hours to complete, the sparse-grid combination technique required only $2.18$ hours.
        This demonstrates that the sparse-grid methodology is more than three times faster, yielding important computational savings without compromising the geometric integrity or the accuracy of the physical metrics.

        In order to visualize the temporal evolution of the phase separation, Figure \ref{fig:ldg:experiments:ac:2d:dynamics} presents the computed phase fields at selected physical times $t\in\{0, 0.006, 0.012, 0.018\}$, for both methodologies. 
        Visual inspection confirms that the moving boundaries are virtually indistinguishable during all stages of the simulation.

        Finally, Figure \ref{fig:ldg:experiments:ac:2d:radius} compares the temporal evolution of the effective phase radius, computed via both the full-grid and the combination technique schemes, against the theoretical sharp-interface limit, alongside their corresponding absolute errors.
        The temporal evolution of the absolute errors reveals three distinct regimes.
        Initially, \joseanto{high-frequency} oscillations occur as the continuous initialization relaxes onto the discrete grid.
        This is followed by a rapid error growth during the onset of the interface motion, representing the accumulation of spatial and temporal tracking errors.
        Finally, the error exhibits an asymptotic plateau.
        This bounded error indicates that the numerical front attains the correct curvature-driven velocity, settling into a steady-state phase lag relative to the sharp-interface analytical limit.
        Notably, the combination technique preserves this bounded behavior without introducing long-term drift when compared with the full-grid scheme.

        For the $3$D spherical shrinking benchmark, we keep the diffuse interface width at $\epsilon=0.015$ and set the final simulation time to $T=0.005$.
        Under these conditions, the theoretical sharp-interface limit radius is $R_{\mathrm{ana}}(T)=0.2061$.
        We set the angular velocity parameter to $\omega=\pi/T$, meaning that the sphere sweeps a total rotation angle of $180$ degrees over the course of the simulation.

        In order to resolve the domain, we establish a full-grid resolution of $\bm{\alpha}=160$, coupled with quadratic basis functions ($\nu=2$).
        To achieve an equivalent effective resolution with the combination technique, we set the base refinement to $N_b=10$ and the maximum refinement level to $N_L=4$.
        As established in the $2$D case, this base refinement proves sufficient to adequately capture the steep gradients of the diffuse interface.
        Given the spatial discretization, the standard advective CFL time-step restriction is strictly more severe than the stability bound imposed by the explicit treatment of the stiff reaction term.
        Consequently, applying a CFL coefficient of $C=0.8$ yields a highly constrained integration step of $\Delta t = 2.251\times 10^{-6}$.
        In addition, we retain the linear solver configuration employed in the two-dimensional scenario, distributing the computational workload across $128$ MPI processes.
        
        \begin{figure}[htbp]
            \centering
                \begin{subfigure}{0.30\textwidth}
                    \centering
                    \includegraphics[trim={10cm 10cm 10cm 10cm}, clip ,width=\textwidth]{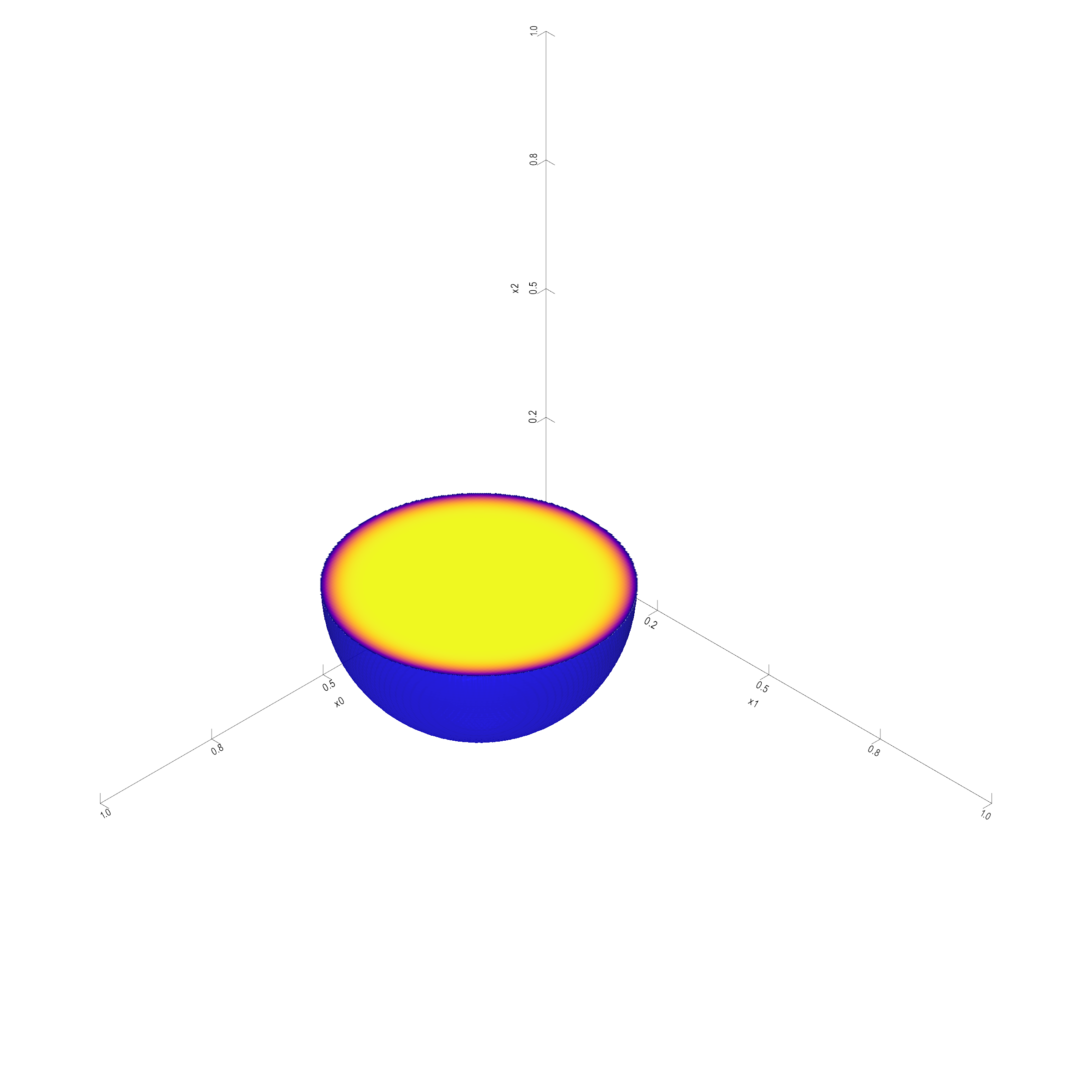}
                \end{subfigure}
                \hfill
                \begin{subfigure}{0.30\textwidth}
                    \centering
                    \includegraphics[trim={10cm 10cm 10cm 10cm}, clip ,width=\textwidth]{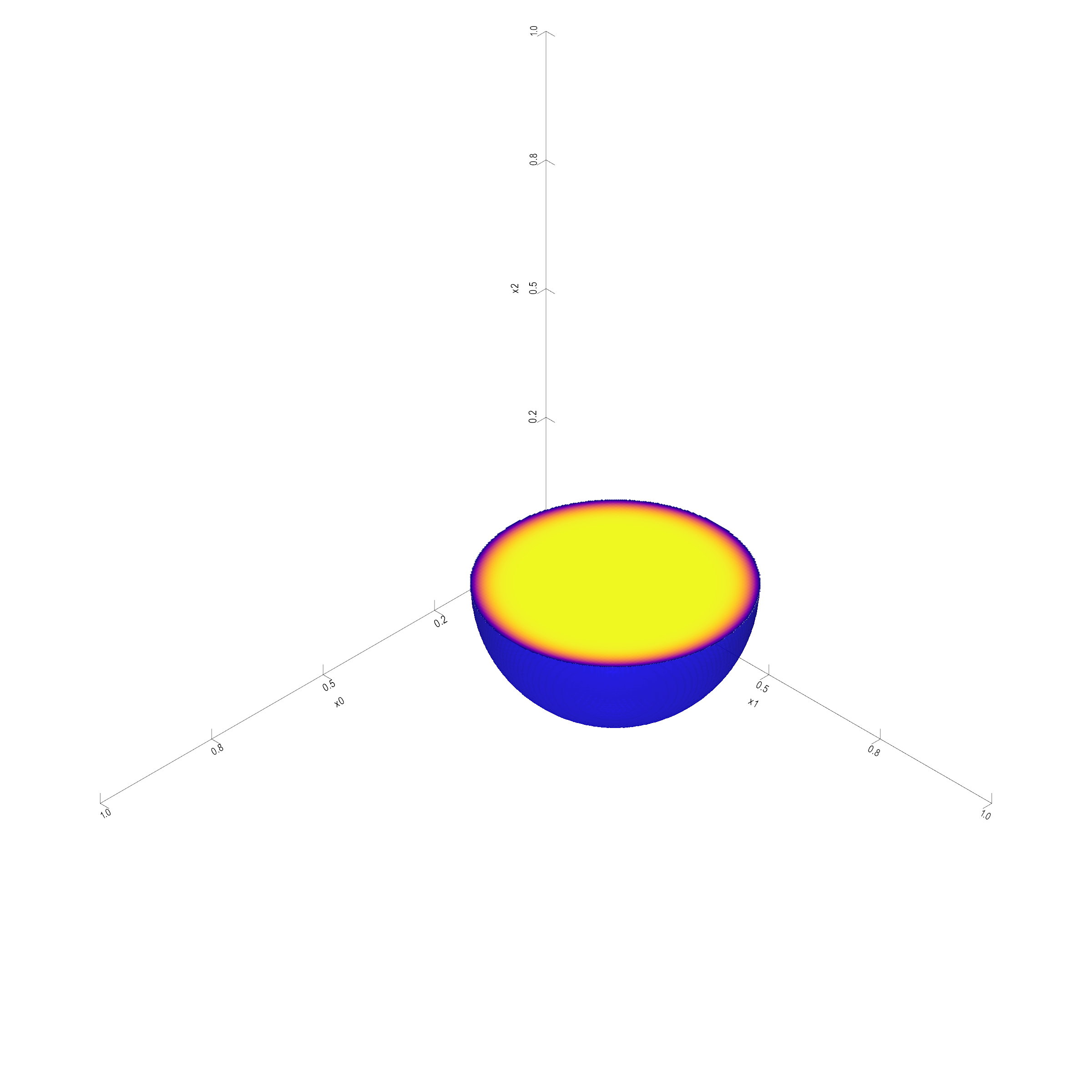}
                \end{subfigure}
                \hfill
                \begin{subfigure}{0.30\textwidth}
                    \centering
                    \includegraphics[trim={10cm 10cm 10cm 10cm}, clip ,width=\textwidth]{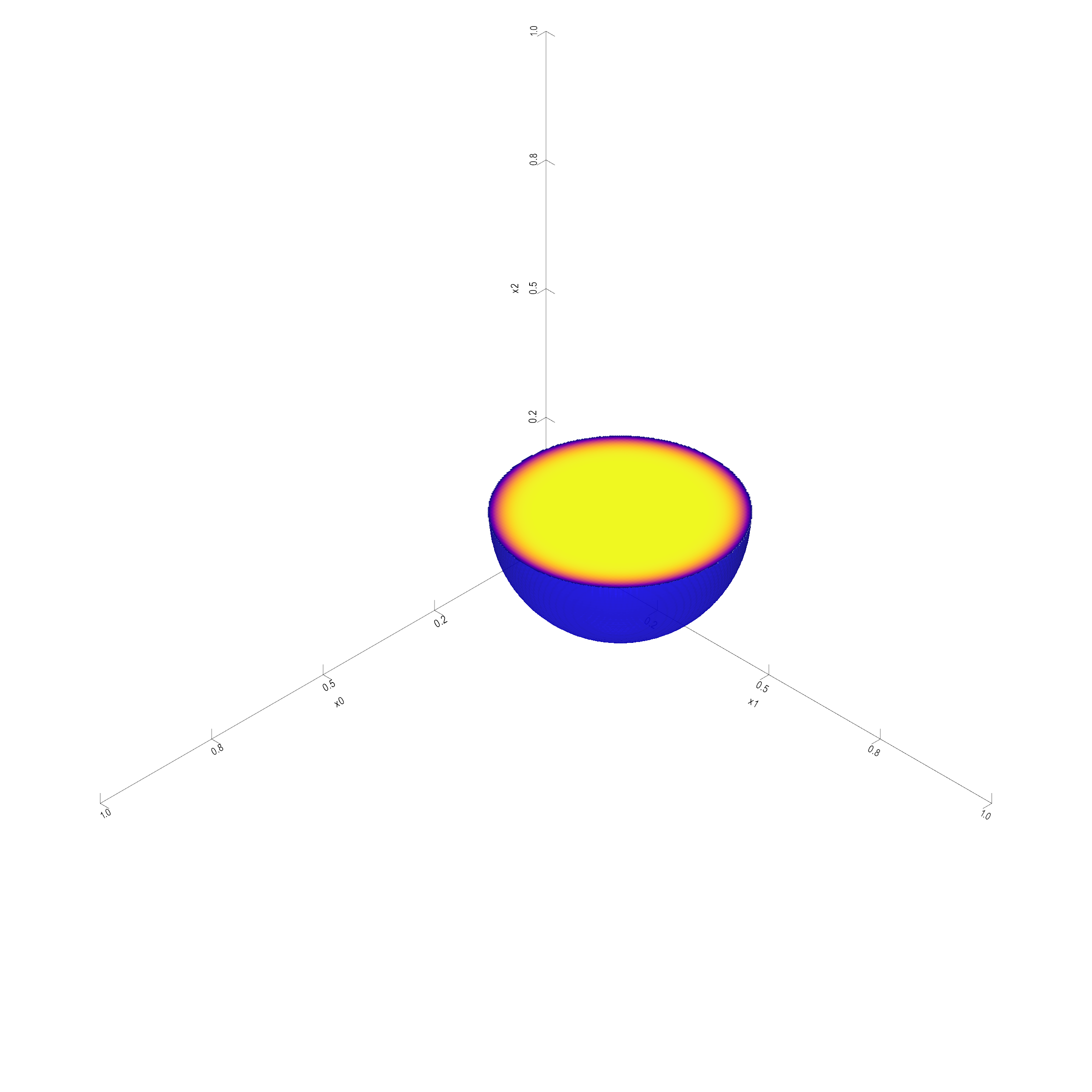}
                \end{subfigure}
        
                \vspace{1mm}
        
                \begin{subfigure}{0.30\textwidth}
                    \centering
                    \includegraphics[trim={10cm 10cm 10cm 10cm}, clip ,width=\textwidth]{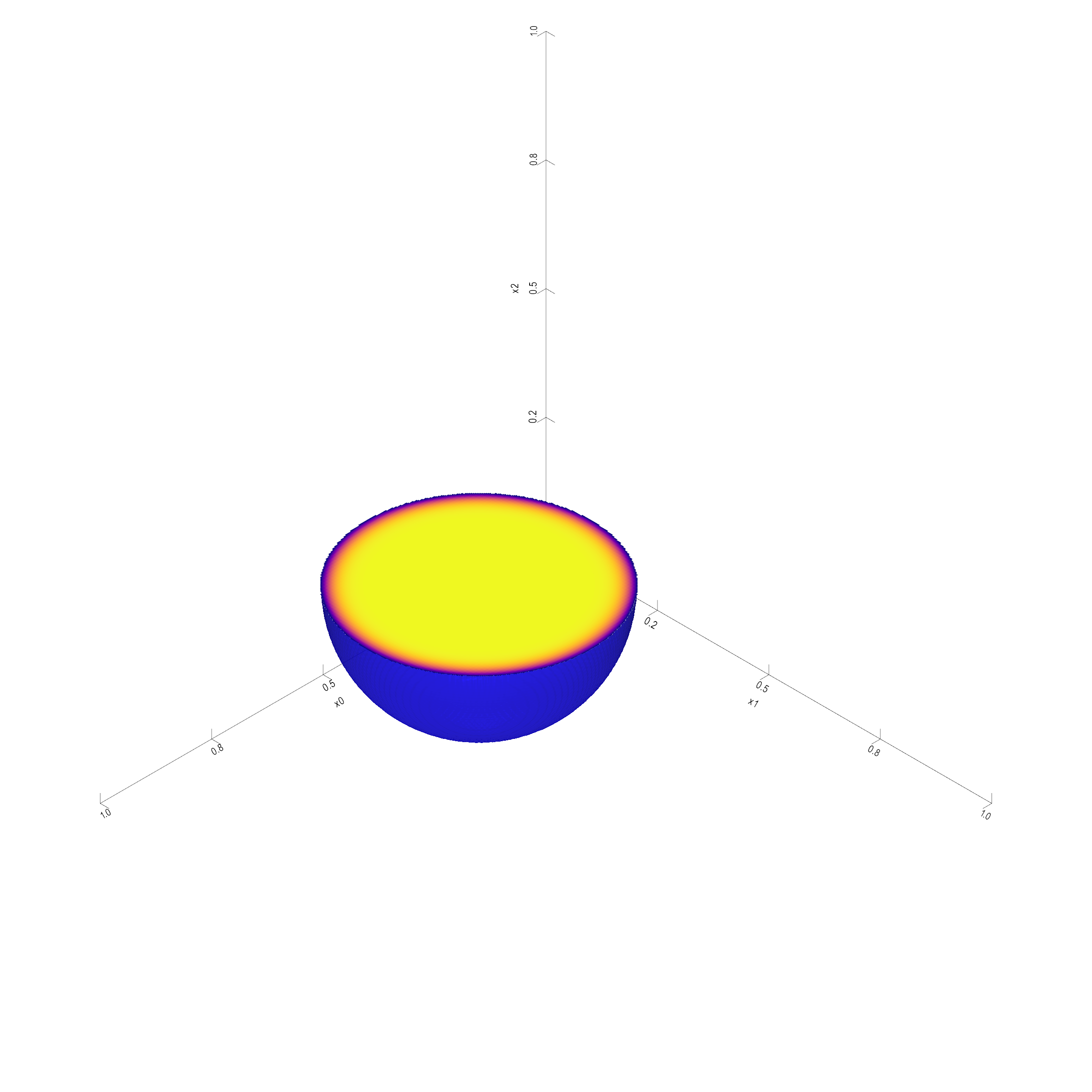}
                    \caption{$t = 0$}
                \end{subfigure}
                \hfill
                \begin{subfigure}{0.30\textwidth}
                    \centering
                    \includegraphics[trim={10cm 10cm 10cm 10cm}, clip ,width=\textwidth]{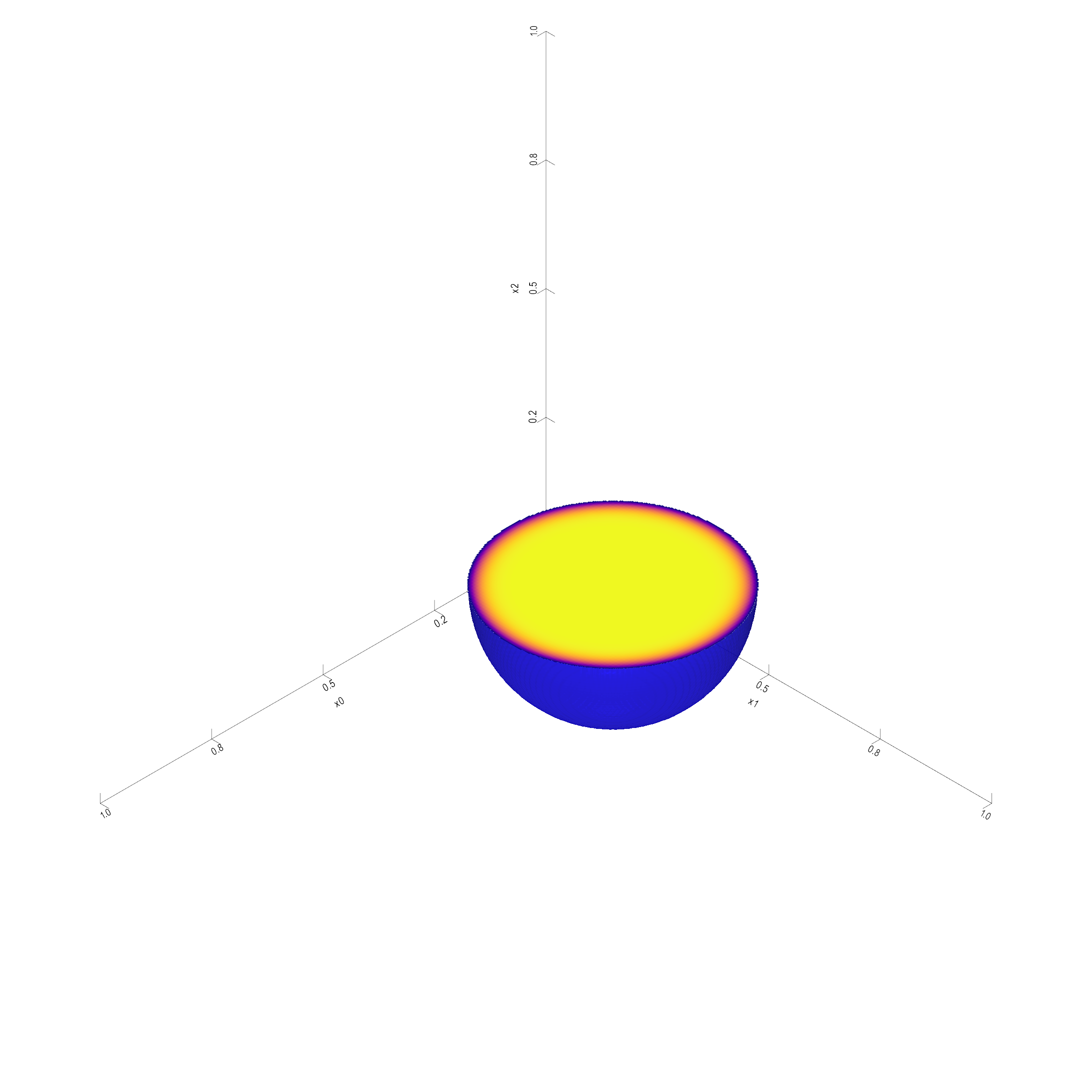}
                    \caption{$t = 0.0025$}
                \end{subfigure}
                \hfill
                \begin{subfigure}{0.30\textwidth}
                    \centering
                    \includegraphics[trim={10cm 10cm 10cm 10cm}, clip ,width=\textwidth]{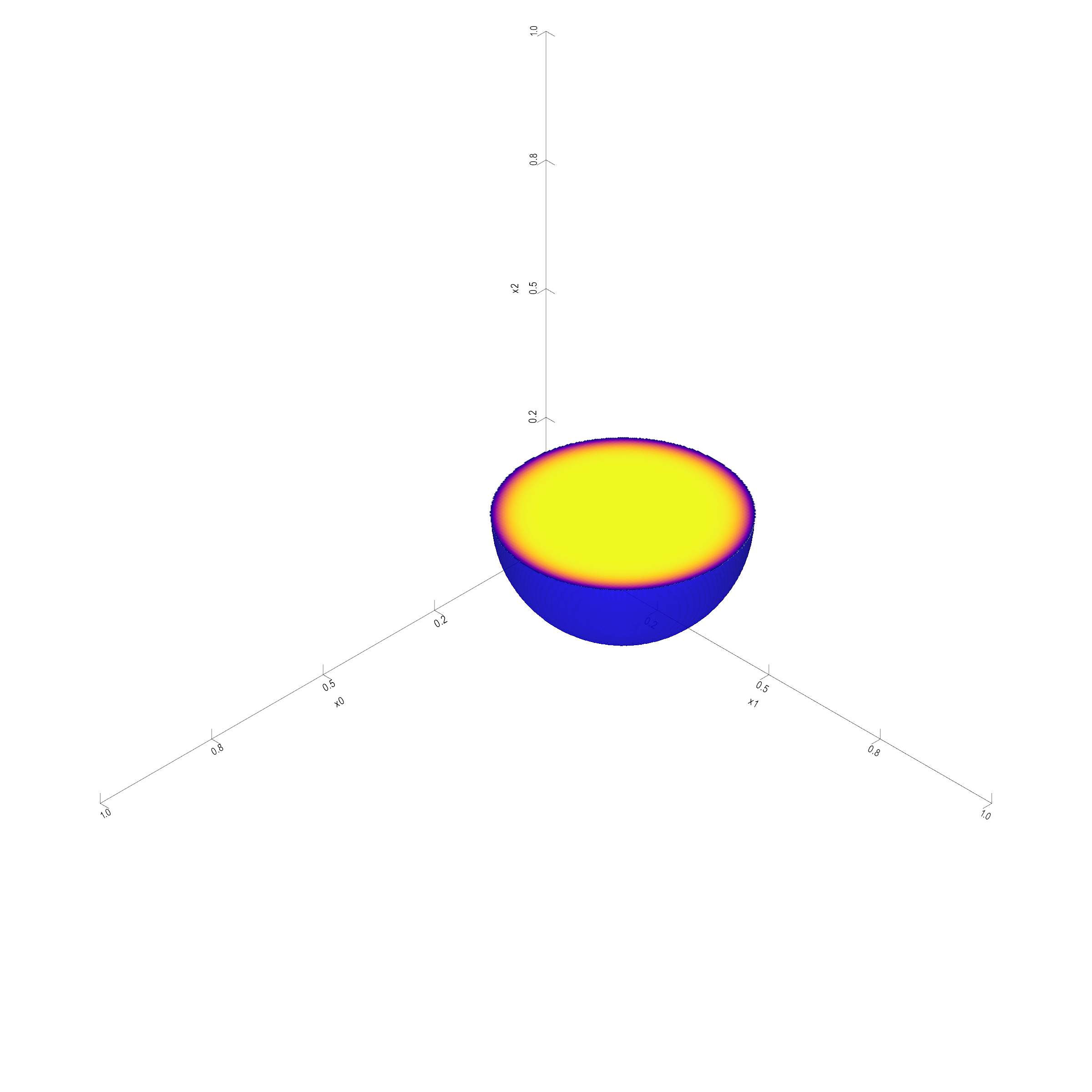}
                    \caption{$t = 0.005$}
                \end{subfigure}
        
            \caption{Evolution snapshots of the $3$D advective Allen-Cahn simulation at selected physical times $t\in\{0,\,0.0025,\,0.005\}$. The top row displays the phase evolution computed using the full-grid IMEX-LDG scheme, while the bottom row illustrates the corresponding results obtained via the sparse-grid combination technique. To visualize the internal dynamics, the upper hemisphere of the domain is clipped. The volumetric color mapping represents the phase value $u$, where the blue contour delineates the outer phase boundary ($u=0$), transitioning through red across the diffuse interface, and culminating in yellow for the pure internal phase ($u=1$).}
            \label{fig:ldg:experiments:ac:3d:dynamics}
        \end{figure}

        The final \joseanto{simulation data} yield conclusions highly consistent with the $2$D benchmark, but with an even more pronounced computational advantage.
        Both the full-grid and the sparse-grid scenarios exhibit a $42.55\%$ mass reduction over the course of the simulation, decaying from an initial phase volume of $0.0654$ to $0.0376$.
        The effective radii computed by the two distinct numerical schemes present a negligible relative difference of only $9.62\times 10^{-6}$, confirming that the physical outputs are identical.
        When compared to the theoretical sharp-interface limit, both methodologies demonstrate a relative error of just $0.83\%$.

        Once again, the critical distinction between the methodologies lies in their computational cost, a gap that widens significantly in three dimensions.
        While the standard full-grid approach required $11$ hours to complete the simulation, the combination technique required only $3.24$ hours. This represents a remarkable speedup, proving that the combination technique alleviates the curse of dimensionality while preserving the physical accuracy of the Allen-Cahn dynamics.

        Figure \ref{fig:ldg:experiments:ac:3d:dynamics} presents three snapshots illustrating the temporal evolution of the phase separation at physical times $t\in\{0.0,\,0.0025,\,0.005\}$ for both methodologies. 
        To visualize the three-dimensional dynamics, we generate a volumetric rendering where the upper hemisphere is clipped, and negative phase values are filtered out.
        This reveals the solid internal structure of the shrinking region while exposing the thickness of the diffuse interface.
        As anticipated from the quantitative analysis, the spatial profiles generated by the full-grid and sparse-grid approaches are visually indistinguishable.
        In both renderings, the interface is cleanly captured, appearing as \joseanto{a distinct red circular annulus} tracking the exterior contour of the phase boundary.
        
        \begin{figure}[htbp]
            \centering
            \begin{subfigure}{0.48\textwidth}
                \centering
                \includegraphics[width=\textwidth]{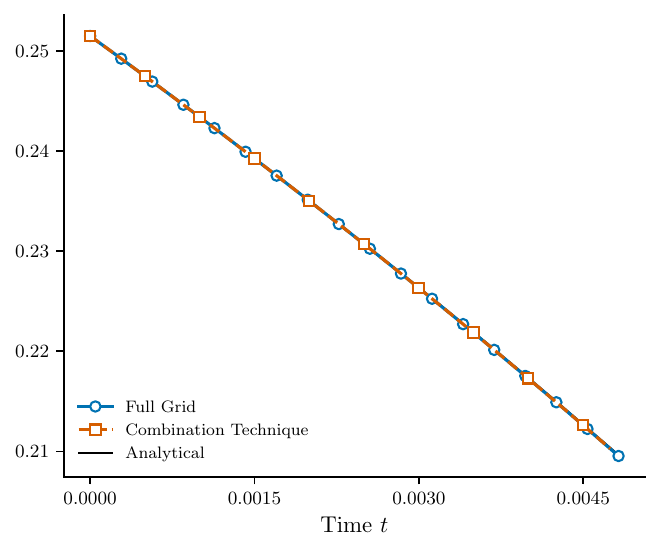}
                \caption{Radius evolution $R(t)$.}
                \label{fig:ldg:experiments:ac:3d:radius:val}
            \end{subfigure}
            \hfill 
            \begin{subfigure}{0.48\textwidth}
                \centering
                \includegraphics[width=\textwidth]{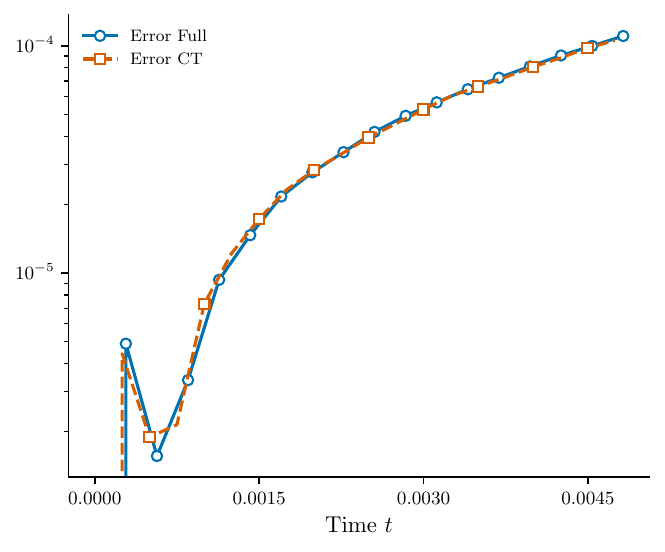}
                \caption{Absolute error $|R_{\mathrm{sim}} - R_{\mathrm{ana}}|$.}
                \label{fig:ldg:experiments:ac:3d:radius:error}
            \end{subfigure}
        
            \caption{Temporal evolution of the effective phase radius (left) and the corresponding absolute error relative to the theoretical sharp-interface limit (right) for the $3$D advective Allen-Cahn problem.}
            \label{fig:ldg:experiments:ac:3d:radius}
        \end{figure}

        Figure \ref{fig:ldg:experiments:ac:3d:radius} presents the simulated effective phase radii alongside the theoretical sharp-interface limit.
        Consistent with the $2$D case, the evolutionary trajectories by both numerical methodologies are almost identical.
        Because the total simulated time frame is shorter than in the previous scenario, the initial transient \joseanto{regime} is captured with higher resolution, providing a clearer visualization of the initial oscillation, which \joseanto{is a direct product} of the continuous analytical initial condition relaxing onto the discrete computational grid.
        Following this brief accommodation \joseanto{phase}, the dynamics seamlessly transition into the subsequent error regimes identified in the $2$D analysis.

    \subsection{Black-Scholes model under credit risk adjustments} \label{ldg:experiments:bs}
        The remainder of this section is devoted to demonstrating the capabilities of the proposed IMEX-LDG sparse-grid combination technique for pricing a financial instrument subject to counterparty credit risk, i.e., taking into account the possibility of default by either party involved in the contract.
        This financial application serves an ideal benchmark, as it perfectly aligns with the core computational challenges addressed throughout the work.
        First, the multidimensional Black-Scholes framework yields a scalar convection-diffusion PDE that can be cast into the conservative formulation \eqref{eq:ldg:pde}.
        Second, the spatial dimensions of the computational domain correspond to the number of financial assets composing the derivative.
        Because portfolios usually feature derivatives depending on baskets of multiple assets, this context generates high-dimensional problems.

        Specifically, we focus on the valuation of a European basket put option.
        This contract grants the holder the right to sell a weighted portfolio of underlying stocks at a predetermined strike price $K_{\text{opt}}$ upon a specific maturity date $T$.
        For a comprehensive theoretical background on multidimensional derivative pricing and credit risk valuation adjustments (XVA), we refer the reader to, for example, \cite{hull2012,Kjaer2011}.

        Given the deterministic payoff function, the present value of the option is governed by the classic multidimensional Black-Scholes PDE equation.
        Formulated as a forward-in-time problem, the governing equation is given by
        \begin{equation}
            \label{eq:ldg:experiments:bs:eq_non_cons}
            \partial_t u-\dfrac{1}{2}\sum_{i=1}^dr_ix_i\partial_{x_i}u = \dfrac{1}{2}\sum_{i,j=1}^d\rho_{ij}\sigma_i\sigma_jx_ix_j\partial_{x_i x_j}u-ru-g_{\text{xva}}(u), 
        \end{equation}
        where $u(t,\xv)$ is the option value, $x_i$ is the price of the $i$-th asset, $\sigma_i$ is its constant volatility, $\rho_{ij}$ is the correlation coefficient, and $r_i$ accounts for the spread between the market risk-free rate and the dividend yield of the $i$-th asset.
        The nonlinear reaction term $g_{\text{xva}}(u)$ incorporates the stiff XVA adjustments into the valuation.
        Specifically, this is defined as,
        \begin{equation}
            \label{eq:ldg:experiments:bs:risk}
            g_{\text{xva}}(u) = (1-R_B)\lambda_B\min\{u,0\}+(1-R_C)\lambda_C\max\{u,0\}+s_F\max\{u,0\}, 
        \end{equation}
        where $\lambda_B,\,\lambda_C>0$ denote the probabilities of default of the involved parties, $R_B,R_C\in [0,1]$ represent their respective recovery rates in the event of default, and $s_F$ is a funding spread.

        To cast equation \eqref{eq:ldg:experiments:bs:eq_non_cons} into the conservative formulation required by the LDG spatial discretization \eqref{eq:ldg:pde}, we apply a sequence of domain transformations.
        \joseanto{First, we} introduce the logarithmic spatial transformation $y_i=\log(x_i)$.
        This maps the semi-infinite domain $\R_+^d$ to $\R^d$ and removes the spatial dependencies from the differential coefficients.
        This transformation not only \joseanto{normalizes} the domain but also improves the conditioning of the resulting matrix system, thereby significantly enhancing the efficiency of the linear solver.
        
        Subsequently, we apply an orthogonal coordinate rotation $\bm{z}=Q^T\bm{y}$.
        The orthogonal matrix $Q$ is chosen to diagonalize the covariance matrix, effectively decoupling the principal diffusion directions and eliminating all mixed-cross derivative terms.
        This rotational step is crucial when tackling high-dimensional problems, as it circumvents the substantial computational overhead associated with evaluating cross-derivatives within the LDG scheme.
        Furthermore, the elimination of these mixed terms is advantageous for the sparse-grid combination technique, sensitive to cross-directional dependencies.

        Under these transformations, the convection-diffusion-reaction problem naturally adopts the conservative form:
        \begin{equation}
            \label{eq:ldg:experiments:bs:eq}
            \partial_t u + \nabla_{\bm{z}}\cdot\left(\bm{c}_v\, u\right) = \nabla_{\bm{z}}\cdot\left(D\,\nabla_{\bm{z}}u\right) - ru-g_{\text{xva}}(u),\quad \text{in } (0,T)\times\Omega_{\bm{z}},
        \end{equation}
        where $u(t,\bm{z})$ is the transformed option value.
        The orthogonal transformation guarantees that the diffusion matrix $D$ is diagonal, although anisotropic, while the advective velocity vector $\bm{c}_v$ incorporates the drift corrections.
        For the sake of clarity and financial interpretability, the remainder \joseanto{of the exposition} is presented in terms of the original asset price variables $\xv$.
        However, it is implicitly understood that all numerical computations are performed in the fully transformed, diagonalized domain $\bm{z}$.

        The terminal payoff of the financial instrument serves as the initial condition for the forward-in-time Black-Scholes PDE.
        Specifically, we consider an arithmetic average basket put option payoff, defined as:
        \begin{equation}
            \label{eq:ldg:experiments:bs:ic}
            u_0(\xv) = \max\left\{K_{\text{opt}} - \dfrac{1}{d}\sum_{i=1}^dx_i, 0\right\}.
        \end{equation}
        Note that this payoff function introduces a non-differentiable hyperplane (`kink' in financial jargon) into the domain. 
        This lack of global smoothness poses \joseanto{a numerical} challenge for the combination technique, as it fundamentally violates the bounded mixed-derivative assumption required for optimal sparse-grid convergence.
        Naturally, this initial condition is properly mapped onto the computational grid using the aforementioned transformations.

        To complete the formulation of the initial-boundary value problem, we prescribe numerical outflow boundary conditions across the entire boundary of the computational domain, which are seamlessly enforced \joseanto{through} the LDG spatial discretization.
        This approach is mathematically justified by the asymptotic behavior of the option contract.
        Because the terminal payoff is constructed from a combination of linear hyperplanes, the option value becomes asymptotically linear in the far-field limit, rendering the diffusive fluxes negligible and allowing the solution to propagate freely outward.

        For our numerical experiment, we consider a European basket put option comprising $d=5$ underlying assets.
        The contract is characterized by a maturity of $T=1$ year and a strike price of $K_{\text{opt}}=15$, operating within a market featuring a risk-free interest rate of $2\%$.
        For the sake of simplicity, we assume a symmetric asset configuration where all underlying stocks exhibit an identical volatility of $\sigma_i=0.25$, and yield no dividends, leading to $r_i=0.02$.
        Also, the correlation between all assets is identical, $\rho_{ij}=0.5,\,i\neq j$.
        To incorporate the XVA adjustments, we utilize the default intensities $\lambda_B=0.02$ and $\lambda_C=0.05$, corresponding recovery rates $R_B=R_C=0.4$, and a funding spread defined as $s_F=\lambda_B(1-R_B)$.
        
        From a computational standpoint, the semi-infinite pricing domain must be truncated to safely impose the outflow boundary conditions.
        For this specific problem, we restrict the original price domain to the hypercube $\Omega=(1,45)^5$. 
        Recall that the numerical solver operates under the logarithmic and rotational coordinate transformations.
        To guarantee that the rotational mapping does not clip the region of interest, the transformed computational domain $\Omega_{\bm{z}}$ is defined as the minimal bounding box completely enclosing the rotated space.
        The spatial discretization and integration procedures are strictly performed within this extended $\Omega_{\bm{z}}$ domain.

        The combination technique is configured with a base refinement of $N_b=6$ and a maximum refinement level of $N_L=6$, utilizing quadratic basis functions $\nu=2$.
        This setup achieves an effective target full-grid resolution of $\bm{\alpha}=384$.
        Under this setting, \joseanto{computing} the solution involves the numerical resolution of $456$ sub-grids.
        Temporal integration is governed by the CFL condition \eqref{eq:ldg:imex:cfl} with \joseanto{a CFL} coefficient of $C=0.75$, yielding  a time step of $\Delta t=0.0853$.
        At each time step, the linear systems are solved using a matrix-free approach without preconditioning, governed by a relative residual stopping criterion of $5\times 10^{-6}$.
        In addition, the workload was distributed across $1024$ MPI processes.
         
        In order to validate the accuracy of the proposed methodology, we benchmark the numerical results against a reference solution.
        First, the risk-free price of the basket options is computed via standard Monte Carlo integration.
        Under the present configuration, an exact analytic mapping can be applied to derive the corresponding price under the XVA adjustments directly from this risk-free baseline.
        We restrict our comparative analysis to a structured $2$D slice of the domain, specifically mapping the $(x_1,x_2)$-plane while clamping the remaining asset prices at the strike region ($x_3=x_4=x_5=15$).
        The Monte Carlo reference solution is generated with a sufficient number of paths to guarantee an accuracy of $3\times10^{-4}$ at $95\%$ confidence level.
        
        \begin{figure}[htbp]
             \centering
             \begin{subfigure}[b]{0.48\textwidth}
                \centering
                 \includegraphics[width=\textwidth]{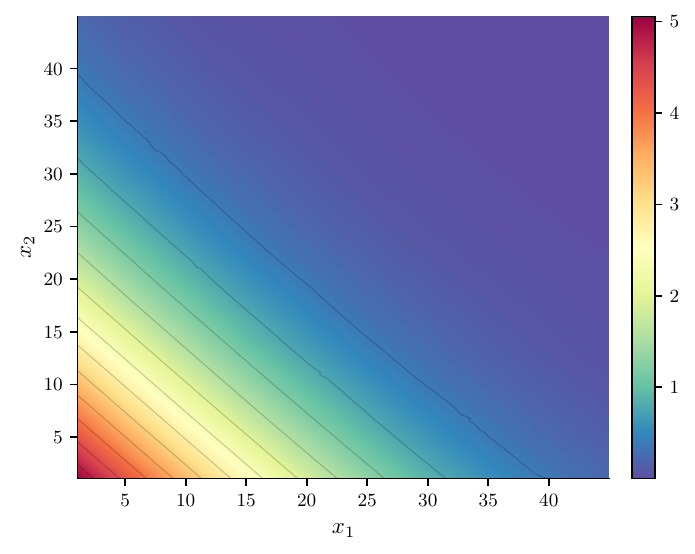}
                 \caption{Numerical solution.}
                  \label{fig:ldg:experiments:bs:5d:xva:sol} 
             \end{subfigure}
             \hfill
             \begin{subfigure}{0.48\textwidth}
                \centering
                \includegraphics[width=\textwidth]{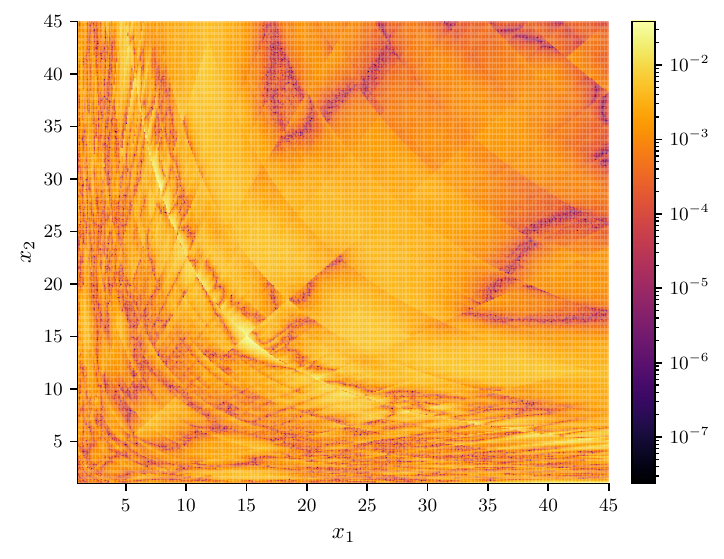}
                 \caption{Absolute error.}
                  \label{fig:ldg:experiments:bs:5d:xva:err} 
             \end{subfigure}
             \caption{Two-dimensional slice of the numerical solution (left) and absolute error (right) for the $5$D arithmetic basket put option with XVA adjustments. The multidimensional domain is intersected at the strike level of the non-varying dimensions, $x_3=x_4=x_5=15$.}
             \label{fig:ldg:experiments:bs:5d:xva}
        \end{figure}

        The numerical results demonstrate the practical viability of the proposed methodology.
        The total simulation time amounted to $9.3$ hours.
        From an industrial perspective, this runtime confirms that the $5$D pricing problem can be resolved within \joseanto{a standard} overnight batch processing window.

        To quantitatively assess the accuracy of the combination technique, we evaluate the discrete $L^p$ error norms as defined in \eqref{eq:ldg:experiments:discrete_error} against the Monte Carlo reference.
        The computed metrics yield $\text{E}_1=1.899\times 10^{-3}$, $\text{E}_2=2.655\times 10^{-3}$ and the maximum error $\text{E}_{\infty}=3.837\times 10^{-2}$.
        These global norms confirm a highly robust overall approximation, with the mean absolute and mean \joseanto{square} errors consistently residing in the $10^{-3}$ regime.
        However, the maximum \joseanto{pointwise} error exhibits a localized peak in the order of $10^{-2}$.
        The spatial distribution and underlying cause of this peak can be comprehensively analyzed \joseanto{through} the cross-sectional colormaps of the numerical solution and absolute error presented in Figure \ref{fig:ldg:experiments:bs:5d:xva}.

        In the numerical solution plot, the superimposed level curves indicate a highly resolved pricing surface, maintaining strict orthogonality with respect to the principal diagonal direction.
        The absolute error profile, computed against the Monte Carlo reference, reveals a curious geometric artifact.
        The characteristic sparse-grid error pattern, previously identified in the $5$D linear and viscous Burgers' benchmarks, reappears here, but it is distinctly deformed by the logarithmic and rotational coordinate mappings applied to the domain.

        The error map visually confirms that the peak numerical error discussed previously, the bright regions reaching the $10^{-2}$ magnitude, is heavily localized at the exact central intersection $x_1=x_2=15$.
        This spatial distribution empirically demonstrates how the non-differentiable kink at the strike price locally degrades the precision of the approximation.
        However, safely away from this singular region, the method rapidly recovers excellent precision.
        
\section{Conclusions} \label{ldg:conclusions}
    In summary, the present work has developed and analyzed general high-order IMEX-LDG numerical schemes for nonlinear parabolic PDEs in high dimensions.
    Regarding the spatial semi-discretization, the LDG formulation, combined with upwind fluxes, provides robust upwinding of convective terms, while preserving stability and high-order accuracy.
    In the temporal domain, the IMEX-RK strategy separates the stiff diffusive contributions, which are treated implicitly, from the explicit evaluation of the convective and source terms. 
    This splitting provides a twofold advantage: it bypasses the severe time-step restrictions induced by the diffusive terms, and circumvents the need to solve computationally expensive nonlinear systems at each time stage. 

    To accelerate the simulations and mitigate the curse of dimensionality, the methodology relies on the sparse-grid combination technique.
    This approach significantly reduces the number of degrees of freedom compared to full tensor grids while keeping comparable accuracy.
    Furthermore, the entire framework has been implemented within an MPI-based parallel code.
    By exploiting the natural decomposition into explicit and implicit components, the solver achieves competitive computational times in high-dimensional settings.

    The developed numerical schemes have been extensively validated \joseanto{through} several numerical experiments.
    First, systematic convergence studies were performed on two academic test problems. 
    Using a linear parabolic equation, we demonstrated the accuracy and efficiency of the methodology.
    Subsequently, a viscous Burgers' type model with an adjusted source term was employed \joseanto{to yield} analytical solutions.
    This example highlighted the scheme's ability to \joseanto{accurately} resolve nonlinear physical fluxes, while also identifying the issue of under-integration when combining nonlinearities with high-degree polynomial bases.

    Second, the methodology was validated in two application scenarios.
    In the Allen-Cahn problem subject to rotational transport, the sparse-grid yielded a drastic reduction in computational time, successfully capturing the same physical phenomena as the full-grid reference at a fraction of the cost.
    Finally, applying the framework to a high-dimensional problem in quantitative finance demonstrated its viability as a robust alternative for the overnight assessment of option contracts.
    
    These experiments confirm the predicted high-order convergence rates and allow for a detailed comparison of accuracy and efficiency trade-offs between full and sparse grids. 
    Overall, the results \joseanto{establish} that sparse-grid IMEX-LDG schemes can deliver comparable precision at a substantially lower computational cost, demonstrating that the methodology provides a powerful and flexible tool for solving high-dimensional nonlinear parabolic PDEs arising in applications like phase-field models or quantitative finance.
        
\section*{Acknowledgements}
This work is supported by the project PCI2024-155061-2 and by the European Union - NextGenerationEU program and by grant PID2022-137637NB-C21 funded by MCIN/AEI/\allowbreak 10.13039/\allowbreak 501100011033 and ``ERDF A way of making Europe'', and also by the Grant PPRO-FQM216-G-2023 funded by Consejer{\'i}a de Universidad, Investigaci{\'o}n e Innovaci{\'o}n and by ERDF Andalusia Program 2021-2027. Manuel J. Castro thanks IMUS-Maria de Maeztu grant CEX2024-001517-M - Apoyo a Unidades de Excelencia Mar{\'i}a de Maeztu for supporting this research, funded by MICIU/\allowbreak AEI/\allowbreak 10.13039/\allowbreak 501100011033.

\bibliographystyle{plain}
\bibliography{refs}

\end{document}